\input amstex
\documentstyle{amsppt}

\input epsf 		
\pageheight{9truein}
\pagewidth{6.3truein}
\NoBlackBoxes
\TagsOnLeft
\def\complex{{\Bbb C}}
\def\HH{{\Bbb H}}
\def\que{{\Bbb Q}}
\def\real{{\Bbb R}}
\def\zed{{\Bbb Z}}
\def\S{{\Cal S}}
\def\cpt{{\complex P^2}}
\def\id{\operatorname{id}}
\def\Im{\operatorname{Im}}
\def\mod{\operatorname{mod}}
\def\rel{\operatorname{rel}}
\def\rint{\operatorname{int}}
\def\chix{{\raise.5ex\hbox{$\chi$}}}

\def\varep{\varepsilon}
\def\hatH{\widehat H}
\def\hatX{\widehat X}
\def\tD{\widetilde D}
\def\tF{\widetilde F}
\def\tM{\widetilde M}
\def\tX{\widetilde X}
\def\oM{\overline M}

\def\Spin{\text{\rm Spin}} 
\def\spinc{{\text{\rm spin}^c}} 
\def\Spinc{{\text{\rm Spin}^c}} 
\def\STop{{\text{STop}}}
\def\SpinTop{{\text{SpinTop}}}
\def\llbrack{\lbrack\!\lbrack}
\def\rrbrack{\rbrack\!\rbrack}
\def\cdubfig#1#2#3#4#5#6{$$\matrix 
	\epsfysize=#1truein \epsfbox{#2} &\qquad 
	& \epsfysize=#4truein \epsfbox{#5}\cr
	\noalign{\vskip6pt}
	\text{\smc Figure #3}&\qquad &\text{\smc Figure #6} 
	\endmatrix$$}
\def\cfig#1#2#3{
	\midinsert {\epsfysize=#1 truein  
	\centerline{\epsfbox{#2}} }
	\botcaption{Figure #3} 
	\endcaption
	\endinsert} 
\def\cxfig#1#2#3{
	\midinsert {\epsfxsize=#1 truein  
	\centerline{\epsfbox{#2}} }
	\botcaption{Figure #3} 
	\endcaption
	\endinsert} 

\topmatter
\title Handlebody Construction of Stein Surfaces\endtitle
\author Robert E. Gompf*\endauthor 
\leftheadtext{ROBERT E. GOMPF}
\date December 3, 1996\enddate
\address Department of Mathematics, The University of Texas at Austin, 
Austin, TX 78712\endaddress
\email gompf\@math.utexas.edu\endemail 
\thanks * Partially supported by NSF grant DMS 9301524.\endthanks 
\abstract 
The topology of Stein surfaces and contact 3-manifolds is studied by 
means of handle decompositions. 
A simple characterization of homeomorphism types of Stein surfaces is 
obtained --- they correspond to open handlebodies with all handles of 
index $\le2$. An uncountable collection of exotic $\real^4$'s is shown 
to admit Stein structures. 
New invariants of contact 3-manifolds are produced, including a complete 
(and computable) set of invariants for determining the homotopy class 
of a 2-plane field on a 3-manifold. 
These invariants are applicable to Seiberg-Witten theory. 
Several families of oriented 3-manifolds are examined, namely the Seifert 
fibered spaces and all surgeries on various links in $S^3$, and in each 
case it is seen that ``most'' members of the family are the oriented 
boundaries of Stein surfaces. 
\endabstract 
\endtopmatter

\document
\subhead 0. Introduction\endsubhead 

It is becoming evident that low-dimensional topology is intimately connected 
with the topology of complex, symplectic and contact manifolds. 
The differential topology of closed 4-manifolds is entwined with that 
of complex surfaces and symplectic 4-manifolds, as is clear from Donaldson 
theory (for example \cite{Do}) and more recent developments involving the 
Seiberg-Witten equations and Gromov invariants (for example, Taubes \cite{T}). 
Cutting and pasting leads us to consider manifolds with boundary, whose 
corresponding geometric analogs are compact complex 
and symplectic manifolds with boundaries that inherit compatible tight 
contact structures. 
The topology of 3-manifolds with tight contact structures is subtle and still 
mysterious, as is the question of when they bound complex or 
symplectic 4-manifolds. 
There are some similarities, however, with the study of taut foliations 
on 3-manifolds, which has had a profound influence on 3-manifold topology.  
Considering the interiors of complex or symplectic manifolds with contact 
boundaries, we are led to the notion of a Stein manifold, which is perhaps 
the most natural notion of an open complex or exact symplectic manifold 
that is ``nice at infinity.'' 
It is natural to ask which open 4-manifolds admit Stein structures. 
Eliashberg \cite{E2} proved that Stein manifolds can be characterized 
in terms of handle decompositions. 
In the present article, we will use Eliashberg's theorem to study Stein 
surfaces (real dimension~4), producing a simple characterization of their 
homeomorphism types and some new examples such as Stein structures on 
exotic $\real^4$'s. 
We will also use handlebodies to study contact 3-manifolds and their 
holomorphic fillings, i.e., compact complex surfaces with contact 
boundaries. 

We begin by discussing Stein manifolds. (See \cite{E2}, \cite{E6} 
for more details and 
Section~1 for some key definitions in low dimensions.) 
{\it Stein manifolds\/} are complex manifolds 
(necessarily noncompact)  that admit proper holomorphic embeddings 
in $\complex^N$ for sufficiently large $N$. 
They carry exact K\"ahler structures, and an equivalent notion can be 
formulated in terms of symplectic structures \cite{EG}. 
A complex manifold $X$ is Stein if and only if it admits an 
``exhausting strictly plurisubharmonic function,'' 
which is essentially characterized as being a proper function 
$f:X\to\real$ that is bounded below and can be assumed a Morse function, 
whose level sets $f^{-1}(c)$ are ``strictly 
pseudoconvex'' (away from critical points), where $f^{-1}(c)$ is oriented 
as the boundary of the complex manifold $f^{-1}(-\infty,c]$. 
Strict pseudoconvexity implies, and for $X$ of real dimension~4 is equivalent 
to asserting, that $f^{-1}(c)$ inherits a contact structure determining 
its given orientation. 

We will call a compact, complex $X$ with boundary a {\it Stein domain\/} 
if it admits a strictly plurisubharmonic function 
such that the boundary $\partial X$ is a level set. 
Then $\partial X$ will be strictly pseudoconvex and $\rint X$ will be 
a Stein manifold (with the same level set structure as $X$). 
Conversely, if $f:X\to\real$ is an exhausting 
strictly plurisubharmonic function on 
a Stein manifold, then $f^{-1}(-\infty,c]$ will be a Stein domain 
for any regular value~$c$. 
Thus, we can almost interchangeably talk about Stein domains  
and Stein manifolds in the presence of (exhausting) 
strictly plurisubharmonic Morse 
functions  with only finitely many critical points. 
We will call a Stein manifold or domain of complex dimension~2 a 
{\it Stein surface\/} ({\it with boundary\/}). 
Any compact, complex surface with  nonempty, strictly pseudoconvex 
boundary can be made Stein in this sense by deforming it and blowing down 
any exceptional curves \cite{Bo}. 

Eliashberg's theorem \cite{E2} characterizing Stein manifolds of 
complex dimension ${n >2}$ can be stated as follows. (For an independent
treatment of handle addition in the symplectic category, see \cite{W}.)

\proclaim{Theorem 0.1 {\rm (Eliashberg)}} 
For $n>2$, a smooth, almost-complex, open $2n$-manifold admits a Stein 
structure if and only if it is the interior of a (possibly infinite) 
handlebody without handles of index $>n$. 
The Stein structure can be chosen to be homotopic to the given almost-complex 
structure, and the given handle decomposition is induced by an exhausting 
strictly plurisubharmonic Morse function. 
Similarly, a smooth, almost complex, compact $2n$-manifold $(n>2)$ with a 
handle decomposition without handles of index $>n$ admits a homotopic 
Stein domain structure with a suitable plurisubharmonic function 
inducing the handle decomposition. 
\endproclaim 

Thus, in high dimensions, an almost-complex structure and a handle 
decomposition with no  handles of index above the middle dimension are 
sufficient to guarantee the existence of a Stein structure.  
Implicit in the same paper is a theorem in the ${n=2}$ case (which is also 
due to Eliashberg, cf.\ \cite{E5}, but not explicitly 
published). 
We state it as Theorem~1.3 below. 
When ${n=2}$, the required almost-complex structure always exists (for 
oriented handlebodies as in the theorem). 
The main difference between this case and the ${n>2}$ case, however, is that 
when ${n=2}$, there is a serious restriction on the allowable framings 
by which the 2-handles are attached. 
In practice, it is a delicate matter to determine whether a given 
4-manifold admits a handle structure with allowable framings, and if so, which 
almost-complex structures can be so realized. 
Thus, while Eliashberg has characterized Stein manifolds in all dimensions 
in terms of differential topology, the case of 4-manifolds and Stein 
surfaces still presents great challenges in applications. 
It is these applications that the present article addresses. 

One application is the characterization of Stein surfaces (without boundary) 
up to homeomorphism (Section~3). 
Freedman's work on topological 4-manifolds \cite{F}, \cite{FQ} shows that 
if we forget about smooth structures, then the theory of 4-manifolds 
(with small $\pi_1$) 
becomes similar to that of higher dimensional manifolds, hence, 
easier to understand. 
In that spirit, we prove (Theorem~3.1) that Eliashberg's theorem in 
high dimensions applies up to homeomorphism to 4-manifolds. 
That is, an open, oriented topological 4-manifold $X$ is 
(orientation-preserving) homeomorphic to a Stein surface if and only if 
it is the interior of a topological (or smooth) handlebody without 
handles of index $>2$, and if so, then any almost-complex structure can be 
so realized. 
The given handle structure will not necessarily come from a 
plurisubharmonic (or even smooth Morse) function on the Stein surface, however. 
As is typical of examples obtained from Freedman theory, these Stein surfaces 
will tend to have smooth structures that are in some sense ``exotic.'' 
We cannot expect them to admit proper Morse functions with finitely many 
critical points. 
As an example (Theorem~3.3), $\cpt$ minus a point admits an uncountable 
family of diffeomorphism types of Stein exotic smooth structures, none of 
which admit proper Morse functions with finitely many critical points, and none
of which contain a smoothly embedded sphere representing a generator 
of the homology. 
Similarly, any $\real^2$-bundle over $S^2$ admits a Stein exotic smooth 
structure containing no generating smoothly embedded spheres. 
In contrast, the standard smooth structure on $\cpt -\rint B^4$ 
(or on $S^2\times D^2$) cannot be realized by a compact Stein surface (or even 
a convex symplectic manifold), by uniqueness of fillings \cite{Gro}, \cite{E3}. 
As a further application, we show that $\real^4$ admits uncountably many 
exotic smooth structures that can be realized as Stein surfaces 
(Theorem~3.4). 
None of these admit proper Morse functions with finitely 
many critical points, but 
subject to that constraint, one example has a handle decomposition that is 
remarkably simple. 
(It is the example of a simple exotic $\real^4$ constructed by Bi\v zaca and 
the author in \cite{BG}.) 

One of the main techniques of this paper is to describe handle 
decompositions of Stein surfaces explicitly using Kirby calculus. 
While this method has already been applied in simple cases without 
1-handles \cite{E5}, the general case is more delicate. 
In Section~2, we establish a standard form for any handle decomposition 
obtained from a strictly plurisubharmonic function on  a 
compact Stein surface. 
We do this via a standard form for Legendrian links in 
the connected sum $\# nS^1\times S^2$ 
that allows us to define and compute the rotation number and Thurston-Bennequin 
invariant of each link component (even those that are nontrivial in $H_1$). 
We also provide a complete reduction from Legendrian link theory in 
$\# nS^1\times S^2$ to a theory of diagrams by introducing a complete set 
of ``Reidemeister moves.'' 
These diagrams allow us to construct Stein surfaces by drawing pictures. 
For example, we obtain the above exotic $\real^4$'s in this manner. 
A Legendrian link diagram in $\# nS^1\times S^2$ also determines a (positively  
oriented) contact 3-manifold $(M,\xi)$, namely the oriented boundary of the 
corresponding  compact Stein surface. 
We say that $(M,\xi)$ is obtained by {\it contact surgery\/} on the 
Legendrian link, and that $(M,\xi)$ is {\it holomorphically fillable\/}. 
In Section~5, we construct several families of examples. 
We realize ``most'' oriented Seifert fibered 
3-manifolds by contact surgery, including all with (possibly nonorientable) 
base $\ne S^2$, and both orientations on many Brieskorn homology spheres 
(Theorem~5.4 and Corollary~5.5). 
We show that any Seifert fibered space can be  realized in this manner 
after possibly reversing orientation. 
For hyperbolic examples, we realize ``most'' rational surgeries on the 
Borromean rings  (Theorem~5.9). 

We also introduce new invariants for distinguishing contact structures. 
We define a complete set of invariants for determining the homotopy 
class of an oriented 2-plane field on an oriented 3-manifold $M$. 
These invariants are readily computable for the boundary of a compact 
Stein surface presented in standard form. 
An explicit formula for the 2-dimensional obstruction (which measures the 
associated spin$^c$-structure) is given by Theorem~4.12. 
The 3-dimensional obstruction (which, by recent work of Kronheimer 
and Mrowka \cite{KM}, distinguishes the grading in Seiberg-Witten-Floer Theory) 
is given by Definitions~4.2 and 4.15. 
The construction of these invariants is surprisingly delicate. 
While a choice of trivialization on the tangent bundle $TM$ reduces the 
problem to the homotopy classification of maps $M\to S^2$ (which was 
solved by Pontrjagin around 1940 \cite{P}), the resulting obstructions depend 
on the choice of trivialization, making them hard to work with directly. 
Our invariants depend (in the worst case) only on a spin structure and 
a framing on a certain 1-cycle in $M$, data that one can easily follow 
through Kirby calculus computations.
As an application, we show (Corollary~4.6) that the rotation number (up to 
sign) and Thurston-Bennequin invariant of a Legendrian knot in $S^3$ are 
both invariants of the contact 3-manifold $M$ obtained by 
contact surgery on the knot. 
In particular, the sign of the surgery coefficient is determined by the 
contact structure on $M$. 
We then observe that we can easily construct 
families of different (nonhomotopic and noncontactomorphic) holomorphically 
fillable contact structures,  on a fixed 3-manifold, that cannot be 
distinguished by the Chern (= Euler) class of $\xi$. 
This provides new counterexamples to Conjecture~10.3 of \cite{E3}. 
(Such examples were already known in the weaker case of symplectically 
fillable structures \cite{Gi2}, although these were all homotopic.) 
We give several corollaries about homotopy classes of 2-plane fields 
(or equivalently,  nowhere zero vector fields or ``combings'') on 3-manifolds. 
We also show (Corollary~4.19) that any contact structure  respecting 
the unusual orientation on the 
Poincar\'e homology sphere must have an overtwisted universal cover. 
Although Conjecture~5.6 
asserts that this oriented manifold should admit no fillable structures, 
we show in Proposition~5.1 that it is common 
(among lens spaces, for example) for holomorphically fillable contact 
structures to have finite covers that are overtwisted. 
In fact, we exploit this phenomenon in Example~5.2 to obtain 
a more subtle way of distinguishing tight contact structures. 
We exhibit a 3-manifold with a pair of  
holomorphically fillable contact structures that are homotopic as plane 
fields, and distinguish these by whether 
the corresponding contact structures on a certain 2-fold covering space 
are tight or overtwisted. 
Recent advances in gauge theory such as \cite{KM} are leading to other methods 
for distinguishing homotopic contact structures; see \cite{AM} and \cite{LM}. 

The author wishes to thank Yasha Eliashberg for many indispensable 
conversations, and to acknowledge the Isaac Newton Institute for 
Mathematical Sciences in Cambridge, England 
for their support during their 1994 program 
on symplectic and contact topology, at whose lectures by Eliashberg the 
author was first properly introduced to contact topology and Stein surfaces. 

\subhead 1. Legendrian links\endsubhead 

In this section, we review some standard theory of contact 3-manifolds, 
Legendrian links and their relation with Stein surfaces. 
Throughout the paper, we will be working with $C^\infty$ oriented 2-plane 
fields $\xi$ on oriented (usually closed) 3-manifolds $M$. 
Such a 2-plane field can be written as the kernel of a nowhere zero 1-form 
$\alpha$ on $M$ that is unique up to multiplication by nonzero scalar 
functions. 
It is easily verified that the integrability of $\xi$ is equivalent to 
the condition that $\alpha\wedge d\alpha$ be identically zero. 
We call $(M,\xi)$ a {\it contact manifold\/} if $\xi$ is completely 
nonintegrable in the sense that $\alpha\wedge d\alpha$ is nowhere zero. 
(This is clearly independent of the choice of $\alpha$.) 
Then $\alpha\wedge d\alpha$ determines an orientation on $M$ 
(independent of $\alpha$ and the orientation of $\xi$). 
The contact structure $\xi$ is called {\it positive\/} if this orientation 
agrees with the given one on $M$ and {\it negative\/} otherwise. 
Except where otherwise indicated, we will (without loss of generality) 
only deal with positive contact structures.  
A {\it contactomorphism\/} is a diffeomorphism preserving contact 
structures. 
According to Gray's Theorem \cite{Gr}, contact structures are all 
locally contactomorphic, and on closed 
manifolds they are deformation invariant in the sense that if a homotopy 
$\xi_t$, $0\le t\le1$, of 2-plane fields on $M$ consists entirely of contact 
structures, then there is an isotopy $\varphi_t$ of $M$ with $\varphi_0 
= \id_M$, $(\varphi_t)_* \xi_0 = \xi_t$ for each $t$, and $\varphi_t = 
\varphi_0$ wherever $\xi_t$ is independent of $t$. 
When this occurs, we say that $\xi_0$ and $\xi_1$ are {\it isotopic\/}. 
Note that contact structures that are homotopic as plane fields 
need not be homotopic through contact structures. 
It follows easily that an isotopy $\psi_t$ of $M$ that preserves $\xi$ 
for all $t$ on a subset $N$ can be changed $\rel N$ 
to a {\it contact isotopy\/}, or 
isotopy of $M$ through contactomorphisms. 
(Apply Gray's Theorem to the family $\psi_t^* \xi$.)

A link $L:\coprod_{i=1}^n S^1 \hookrightarrow M$ in a contact 3-manifold 
$(M,\xi)$ is called {\it Legendrian\/} \cite{Ar} if its tangent vectors all lie 
in $\xi$.  
Two such links are considered equivalent if they are isotopic through a 
family of Legendrian links, or equivalently, if they are contact isotopic.  
Any link in a contact 3-manifold is $C^0$-small isotopic to a (nonunique) 
Legendrian link. 
(Simply replace each arc transverse to $\xi$ by a (left-handed) 
Legendrian spiral.) 
Any diffeomorphism $\varphi:K\to K'$
between Legendrian knots extends to a contactomorphism 
on some neighborhoods of the knots. 
A Legendrian link comes equipped with a {\it canonical framing\/} of its 
normal bundle (up to fiber homotopy and orientation reversal), 
which is induced by any vector field 
transverse to $\xi$, or equivalently, by a vector field in $\xi|L$ 
transverse to $L$. 
This framing is preserved by contactomorphisms. 
For any nullhomologous knot $K$, there is a 
canonical bijection from (normal) framings of $K$ to the integers, sending 
each framing $f$ to the linking number of $K$ with its push-off 
determined by $f$, and in the Legendrian case 
the integer corresponding to the canonical 
framing  of $K$ is called the {\it Thurston-Bennequin invariant\/} $tb(K)$. 
For any Legendrian knot $K$, we can find a $C^0$-small isotopy (necessarily 
changing its Legendrian knot type) that adds any number of left (negative) 
twists to the canonical framing --- or in the nullhomologous case, 
decreases $tb(K)$ by any integer. 
(Simply add a spiral to $K$.) 
It is not always possible to add right twists (increase $tb(K)$), however. 
If $(M,\xi)$ admits a topologically unknotted Legendrian knot $K$ with 
$tb(K)=0$, then $\xi$ is {\it overtwisted\/}. 
(The reader can take this as a definition, or note that any disk bounded by $K$ 
is isotopic to an overtwisted disk.) In this case,
we can add right twists to any 
canonical framing. 
Otherwise, $\xi$ is called {\it tight\/}, and any nullhomologous knot will 
have a Legendrian representative with maximal $tb$. 
For the unknot, the maximum will be $-1$. 
There are knots in $S^3$, with its standard tight structure, for which the 
maximal $tb$ is arbitrarily large or small. (See \cite{R1};  
for $q=TB$, see \cite{R2}). 
For any nullhomologous knot in a tight contact manifold there is a bound 
of $tb(K) \le -\chix (F)$ for any embedded, orientable, 
connected surface $F$ bounded 
by $K$ \cite{E5}. 
(See below for a sharper statement.) 

The most interesting contact structures on 3-manifolds are the tight ones 
--- they are somewhat analogous to taut foliations. 
While the classification of overtwisted contact structures on a closed 
3-manifold is simple (there is a unique such structure in each homotopy 
class of 2-plane fields \cite{E1}), the occurrence of tight structures is 
poorly understood. 
Almost the only known obstruction to the existence of such structures $\xi$ 
is  a bound on $c_1(\xi) \in H^2 (M;\zed)$, the Chern class (or equivalently, 
Euler class) of $\xi$ as a complex line (real oriented 2-plane) bundle. 
Namely, if $F$ is a closed, connected, oriented surface in $M$, then 
$|\langle c_1(\xi),F\rangle| \le -\chix (F)$ for $F\ne S^2$, and the 
left-hand side vanishes for $F=S^2$ \cite{E4}. 
A primary source of tight contact structures is as follows: 
Consider a compact, complex surface $X$ with boundary $M$. 
Each tangent space $T_xM$ to $M$ will contain a unique 1-dimensional 
complex subspace of $T_xX$ (namely $T_xM\cap iT_xM$). 
If these complex lines comprise a contact structure $\xi$ on $M$ 
determining the boundary orientation, then $X$ has a {\it strictly 
pseudoconvex\/} boundary, and $(M,\xi)$ is called {\it holomorphically 
fillable\/}. 
Eliashberg \cite{E3} proved that any holomorphically fillable contact 
3-manifold is tight. 
This theorem and its symplectic generalization \cite{E3} are the main tools 
available for proving tightness of contact structures. 
Note that if $M=\partial X$ and $\xi$ is the complex line field on $M$ 
induced by a complex (or almost-complex) structure $J$ on $X$, then 
the Chern class $c_1(\xi) \in H^2(M;\zed)$ is the restriction of that 
of $J$ on $X$, $c_1(J)\in H^2(X;\zed)$. 
This is because the complex bundle $TX|M$ splits as the sum of $\xi$ and 
a trivial complex line bundle. 

There is one additional invariant known for Legendrian links. 
Let $L$ be a nullhomologous, oriented Legendrian link in $(M,\xi)$. 
(In practice, $L$ will be a knot.) 
Let $F$ be a {\it Seifert surface\/} for $L$, i.e., a compact, oriented surface 
embedded in $M$ (without closed components) whose oriented boundary is $L$. 
We define the {\it rotation number\/} $r(L,F)$ of $L$ with respect to $F$ 
to be the relative Chern number $\langle c_1(\xi,\tau),F\rangle$ of $\xi$ 
relative to a tangent vector field $\tau$ along $L$, evaluated on $F$. 
(Note that $c_1(\xi,\tau)\in H^2 (M,L;\zed)$.) 
That is, we compute $r(L,F)$ by trivializing the 2-plane bundle $\xi|F$ 
and counting (with sign) how many times $\tau$ rotates in $\xi$ 
with respect to the trivialization as we travel around $L$. 
Clearly, $r(L,F)$ only depends on $F$ through its homology class  in 
$H_2(M,L;\zed)$, and if $c_1(\xi)=0\in H^2(M;\zed)$ then $r(L) = r(L,F)$ 
is independent of $F$. 
Note that $r(L,F)$ reverses sign if we reverse the orientation of either 
$L$ (hence, $F$) or $\xi$. 
By adding left twists to $L$, one can realize any preassigned value of 
$r(L,F)$ at the expense of decreasing $tb$, as will be clear from 
our pictures below. 
In a tight contact manifold, the invariants 
$r$ and $tb$ classify those Legendrian knots that are 
topologically unknotted \cite{EF}, and  
for arbitrary nullhomologous Legendrian knots $K$ our previous bound on 
$tb(K)$ is sharpened by the inequality $tb(K)+|r(K)| \le -\chix(F)$ \cite{E5}. 

Our basic examples of contact 3-manifolds will be $\real^3$, $S^3$ and 
the connected sum 
$\# nS^1\times S^2$, each of which admits a unique (up to isotopy) tight 
contact structure compatible with its standard orientation 
\cite{Be}, \cite{E4}. 
The structures on $S^3$ and $\#nS^1\times S^2$ are uniquely (up to blowups) 
holomorphically fillable --- $S^3$ as the boundary of a round ball $B^4$ in 
$\complex^2$, and $\#nS^1\times S^2$ as the boundary of $B^4$ union $n$ 
1-handles \cite{E3}.  
To see the tight contact structure $\xi$ on $S^3$, we delete a point to 
obtain the tight structure on $\real^3$. 
We will always represent this 
as the kernel of the 1-form $\alpha = dz+xdy$ on $\real^3$. 
(Note that $\alpha \wedge d\alpha = dx\wedge dy\wedge dz$.) 
We orient $\xi$ via the nowhere zero form $d\alpha|\xi= dx\wedge dy |\xi$. 
(Note that the contactomorphism $(x,y,z)\mapsto (x,-y,-z)$ reverses 
orientation on $\xi$.) 
We visualize $(\real^3,\xi)$ by projecting into the $y$-$z$ plane. 
Then the plane $\xi_{(x,y,z)}$ at a point $(x,y,z)$ projects to a line  
at $(y,z)$ whose slope is $-x$. 

Legendrian link theory in $S^3$ or $\real^3$ now reduces without loss of 
information to the theory of the corresponding {\it front projections\/} 
in $\real^2$, as developed by Arnol'd \cite{Ar}. 
A Legendrian knot in $\real^3$ projects to a closed curve $\gamma$ in $\real^2$ 
that may have cusps  and transverse self-crossings but has 
no vertical tangencies. 
(See Figure~1.) 
Any such curve comes from a unique Legendrian knot in $\real^3$, which may be 
reconstructed by setting $-x(t)$ equal to the slope of $\gamma$ at $t$, 
with cusps corresponding to points where the knot  is parallel 
to the $x$-axis. 
Thus, at self-crossings, the curve of most negative slope always crosses 
in front. 
For example, Figure~1 represents the right-handed trefoil knot as in Figure~2. 
We will continue to draw overcrossings to avoid confusion, even though it is 
not actually necessary. 
Beware that the literature contains diagrams using the opposite convention, 
i.e., a left-handed coordinate system. 
In a front projection of a generic Legendrian link, the only singularities 
are transverse double points and cusps isotopic to the curves 
$z^2 = y^3$ or $-y^3$. 
In analogy with the Reidemeister moves of ordinary link theory, the moves 
shown in Figure~3 (together with their images under $180^\circ$ rotation 
about any coordinate axis in $\real^3$ and isotopies of $\real^2$ 
introducing no vertical tangencies) suffice for realizing any 
equivalence of Legendrian links  \cite{S}. 
\cdubfig{1.75}{rg-fig01.eps}{1}{1.25}{rg-fig02.eps}{2}
\cfig{3.25}{rg-fig03.eps}{3}

The invariants $tb(K)$ and $r(K)$ for an oriented Legendrian knot $K$ in 
$S^3$ or $\real^3$ are easy to compute from a front projection. 
(Note that these are always well-defined, since 
any knot is nullhomologous in $S^3$ or $\real^3$, and $c_1(\xi)=0$.) 
We begin with $tb(K)$. 
Any smooth planar diagram of a knot $K$ in $S^3$ determines an obvious 
{\it blackboard framing\/} via a normal vector field to the immersed curve 
in $\real^2$. 
This framing is not isotopy invariant. 
In fact, the integer corresponding to the blackboard framing is equal to the 
{\it writhe\/} $w(K)$ of $K$, which is the number of self-crossings of $K$ 
counted with sign. 
(Signs are determined by orienting $K$ arbitrarily and comparing with 
Figure~4, up to rotation.) 
\smallskip

{\vbox{\epsfysize=1.25truein \centerline{\hskip1.25truein\epsfbox{rg-fig04.eps}}
	\centerline{\smc Figure 4} }}

\noindent
If $K$ is Legendrian, we can compute its canonical framing by observing that  
the vector field ${\partial\over\partial z}$ on $\real^3$ is everywhere 
transverse to $\xi$, so it gives the canonical framing of $K$. 
This framing will agree with the blackboard framing (obtained by smoothing 
the projection of $K$ without adding crossings) 
except for a half left twist at each cusp (Figure~5). 
\cfig{3.0}{rg-fig05.eps}{5}
\noindent
Thus, if $\lambda (K)$ (resp. $\rho (K)$) denotes the number of cusps with 
vertex on the left (resp. right) (Figure~5), we obtain that the canonical 
framing differs from the blackboard framing by $\frac12 (\lambda(K) + 
\rho (K))$ left twists. 
But clearly, $\lambda(K) = \rho(K)$, so we obtain 
$$tb(K) = w(K) - \frac12 (\lambda (K) + \rho (K)) = w (K) - \lambda (K)\ .
\tag 1.1$$
For example, the trefoil in Figure 2 has $tb=1$. 

To compute the rotation number $r(K)$, we observe that the vector field 
$\partial\over\partial x$ trivializes the 2-plane bundle $\xi$ on $\real^3$, 
so it restricts to a trivialization of $\xi|F$ for any Seifert surface $F$. 
Thus, it suffices to count (with sign) how many times the tangent vector 
field $\tau$ of $K$ crosses $\partial\over\partial x$ as we travel around $K$. 
Let $\lambda_+(K)$ (resp. $\lambda_-(K)$) be the number of left cusps at 
which $K$ is oriented upward (resp. downward), and define $\rho_\pm (K)$ 
similarly (Figure~6). 
Let $t_\pm = \lambda_\pm + \rho_\pm$ be the total number of upward $(t_+)$ 
and downward $(t_-)$ cusps. 
Since each downward left cusp represents a positive crossing of $\tau$ 
past $\partial\over\partial x$ (counterclockwise with respect to 
$dx\wedge dy$) and each upward right cusp represents 
a negative crossing, we obtain 
$$r(K) = \lambda_- - \rho_+ = \rho_- - \lambda_+ = \frac12 (t_- - t_+)\ ,
\tag 1.2$$ 
where the second equality is obtained by using $-{\partial\over\partial x}$ 
instead of $\partial\over\partial x$ and the third is obtained by averaging 
the first two. 
Note that reversing the orientation of $K$ reverses the sign of 
$r(K)$, as required. 
If we add $k>0$ upward (downward) zig-zags to $K$ as in Figure~7 (with 
$k=3$), the effect will be to decrease $tb(K)$ by $k$ and decrease (increase) 
$r(K)$ by $k$. 

\cfig{2.0}{rg-fig06.eps}{6}
\cfig{2.0}{rg-fig07.eps}{7}

The main purpose of this discussion was to prepare for the statement and 
applications of Eliashberg's Theorem in dimension~4. 
Recall \cite{K}, \cite{GS}
that a handle decomposition of a compact, oriented 4-manifold $X$ 
with all handles of index $\le2$ is equivalent to a  framed link 
in $\#nS^1\times S^2 = \partial$(0-handle $\cup$ 1-handles). 
For each link component, we attach a 2-handle $D^2\times D^2$ along 
$S^1\times D^2$, by identifying the link component with $S^1\times 0$ and 
its framing with the product framing on $S^1\times 0\subset S^1\times D^2$. 
If the attaching circle is nullhomologous in $\#nS^1\times S^2$, then its 
framing is specified by an integer, which is also the self-intersection number 
$\alpha^2$ ($=\alpha\cdot\alpha$) of the homology class 
$\alpha \in H_2(X;\zed)$ determined 
(up to sign) by the 2-handle. 
For noncompact 4-manifolds we obtain a similar description, although it can be 
more complicated due to the end structure  of the union of 0- and 1-handles. 
For Stein surfaces, we have the following theorem, which is implicit 
in \cite{E2}. 
(See also \cite{E5}.) 

\proclaim{Theorem 1.3 {\rm (Eliashberg)}} 
A smooth, oriented, open 4-manifold $X$ admits a Stein structure if and 
only if it is the interior of a (possibly infinite) handlebody such that 
the following hold: 
\roster 
\item"(a)" Each handle has index $\le2$, 
\item"(b)" Each 2-handle $h_i$ is attached along a Legendrian curve $K_i$ in 
the contact structure induced on the boundary of the underlying $0$- and 
$1$-handles, and 
\item"(c)" The framing for attaching each $h_i$ is obtained from the canonical 
framing on $K_i$ by adding a single left (negative) twist. 
\endroster
A smooth, oriented, compact 4-manifold $X$ admits a Stein structure if and 
only if it has a handle decomposition satisfying {\rm a}, {\rm b} and {\rm c}. 
In either case, any such handle decomposition comes from a strictly 
plurisubharmonic function (with $\partial X$ a level set). 
\endproclaim 

Thus, the previous discussion gives an explicit procedure for constructing 
Stein surfaces with strictly plurisubharmonic Morse functions lacking index~1 
critical points --- 
Simply add 2-handles to a Legendrian link in $S^3 = \partial B^4$ such 
that the framing on each $h_i$ is given by $tb(K_i)-1$. 
We deal with general Stein surfaces in the next section. 

Eliashberg proves this theorem by explicit holomorphic gluing. 
Each 2-handle is given as a neighborhood of $D^2\times 0\subset i\real^2
\times \real^2=\complex^2$. 
(Note that the equality reverses the natural orientations.) 
The attaching circle $S^1\times 0$ is glued to the given Legendrian curve. 
In particular, the unit tangent vector field $\tau$ to $S^1\times 0$ 
is mapped into $\xi$. 
Since $\xi$ now represents a complex line field, $i\tau$ must also map 
into $\xi$, so $i\tau$ goes to the canonical framing. 
However, $i\tau$ differs from the product framing on $S^1\times \real^2$ 
by one twist, since $\tau :S^1\to S^1$ has degree~1. 
This accounts for the difference of one twist between the canonical framing 
and the attaching framing of the 2-handle. 
(The twist is left-handed because of the above-mentioned orientation 
reversal.) 

The Chern class of a Stein surface $X$ without 1-handles is easy to 
compute from a Legendrian link diagram (cf.\ \cite{E5}). 
The homology $H_2(X;\zed)$ will be free abelian. 
To choose a basis, one orients the framed link. 
Then each link component is the oriented boundary of a Seifert surface 
whose union with the (suitably oriented) core of the 2-handle 
represents a basis element. 
The Chern class $c_1(X)$ can now be characterized as the unique element 
of $H^2(X;\zed)$ whose value on each such basis element is the rotation  
number of the corresponding oriented Legendrian link component. 
We will prove a more general statement in the next section 
(Proposition~2.3).

\subhead 2. 1-handles\endsubhead 

In order to represent arbitrary Stein surfaces by means of Legendrian link 
diagrams, we must find a way to represent 1-handles. 
In the usual Kirby calculus for representing handlebodies via framed links 
\cite{K}, \cite{GS},
one can represent a 1-handle by drawing its attaching region. 
This will be a pair of 3-balls $B_1$ and $B_2$ in $S^3 = \partial B^4$ 
(represented by $\real^3$), which are taken to be identified by a 
diffeomorphism $\varphi : B_1\to B_2$ 
that reverses orientation if the 4-manifold is orientable. 
Without loss of generality, we can take $B_1$ and $B_2$ to be round, and 
choose $\varphi$ conveniently. 
The standard convention is to take $\varphi$ to be reflection through 
the plane perpendicularly bisecting the line segment connecting the 
centers of $B_1$ and $B_2$. 
It is sometimes convenient to put the balls in a standard position and 
assume (without loss of generality) that the attaching curves of the 
2-handles lie in the region between the balls.  
For example, this is convenient if we wish to specify framings by integers,  
as will be discussed below. 

We will find a similar description of 1-handles in the setting of Stein 
surfaces, and establish a standard form for Legendrian links in 
$\#n S^1\times S^2$ that allows us to conveniently define and compute 
Thurston-Bennequin invariants. 
We must observe, however, that the diffeomorphism $\varphi$ is necessarily 
more complicated than in the smooth case. 
To understand this, we must consider the {\it characteristic foliation\/} 
of a surface $F$ in a contact 3-manifold $(M,\xi)$, which is the singular 
foliation on $F$ induced by the singular field of lines $\xi_x\cap TF_x$ 
on $F$.  
The characteristic foliation of a small round sphere in $\real^3$, centered in 
the $y$-$z$ plane, has two singular points, one at each pole. 
In between, the foliation is roughly a left-handed spiral with infinitely 
many  turns near each pole. 
The spirals are somewhat twisted, however, so that they always cross the 
$y$-$z$ plane orthogonally to it. 
If gluing $B_1$ and $B_2$ by $\varphi$ is to produce a manifold with 
contact boundary, then $\varphi$ must preserve the characteristic foliations 
on $\partial B_j$. 
If the contact plane field is to be orientable, then $\varphi$ must 
interchange the poles of the spheres $\partial B_j$ (so that an inward 
$\partial\over\partial z$ maps to an outward one). 
Since $\varphi$ reverses orientation,  the spheres $\partial B_j$ 
cannot both be round, for the left-handed spiral of $\partial B_j$ at each 
pole would necessarily map to a right-handed spiral. 
This is easily remedied by a $C^1$-small perturbation    near each pole 
to make $\partial B_j$ look locally like the saddle $z=-\frac12 xy$ at the 
origin, up to translation (cf.\ \cite{E4}). 
The characteristic foliations will then be radial near the poles, so it will 
be possible to define a map $\varphi$ preserving the characteristic foliations. 
Note that an additional twist will be required to match up the spirals 
away from the poles. 

Fortunately, we will see that in the end it is not necessary to understand 
$\varphi$ in detail. 
We will align $B_1$ and $B_2$ along the same horizontal line, show that we can 
assume $\varphi$ identifies the two points $p_j\in B_j$ that are closest 
to each other, and see that it suffices to keep track of $\varphi$ in a 
neighborhood of each $p_j$. 
Note that at each $p_j$, the characteristic foliation will  have tangent 
$\pm {\partial\over\partial x}$. 
If we orient these flow lines to point toward us 
$(+{\partial\over\partial x})$, they will continue toward opposite poles 
of $B_1$ and $B_2$. 
Thus, $\varphi$ will preserve the direction $\partial\over\partial x$ at $p_j$, 
which will allow us to approximate $\varphi$ by 
the usual reflection in a neighborhood of $p_1$. 
(Other conventions are possible --- for example, we could align $B_1$ and 
$B_2$ vertically and identify their nearest poles by translating the 
saddle-shaped regions. 
This seems to complicate the resulting Legendrian link diagrams, however, 
since one must deal carefully with Legendrian curves passing 
through the saddles.) 

\definition{Definition 2.1} 
A {\it Legendrian link diagram\/} 
in {\it standard form\/},  with $n\ge0$ 1-handles,  
is given by the following data (see Figure~8): 
\roster
\item"1)" A rectangular box parallel to the axes in $\real^2$,  
\item"2)" A collection of $n$ distinguished segments of 
each vertical side of the box, aligned horizontally in pairs 
and denoted by balls, and  
\item"3)" A front projection of a  generic 
Legendrian tangle (i.e., disjoint union of Legendrian knots and arcs) contained 
in the box, with endpoints lying in the distinguished segments and aligned 
horizontally in pairs.
\endroster
\enddefinition 
\cfig{3.0}{rg-fig08.eps}{8}

Thus, if we attach 1-handles to the pairs of balls, we will obtain a link 
in $\#n S^1\times S^2$. 
Now let $H$ denote a handlebody consisting of 
a 0-handle and $n$ 1-handles, with the canonical Stein structure 
determined by Theorem~1.3 and contact structure $\xi$ on $\partial H$. 
Fix an ordering of the 1-handles, a direction for each 1-handle and a 
homotopy class of nowhere zero vector fields in $\xi$. 
Such vector fields exist because $c_1(\xi)$ is the restriction of 
$c_1(H)\in H^2(H;\zed)=0$. 

\proclaim{Theorem 2.2} 
The boundary of $H$ can be identified with a contact manifold obtained 
from the standard contact structure on $S^3$ by removing smooth balls 
and gluing the resulting boundaries as in Figure~8. 
The identification exhibits $H$ in the usual way as a smooth handlebody, 
and it can be assumed to match the above data for $H$ with corresponding 
preassigned data for the diagram. 
Any Legendrian link in $\partial H = \#nS^1\times S^2$ is contact isotopic to 
one in standard form. 
Two Legendrian links in standard form are contact isotopic in 
$\partial H$ if and only if they are 
related by a sequence of the six moves shown in Figures~3 and 9 (and their 
images under $180^\circ$ rotation about each axis), together with isotopies 
of the box that fix the boundary outside of the balls and introduce no 
vertical tangencies. 
\endproclaim 
\cfig{2.25}{rg-fig09.eps}{9}

The new moves correspond to sliding cusps and crossings over 1-handles, and to 
swinging a strand around an attaching ball of a 1-handle. 
We prove the theorem at the end of this section. 
Following the proof, we show that different homotopy classes of vector fields 
result in different descriptions of $H$. 
These will be contactomorphic, and we describe the contactomorphisms 
explicitly, but no such contactomorphism will be contact isotopic to the 
identity. 
For a sample of the difficulties encountered in proving the theorem, 
consider how to put a link in standard form if 
it contains Figure~10  --- Beware of disallowed crossings. 
\cfig{1.5}{rg-fig10.eps}{10}

For arbitrary knots in $\#nS^1\times S^2$, we would like to establish a 
convention for identifying framings with integers, just as we did using 
linking numbers in the nullhomologous case. 
There is no isotopy-invariant way to do this. 
However, we can adopt a convention from Kirby calculus that assigns integers 
to framed knots that are in standard form. 
These integers will be invariant under isotopies within the box 
of Definition~2.1, i.e., Moves 1-5 (Figures~3, 9), but will change under 
Move~6. 
To establish the convention, simply stretch the box in the plane and glue 
its lateral edges together so that it becomes an annulus, 
identifying endpoints of the tangle in the natural way. 
The tangle becomes an actual knot projected into the plane. 
Now we can identify framings with integers using our previous convention. 
As before, the blackboard framing becomes $w(K)$, the signed number 
of self-crossings of $K$. 
Clearly, this definition has the required isotopy invariance. 
It also agrees with the previous one in the nullhomologous case. 
Move~6 changes the integer corresponding to any framing on $K$ by twice 
the number of times (counted with sign) that $K$ runs over the 1-handle. 
To check this and fix the sign, simply 
consider the blackboard framing, which is invariant under the move. 
Now for $K$ Legendrian, we define $tb(K)$ as before, as the integer 
corresponding to the canonical framing. 
It follows immediately that $tb(K)$ is still given by Formula~1.1. 
Note that $tb(K)$ can change under Move~6, even though the canonical 
framing itself is invariant under contactomorphisms. 

We must also define the rotation number $r(K)$ in more generality. 
Suppose that $K$ is an oriented Legendrian knot in a contact 3-manifold 
$(M,\xi)$ with $c_1(\xi)=0$. 
Thus, as an abstract real 2-plane bundle, $\xi$ is trivial. 
For any preassigned choice of nowhere zero vector field $v$ in $\xi$ 
(up to homotopies through such vector fields) we can define  
the rotation number $r_v(K)$ 
to be the signed number of times that the tangent vector 
field $\tau $ to $K$ rotates in $\xi$ relative to $v$ as we traverse $K$. 
In the nullhomologous case, $v$ determines a trivialization of $\xi|F$ 
on any Seifert surface $F$ for $K$, so our new definition agrees with 
the old one (and hence is independent of $v$) in this case. 
Now for the tight contact structure $\xi$ on $\#nS^1\times S^2$, 
we choose $v$ to be $\partial\over\partial x$ everywhere in the box in 
Definition~2.1. 
This fits together in the obvious way on the 1-handles, and then extends 
uniquely over $\#nS^1\times S^2$. 
(Such vector fields are classified by $H^1(\,\cdot\, ;\zed)$.) 
Now $r(K)= r_v(K)$ is well-defined for any oriented Legendrian knot $K$ in our 
diagram of $\#nS^1\times S^2$, and it is invariant under contact isotopy. 
Note, however, that Theorem~2.2 allows us to identify the diagram with $H$ 
so that $v$ corresponds to any preassigned nowhere zero vector field in $\xi$ 
in $\partial H$. 
Changing the latter vector field will change rotation numbers of Legendrian 
knots in $\partial H$. 
The resulting nontrivial contactomorphisms of $\partial H$ will be described 
explicitly at the end of this section. 
Our previous argument shows that $r(K)$ can still be computed via 
Formula~1.2 when $K$ is in standard form. 
In Section~4 we will show that $tb(K)+ r(K)+1$ is always congruent mod~2 to 
the number of times $K$ crosses 1-handles (Corollary~4.13). 

We are led to the following characterization of compact Stein surfaces 
with boundary. 

\proclaim{Proposition 2.3} 
A smooth, oriented, compact, connected 
4-manifold $X$ admits the structure of a Stein 
surface (with boundary) if and only if it is given by a handlebody on 
a Legendrian link in standard form (Definition~2.1) with the 
$i^{th}$ 2-handle $h_i$ attached to the $i^{th}$ link component $K_i$ 
with framing $tb(K_i)-1$ (as given by Formula~1.1). 
Any such handle decomposition is induced by a strictly plurisubharmonic 
function. The Chern class $c_1(J)\in H^2(X;\zed)$ of such a Stein structure 
$J$ is represented by a cocycle whose value on each $h_i$, oriented as 
at the end of Section~1, is $r(K_i)$ (as given by Formula~1.2). 
\endproclaim 

\demo{Proof} 
The first part of this proposition is just Eliashberg's theorem (1.3) 
augmented by Theorem~2.2 and the subsequent discussion. 
(To obtain a unique $0$-handle, apply the connectedness of $X$ and 
uniqueness of fillings of $S^3$.)
To compute $c_1(J)$, we trivialize $TX$ as a complex bundle over the 
union $X_1$ of the 0- and 1-handles, and measure the failure of the 
trivialization to extend over the 2-handles. 
We trivialize $TX$ over the box given in Definition~2.1 (as a subset 
of $X_1$), by using the vector field $\partial\over\partial x$ (which spans 
$\xi$) and an inward normal to $\partial X_1$. 
These extend uniquely to a frame field $(u,v)$ that trivializes $TX_1$ 
(since $X_1\simeq \vee S^1$). 
Now recall that each 2-handle $h_i$ is a neighborhood of $D^2\times 0\subset 
i\real^2 \times \real^2$. 
We choose a convenient trivialization of its tangent bundle. 
The tangent and outward normal vector fields $\tau$ and 
$\nu$ to $S^1\subset D^2$ 
together trivialize $TD^2|S^1$ over $\real$. 
The frame field $(\tau,\nu)$ differs from the product frame field on 
$D^2\subset i\real^2$ by an element of $\pi_1(SO(2))$. 
Since $D^2$ lies in $i\real^2\subset \complex^2$, 
the field $(\tau,\nu)$ also forms a complex trivialization of $Th_i|S^1$,  
which differs from a product trivialization by the same element of 
$\pi_1(SO(2))$. 
Since $SO(2) \subset SU(2)\subset U(2)$ and $SU(2)$ is simply connected, 
$(\tau,\nu)$ extends to a complex trivialization $(\tau^*,\nu^*)$ over 
all of $h_i$. 
When we attach $h_i$ to $X_1$, $\tau$ is identified with a tangent 
vector field to $K_i$ and $\nu$ becomes inward normal to $\partial X_1$ 
in $TX_1$. 
Clearly, $\tau^*$ and $u$ span a complex line bundle $L$ (which agrees 
with $\xi$ in the box), while $\nu^*$ and $v$ can be fit together to 
span a complementary trivial line bundle. 
Thus, the desired cochain evaluated on $h_i$ is just the relative Chern 
number of $L$, which is given by the rotation number of $\tau$ in $\xi$ 
relative to $\partial\over\partial x$, or $r(K_i)$.
(To check the sign, recall that in the nullhomologous case, $r(K_i)$ 
was defined as a relative Chern number on a Seifert surface for $K_i$, 
so it is $c_1(L)$ evaluated on the corresponding closed surface in $X$.)\qed
\enddemo 

Note that the cochain specified above depends on our choice of the vector 
field used to define $r(K_i)$, and hence, on how we chose to represent 
$X_1$ as a picture. 
However, the cohomology class is well-defined, since changing 
the vector field modifies the cocycle by a coboundary. 
Clearly, the above procedure also yields $c_1(J)$ for Stein surfaces 
$(X,J)$ with infinite topology (where we define the rotation numbers $r(K_i)$ 
relative to any convenient vector field). 
Alternatively, note that while an infinite handle decomposition may have a 
more complicated link diagram than we have considered (due to the 
end-structure of the union of 0- and 1-handles), any finite 
subhandlebody can be put in standard form. 

Before proving Theorem~2.2, it will be convenient to introduce another 
model for $S^1\times S^2$. 
We begin by identifying $S^1\times B^3$ with a closed $\varep$-neighborhood,  
with $\varep=1$, of $0\times i\real$ in $\complex^2$, modulo 
translations by $0\times 2\pi i\zed$. 
The tight contact structure of complex lines tangent to $S^1\times S^2$ 
will be given by the 1-form $u_1\,du_2 - u_2\,du_1 + u_3\,du_4$ on 
$S^1\times S^2$, in the coordinates ($u_1+ iu_2$, $u_3 + iu_4$) on 
$\complex^2$. 
We pull this form back to $\real^3$ with cylindrical coordinates 
$(r,\theta,t)$, via stereographic projection 
from a plane through the origin, $\real^2\times \real\to (S^2-(0,0,-1))
\times \real \subset \complex^2$, 
obtaining the form $d\theta +{1-r^4\over 
4r^2}\,dt$ after rescaling the form by a positive function. 
This pulls back to the standard form $\alpha = dz+x\,dy$ on $\real^3$,  
via the covering map from $\real^3$ to $\real^3$ minus 
the $t$-axis given by setting 
$x= {1-r^4\over 4r^2}$, $y=t$, $z=\theta$. 
Thus, we have exhibited the standard contact structure on 
$S^1\times (S^2 -\{\text{poles}\})$ as being the standard structure 
$(\real^3,\xi)$ modulo translations by $2\pi (\zed\times\zed)$ in the 
$y$-$z$ plane. 
The spheres $p\times S^2$ are given by the annuli $y=$ constant (modulo 
translations in the $z$-direction), each compactified by two points 
(at $x=\pm\infty$). 
Now we can represent Legendrian links in $S^1\times S^2$ in the usual way 
by their front projections  into $\real^2/2\pi\zed^2$. 
The new move corresponding to Move~6 (Figure~9) is given by Figure~11 and 
its image under $180^\circ$ rotation  about the $y$-axis 
(corresponding to pushing the link through the transverse circles at 
$x=\pm\infty$). 
Note that we can use these moves to make a link disjoint from the top 
and bottom edges of the square; the correspondence with Theorem~2.2 $(n=1)$ 
should now be clear. 
\cfig{1.5}{rg-fig11.eps}{11}

\demo{Proof of Theorem  2.2} 
First, we show that an arbitrary 
embedding $f:S^2\to M$ into a tight contact 3-manifold can be perturbed so 
that its image $f(S^2)$ has a neighborhood contactomorphic to that of 
$p\times S^2\subset S^1\times S^2$. 
By \cite{Gi1}, it suffices to perturb $f$ so that $f(S^2)$ has the same 
characteristic foliation as $p\times S^2$. 
Tightness guarantees that the foliation on $f(S^2)$ has no closed leaves 
(which would be Legendrian unknots with $tb=0$). 
By \cite{E3} (see also \cite{E5}), we can make a $C^0$-small perturbation 
of $f$ so that $f(S^2)$ 
has only two singular points in its characteristic foliation.  
By a $C^1$-small perturbation near the singular points, we can arrange the 
characteristic foliation to be radial there (cf.\ \cite{E4}), so that it 
is globally determined (up to diffeomorphisms of $S^2$) by its monodromy 
map between the circles of unit tangent vectors at  the singular points. 
By an additional $C^1$-small perturbation near one singular point,  
we can change this 
monodromy by a diffeomorphism of $S^1$, so that the characteristic foliation 
agrees with that of $p\times S^2$ in $S^1\times S^2$ as required. 
To see this, consider the 
local model of the singularity given by one pole of $p\times S^2$ in 
$S^1\times S^2$, or in our model of the latter,  
the $x$-$z$ plane in the region $x\ge N$, $|y|<\varep$ 
in $(\real^3,\xi)/2\pi\zed^2$. 
In the annulus between the circles $y=0$ and $y=\varep/2$ in the 
$y$-$z$ plane modulo $0\times 2\pi \zed$, draw a flow realizing the given 
diffeomorphism of $S^1$. 
By composing with additional $2\pi$-rotations, we may assume the flow lines 
have arbitrarily negative slope everywhere, and that the slopes increase 
in absolute value with increasing $y$ for fixed $z$. 
Interpreting these flow lines as Legendrian curves in $\real^3$, we obtain 
a surface near the $x$-$z$ plane with the required characteristic 
foliation. 
(To make the perturbation $C^1$-small, identify $p\times S^2$ locally with 
the plane $z=0$ in $(\real^3, dz+r^2d\theta$) and rescale by 
$(r,z)\mapsto (tr,t^2z)$.) 

Now given $H$ with its ordered collection of $n$ 1-handles, let $S_1,\ldots,
S_n$ be the corresponding belt spheres. 
As above, we can assume 
that a neighborhood of $p\times S^2$ maps 
contactomorphically into $\partial H$ with $p\times S^2$ mapping to $S_i$. 
Similarly, we can find such a map into $\real^3$, and the 
composite realizes $S_i$ (with both choices of orientation) as 
the boundary of a contact 3-ball $B\subset \real^3$. 
We can choose $B$ to be nearly round, with $0\in \real^3$ lying on the 
equator of $\partial B$, so that each contactomorphism $\psi_t(x,y,z) = 
(tx,ty,t^2z)$, $0<t\le 1$, sends $B$ into itself fixing $0$. 
Now we can cut open $\partial H$ along each $S_i$ and glue in a pair of 
copies $B_{i1}$ and $B_{i2}$ of $B$, with the second index determined by 
the direction of the 1-handle. 
The result is $S^3$ with its unique tight contact structure. 
(An overtwisted disk could be isotoped off of the balls $B_{ij}$ 
using $\psi_t$, 
contradicting the tightness of $\partial H$.) 
Choose nonsingular points $p_{ij}\in \partial B_{ij}$ such that $p_{i1}$ and 
$p_{i2}$ correspond under the gluing map $\varphi_i:\partial B_{i1}\to 
\partial B_{i2}$ for $\partial H$. 
We can assume (after a contact isotopy of $\real^3$ fixing $\partial B$ 
setwise) that each identification $B_{ij}\approx B$ sends $p_{ij}$ to $0$. 
After a contact isotopy of $S^3$, we can assume that the points $p_{ij}$ 
are arranged along the edges of a box, say $p_{ij} = (0,j,i)\in \real^3$, 
with the outward normals to $B_{ij}$ at $p_{ij}$ directed into the box 
and $(d\varphi_i)_{p_{i1}} = \id_{\real^2}$. 
By conjugating 
the contact isotopy $\psi_t$, $\varep \le t\le 1$, 
with each identification $B_{ij}\approx B$, 
and extending to a contact isotopy 
of $S^3$ via Gray's Theorem, we can arrange for the balls $B_{ij}$ to 
be embedded in $\real^3$ by nearly linear maps. 
(We draw them as round balls although they will actually be ellipsoidal.) 
Since $(d\varphi_i)_{p_{i1}} = \id_{\real^2}$, we can find connected 
neighborhoods $U_{ij}$ of $p_{ij}$ in $\partial B_{ij}$ on which the 
gluing map $\varphi_i : U_{i1}\to U_{i2}$ is arbitrarily closely 
approximated by a translation in $\real^3$. 
If we apply 
a contactomorphism  of $S^3$ supported near $p_{i1}$, 
rotating the tangent space $360^\circ$ about a vertical axis,  before 
shrinking $B_{i1}$ with $\psi_t$, the vector field in $\partial H$ 
corresponding to $\partial\over\partial x$ in the box will be changed 
by a full twist. 
Since nowhere zero vector fields in $\xi$ are classified by 
$H^1(\partial H;\zed)$, the vector field corresponding to $\partial\over 
\partial x$ can be chosen arbitrarily. 
We have now represented $H$ as in Definition~2.1, with the required 
control of the auxiliary data.

Now let $L$ be a Legendrian link in $\partial H$. 
The image of $L$ in $S^3$ will be a tangle that can run into, around, 
behind or in front of the sphere $S_{ij}=\partial B_{ij}$ as 
in Figure~10, although we can 
assume that the front projection of $L$ is disjoint from each $p_{ij}$, 
and hence (after we 
shrink the neighborhoods $U_{ij}$) the projections of $L$ 
and $U_{ij}$ are disjoint. 
For each $i$, let $A_i\subset S_i$ be the arc whose interior is a leaf 
of the characteristic foliation, and whose image in $S_{i1}$ contains $p_{i1}$. 
By general position, we can assume that $L$ is disjoint from each $A_i$. 
Using our contactomorphism on a neighborhood of $S_i$ sending 
$S_i$ to $p\times S^2 \subset S^1\times S^2$, 
we can visualize $S_i-A_i$ as $\real\times 0 \times (0,2\pi)
\subset \real\times (-\varep,\varep) \times \real/2\pi\zed$ in our picture 
of $S^1\times S^2$. 
For any compact subset $C$ of $\partial H-A_i$, there is a contact isotopy 
$F_s$ compactly supported in $\real \times (-\varep,\varep)\times (0,2\pi)$ 
that moves $C$ to a set whose intersection with $S_i$ maps into 
$U_{i1}$ in $S_{i1}$. 
To construct this, start with a suitable isotopy $f_s$ in $(-\varep,\varep)
\times (0,2\pi)$ preserving the direction of ${\partial\over\partial z}$ 
(see Figure~12). 
There is a unique contact isotopy of 
$\real\times (-\varep,\varep)\times (0,2\pi)$ 
projecting to $f_s$, since the $x$-coordinate of each point is determined 
by the slope of the corresponding contact plane, hence controlled by $df_s$. 
To obtain compact support, 
truncate the isotopy for sufficiently large values of $|x|$ 
by applying Gray's Theorem.  
Taking $C=L$, we isotope $L$ so that it 
intersects the spheres $S_{ij}$ only in $U_{ij}$.  
By prechoosing $U_{ij}$ and $\varep$ to be sufficiently small, we can assume 
that all points in our neighborhood of $S_i$ that project onto $U_{ij}$ 
in our original planar diagram of $\partial H$ (Definition~2.1) will lie 
in a preassigned neighborhood of the top and bottom edges of Figure~12, with 
a narrow range of $x$-coordinates. 
Thus, the only points of $L$ that may project onto $U_{ij}$ (after the 
above isotopy of $L$) will lie on small arcs intersecting $S_i$ in Figure~12. 
In the original diagram, these arcs will be parallel, emerging horizontally 
from $S_{ij}$ and then doubling back with nowhere zero slope to follow 
the characteristic folation of $S_{ij}$. 
By increasing the slope if necessary, we can arrange that $L$ project onto 
$U_{ij}$ only at the endpoints of the tangle. 
Let $\ell_i$ denote the 
horizontal line in $\real^3$ passing through the centers of the 
balls $B_{i1}$ and $B_{i2}$, with the segment between the balls deleted.  
By general position, we can assume that $L$ is disjoint from each $\ell_i$. 
It is now easy to construct a contact isotopy of $\real^3$ (for example, 
via a planar projection preserving $\partial\over\partial z$) 
that pushes $L$ away from each $\ell_i$ and then entirely into the box, 
i.e., into standard form. 
Note that it is crucial here that $L$ not run in front of or behind $U_{ij}$. 
(Consider the difficulties inherent in Figure~10.) 
\cfig{2.0}{rg-fig12.eps}{12}

To prove the final assertion of the theorem, 
suppose that $\{L_t \mid 0\le t\le 1\}$ is an isotopy of 
Legendrian links in $\partial H$, with $L_0$ and $L_1$ in standard form. 
We wish to show that $L_0$ and $L_1$ are related by Moves~1-6. 
Let $T\subset (0,1)$ be the set of parameter values $t$ for which $L_t$ 
intersects some Legendrian curve $A_i$ or $\ell_i$. 
The homogeneity of contact structures allows us to construct enough 
deformations of a Legendrian isotopy to apply transversality theory, 
so we may assume that 
$T$ is finite, for $t\in T$, 
$L_t\cap \bigcup_i(A_i\cup \ell_i)$ is a single point, 
and near each such intersection the isotopy has a canonical form. 
After an additional modification, we obtain a small $\varep$ such that 
for each $t\in T$, the links with parameter between $t\pm \varep$ are 
identical except for a canonical push across some $A_i$ or $\ell_i$. 
The complement of an $\varep$-neighborhood of $T$ in $[0,1]$ is a finite 
collection of intervals $[t_1,t_2]$, each of which parametizes 
a Legendrian isotopy from $L_{t_1}$ 
to $L_{t_2}$ in $\partial H$ that is disjoint from each $A_i$ and $\ell_i$. 
The procedure of the previous paragraph simultaneously pushes each isotopy 
into the box, after we suitably enlarge $T$ to allow for intersections 
with $\ell_i$ created by the first isotopy $F_s$. 
(Note that while we previously disallowed links projecting onto any $p_{ij}$, 
a transverse pass across $p_{ij}$ does not affect the construction.) 
Each $L_t$ is sent to a Legendrian link $L'_t$ in the box, and we can assume 
that the endpoints $L'_{t_k}$ of each isotopy are in standard form. 
Now if we identify opposite edges of the box at the 1-handles, we essentially 
have an ordinary front projection of isotopies of Legendrian links,  
together with a set of distinguished vertical arcs where we have glued. 
Thus, we can reduce each isotopy between $L'_{t_1}$ and $L'_{t_2}$ to a 
sequence of Moves~1-3 (cf.\ \cite{S}), 
together with Moves~4 and 5 to account for sliding 
cusps and double points over 1-handles. 
We will complete the proof by showing 
that for each $t\in T$,  
the links $L'_{t-\varep}$ and $L'_{t+\varep}$ are related by 
Moves~1-6, as are the links $L_t$ and $L'_t$ for $t=0,1$.  

Suppose that 
$t\in T$ corresponds to an intersection with some $\ell_i$. 
At the intersection, we can take a local model of the front projection 
of $L_t$ (or $F_1(L_t)$) to 
be a parabola tangent to $\ell_i$,  
and the isotopy from $L_{t-\varep}$ to 
$L_{t+\varep}$ to correspond to a vertical translation. 
Tracing through the above procedure for producing $L'_{t\pm\varep}$, 
we see that these links will differ by Move $6'$ of Figure~13. 
Figure~14 shows how to reduce Move~$6'$ to Moves~1-6. 
The same argument, with the isotopy $F_s$ of the previous paragraph taken 
to be the identity, shows that two links in standard form are related by 
Moves~1--6 if they are related by a Legendrian isotopy in $\partial H$ 
disjoint from each $S_{ij}- U_{ij}$. 
It follows that $L_t$ and $L'_t$ are so related for $t=0,1$. 
\cfig{1.5}{rg-fig13.eps}{13}
\cfig{3.5}{rg-fig14.eps}{14}

For the remaining values $t\in T$, we have $L_t$ intersecting some $A_i$. 
Now we must analyze the effect of the isotopy $F_s$ of Figure~12 on 
pictures in standard form. 
Recall that the top and bottom edges of Figure~12 are identified, and that 
$A_i$ projects to the identified endpoints of the vertical arc 
representing $S_i$. 
Allowing isotopies of $L'_{t\pm\varep}$ disjoint from $S_{ij}-U_{ij}$, 
and reversing the sign of $\varep$ if necessary, we can assume that 
$L_{t-\varep}$ is given near $L_t\cap A_i$ by a horizontal arc through a 
point $p_-$ in $S_i$ near the top edge of Figure~12, and that this curve 
is fixed by $F_s$. 
Similarly, we can obtain $L_{t+\varep}$ from $L_{t-\varep}$ by pushing 
the point of intersection upward through $A_i$, to a point $p_+$ near 
the bottom of the figure. 
Thus the isotopy $F_s$ will fix $p_-$ and move $p_+$ to the top of the  
picture. 
We can assume that the maps from a neighborhood of $L_t\cap A_i$ 
in Figure~12  to 
$U_{i1}$ and $U_{i2}$ in Figure~8 are sufficiently 
well approximated by translations in $\real^3$  that 
$F_1$ aligns the points of $L_{t\pm\varep}\cap S_{ij}$ vertically in 
$S_{i1}$ and $S_{i2}$, with $F_1(p_+)$ on the bottom and $p_-$ on top. 
To translate Figure~12 into a standard picture of $\partial H$, we 
first draw the $z$-axis of Figure~12 as a transverse equatorial circle in 
each $S_{ij}$, projecting to a figure-8. (See Figure~15.) 
The arrows indicate the direction of $F_s$ ($\partial\over\partial z$ 
in Figure~12). 
Note that this orients the equators of both spheres counterclockwise, 
and is consistent with an orientation-reversing gluing map $\varphi_i :
S_{i1}\to S_{i2}$ interchanging the poles. 
(We can assume that $\varphi_i$ identifies the equators by $180^\circ$ 
rotation about the $z$-axis.) 
The isotopy $F_s$ will push $p_+$ counterclockwise once around each 
equatorial curve. 
To determine its effect on $L_{t+\varep}$, we examine the curves in 
Figure~12. 
On each segment where the second derivative has large magnitude, the 
curve will be traveling nearly parallel to the characteristic  foliation 
of $S_i$ (a left-handed spiral in Figure~15). 
Where the first derivatives have large magnitude, the curves travel 
crosswise to the foliation near the poles of $S_i$. 
By checking orientations, we see that these parts of the curves will 
lie above both northern hemispheres in Figure~15. 
Figures~15 and 16 (and the image of Figure~16 under $180^\circ$ rotation  
about the $z$-axis) show the isotopy $F_s$. 
Completing our contact isotopy to obtain $L'_{t\pm\varep}$, we see that 
these two links differ by the move shown in Figure~17. 
(We can remove the extra spirals by Move~1. 
Note that the other curves in Figure~12 lie farther from $S_{ij}$ 
in Figures~15 and 16, as does the rest of $L_{t\pm\varep}$, so these do 
not interfere, and can be pushed away toward the center of the box in 
Figure~17 by Moves~2 and 3.) 
Figure~18 provides the reduction from this move to Moves~1-6.\qed
\enddemo
\cfig{1.75}{rg-fig15.eps}{15}
\cfig{3.0}{rg-fig16.eps}{16}
\cfig{2.2}{rg-fig17.eps}{17}

\remark{Remarks} 
It is easily checked that Move $6'$ (Figure~13) is equivalent to Move~6 
(Figure~9) in the presence of Moves~1-3 (Figure~3). 
Figure~18 generalizes to give a derivation of the move in Figure~19, which 
represents a collection of strands being looped around and under a 1-handle. 
Of course, this contact isotopy preserves the canonical framing of each 
knot $K$, as well as $r(K)$, but the integer $tb(K)$ changes in general. 
A related example is given by Figure~20, which shows how we can convert 
a right crossing to a left one, by a smooth (not contact) 
isotopy that adds a {\it right\/} 
twist to the canonical framing of $K$. 
Such tricks can be useful for putting Stein structures on preassigned 
smooth 4-manifolds. 
\endremark 

\cfig{1.70}{rg-fig18.eps}{18}
\cfig{0.85}{rg-fig19.eps}{19}
Implicit in the proof of Theorem 2.2 are nontrivial self-equivalences of 
$H$ (up to Stein-homotopy \cite{E6}) that are smoothly 
isotopic to the identity. 
We found a contactomorphism $\psi_i$ that twisted the vector field in 
$\partial H$ corresponding to $\partial\over\partial x$, 
by rotating $B_{i1}$ through $360^\circ$. 
The proof of the theorem shows that $\psi_i$ acts on Legendrian links as 
in Figure~21, up to contact isotopy. 
(For example, start with horizontal arcs in Figure~12, then apply a vertical 
Dehn twist on one side.) 
Figure~19  verifies that these moves for $\psi_i$ and $\psi_i^{-1}$ 
are indeed inverses up to contact isotopy. 
In our alternate picture of $S^1\times S^2$, $\psi_1$ is given by a Dehn twist 
on the torus, $\psi_1(x,y,z) = (x-1,y,y+z)$ on $\real^3/0\times 2\pi \zed^2$. 
Now $\psi_i\circ \psi_i$ is smoothly isotopic to the identity, since 
$\pi_1(SO(3))\cong \zed_2$. 
However, no nonzero power of $\psi_i$ is contact isotopic to the identity. 
This is because $\psi_i$ twists any nonzero vector field $v$ in $\xi$,  
and $\pi_1 (SO(2))\cong \zed$. 
For example, if $K$ runs once over the 1-handle, then $\psi_i$ will 
increase the rotation number of $K$ by 1, so $\psi_i^n(K)$ will not 
be contact isotopic to $K$ for $n\ne 0$. 
Clearly, we obtain other automorphisms of $H$, which are nontrivial 
on $\pi_1$, by changing the other auxiliary data. 
Changes of order or directions of the 1-handles correspond to various 
permutations of the balls $B_{ij}$. 
We can reverse the orientation of $\xi$ by rotating 
about the $x$- or $y$-axis. 
Allowing other handle structures on $H$ increases the symmetries 
(by 1-handle slides). 
\cfig{1.75}{rg-fig20.eps}{20}
\cfig{2.5}{rg-fig21.eps}{21}
\newpage 

\subhead 3. Exotic Stein surfaces\endsubhead 

In this section, we give a simple characterization of those manifolds that  
are homeomorphic to Stein surfaces. 
We obtain many new examples of Stein surfaces, most of which are in some 
sense ``exotic'' as smooth manifolds. 
For example, we obtain an uncountable family of nondiffeomorphic Stein 
surfaces, all of which are homeomorphic to $\real^4$ (``exotic'' $\real^4$'s). 

\proclaim{Theorem 3.1} 
An open, oriented, topological 4-manifold $X$ is homeomorphic to a Stein  
surface if and only if it is the interior of a 
(possibly infinite) handlebody $H$ without handles 
of index $>2$. 
If so, then any almost-complex structure on $X$ (up to homotopy)
is induced by an 
orientation-preserving homeomorphism from a Stein surface. 
\endproclaim 

To clarify the last statement, we must specify the meaning of an 
almost-complex structure (up to homotopy) on the topological manifold $X$. 
Any 4-dimensional topological handlebody can be canonically 
smoothed, since the gluing maps are flat topological embeddings of 
3-manifolds, and these are always uniquely smoothable. 
Thus, we can define almost-complex structures on $X$ in the usual 
way relative to a fixed smoothing. 
To ensure that homeomorphisms induce canonical correspondences of such 
structures, however, requires more work. 
Recall (e.g., \cite{KM}, \cite{T}) that an almost-complex structure on 
a smooth 4-manifold determines a $\spinc$-structure. 
In our case $H^i(X;\zed) =0$ for $i>2$, so the correspondence is bijective 
(since $\Spinc(4)/U(2)\approx S^3$ is 2-connected), and it suffices to 
define $\spinc$-structures on oriented topological 4-manifolds. 
Every such manifold has a topological tangent bundle with structure 
group $\STop(4)$, the group of orientation-preserving homeomorphisms of 
$(\real^4,0)$ \cite{KS}. 
Since the inclusion $SO(4)\hookrightarrow \STop (4)$ is double covered by 
an inclusion  $\Spin (4) \hookrightarrow \SpinTop (4)$, we can replace 
$\Spin (4)$ by $\SpinTop(4)$ in $\Spinc (4) = (S^1\times \Spin (4))/\zed_2$, 
and construct a theory of topological $\spinc$-structures that reduces 
to the usual theory in the presence of a smooth structure \cite{G3}. 
Orientation-preserving homeomorphisms will preserve these structures in 
the same manner that diffeomorphisms do, so we obtain a 
homeomorphism-invariant notion of almost-complex structures on $X$. 
Absolute and relative Chern classes $c_1$ will be preserved by 
orientation-preserving homeomorphisms, since they are determined by 
the complex line bundles associated to the corresponding topological 
$\spinc$-structures. 
(A related argument applies to homotopy equivalences \cite{G3}.) 

To enumerate almost-complex structures on $X$, we smooth $X$ as $\rint H$. 
Since $SO(4)/U(2)\approx S^2$ is simply connected, almost-complex structures 
exist, and we can assume that they all agree on the union $X_1$ of 
0- and 1-handles. 
After fixing a complex trivialization $\tau$ there, we obtain a relative 
Chern class in $H^2(X,X_1;\zed)$ for each almost-complex structure. 
These will reduce modulo~2 to the relative Stiefel-Whitney class 
$w_2(X,\tau)\in H^2(X,X_1;\zed_2)$. 
Since $\pi_2 (SO(4)/U(2))\cong \zed$, obstruction theory produces a 
surjection from integer lifts of $w_2(X,\tau)$ to almost-complex structures,  
via $c_1$.  
(The map fails to be injective, since there will be nontrivial self-homotopies 
of the structure over $X_1$, changing $\tau$ by any even number of twists. 
In fact, the set of almost-complex structures on $X$ is affinely 
$H^2(X;\zed)$, with difference cocycles given by half the difference in 
relative Chern classes.) 

\demo{Proof}
Clearly, any manifold homeomorphic to a Stein surface has the required 
handle structure. 
(See Theorem~1.3.) 
For the converse, we start with $X=\rint H$ as in the theorem, and construct 
a Stein surface homeomorphic to $X$. 
We assume (as above) that $H$ is a smooth handlebody, and fix a trivialization 
$\tau$ of the unique (up to homotopy) almost-complex structure on $X_1$, 
the intersection of $X$ with the 0- and 1-handles of $H$. 
We pick an integer lift $c\in H^2(X,X_1;\zed)$ of $w_2(X,\tau)$, and 
arrange for the homeomorphism to map $c$ onto
the Chern class of the Stein surface 
(relative to the corresponding trivialization on the preimage of $X_1$). 
The theorem then follows from the discussion above. 

Consider  a handlebody $H'$ obtained from $H$ by removing each 2-handle 
$h_i$ and regluing it along the same attaching circle, but with an even 
number $2k_i$ of left (negative) twists added to its framing. 
By Eliashberg's Theorem (1.3), $\rint H'$ will be Stein if each $k_i$ is 
sufficiently large. 
Clearly, $H'$ is canonically homotopy equivalent $\rel X_1$ to $H$ and $X$, 
with $w_2(H',\tau)$ mapping to $w_2(X,\tau)$ 
since the numbers of twists were even. 
Thus, $c_1(H',\tau)$ and $c$ differ by an even class in $H^2(X,X_1;\zed)$. 
Since adding 2 left twists to the framing of $h_i$ allows us to change 
$c_1(H',\tau)$ by $\pm2$ on $h_i$ (cf.\ Figure~7 and Proposition~2.3), 
we can assume (for each 
$k_i$ sufficiently large) that $c_1(H',\tau)$ maps to $c$ under the homotopy 
equivalence.

Our homotopy equivalence between $H'$ and $X$ does not preserve the 
intersection pairing, a defect which we will remedy at the expense of 
increasing $\pi_1$, as follows. 
Close inspection of Eliashberg's paper \cite{E2} shows that we can put 
self-plumbings of either sign in each of the 2-handles of $H'$, 
without changing 
the gluing maps, and after smoothing, still have a manifold $H''$ whose 
interior is Stein. 
(A self-plumbing is performed by choosing a pair of disjoint disks 
$D,D'\subset D^2$ and gluing $D\times D^2$ to $D'\times D^2$ in the 
2-handle $D^2\times D^2$, by a diffeomorphism interchanging the factors. 
Since the 2-handles are constructed by Lemmas~3.5.1 and 3.4.3 of \cite{E2}, 
we can assume that part of each handle is given as an arbitrarily 
narrow $\varep$-neighborhood of $D^2\times 0$ in $i\real^2\times \real^2$. 
We can then do the self-plumbings, smoothing by Lemma~3.4.4 to preserve the 
Stein property.) 
Alternatively, we may do each plumbing explicitly in a link picture by 
adding a 1-handle as in Figure~22 (ignoring the dashed curves) 
(cf.\ \cite{C}, \cite{GS}). 
Suppose we obtain $H''$ from $H'$ by adding $k_i$ self-plumbings (counted 
with sign) to each $h_i$, transforming it into a {\it kinky handle\/} $h''_i$. 
There is a canonical local diffeomorphism $\varphi :H'\to H''\rel X_1$ 
that induces isomorphisms on $H_2$ and $H^2$. 
Since $\varphi$ preserves the almost-complex structures, the Chern 
class $c_1(H'',\tau)$ corresponds to $c\in H^2(X,X_1;\zed)$ under the 
obvious isomorphism. 
Furthermore, the given isomorphism $H_2(H'';\zed)\cong H_2(X;\zed)$ preserves 
the intersection pairing. 
In fact, by the most natural way of defining the attachment of kinky 
handles along framed circles, the handles $h''_i$ are actually attached along 
the framed curves defining the original handlebody $H$ \cite{C}. 
To understand  this, 
recall that for closed surfaces generically immersed in 4-manifolds, the 
homological intersection number differs from the normal Euler number by 
twice the signed number of self-intersections. 
Thus, the above convention for attaching kinky handles, 
correcting the framing by twice the signed number of self-plumbings, 
guarantees that the 
resulting intersection pairing only depends on the framed link, and 
not on the numbers of self-plumbings of the kinky handles. 
Alternatively, observe that if we add self-plumbings to $h_i$ using 
Figure~22, each $\pm$ self-plumbing will increase the Thurston-Bennequin 
invariant of the attaching circle by $\pm2$ and leave its rotation number 
unchanged, as required. 

\cfig{2.75}{rg-fig22.eps}{22}

Finally, we eliminate the unwanted extra $\pi_1$ by extending the kinky 
handles $h''_i$ to Casson handles. 
According to Casson \cite{C}, each kinky handle $h''_i$ has a canonical 
framed link in its boundary, such that if we add 2-handles along the 
framed link, $h''_i$ will be transformed into a standard 2-handle, and the  
above-mentioned natural framing on the attaching circle of $h''_i$ will 
correspond to the product framing on the 2-handle. 
Thus, attaching 2-handles to all kinky handles $h''_i$ in this manner 
would transform $H''$  back into $H$. 
In general, we cannot do this in the Stein setting, but our previous 
argument shows that we can add kinky handles to these framed links to obtain 
a new manifold whose interior is Stein. 
(In fact, the relevant circles appear in Figure~22 as dashed curves with 
the zero framing, so we can add any kinky handles with more positive than 
negative self-plumbings.) 
We iterate the construction, adding a third layer of kinky handles onto the 
second layer, and continue, to construct a manifold $\hatH$ with infinitely 
many layers of kinky handles. 
Clearly, $\rint \hatH$ is Stein. 
But for each infinite stack of kinky handles (starting with some $h''_i$), 
the interior union the attaching region is (by definition) a {\it Casson 
handle\/} \cite{C}. 
Freedman \cite{F} proved that any Casson handle is homeomorphic to an 
open 2-handle $D^2\times \real^2$ such that the natural framings of the 
attaching circles correspond. 
Thus, the Stein surface $\rint\hatH$ is 
homeomorphic $\rel X_1$ to $\rint H =X$. 
The restriction map $H^2(\hatH,X_1;\zed)\to H^2(H'',X_1;\zed)$ 
is an isomorphism 
preserving $c_1$, so the given homeomorphism $\rint\hatH \approx X$ maps 
$c_1(\hatH,\tau)$ onto $c$ as required.\qed
\enddemo 

A related observation applies in the compact setting.  
If $H$ is a 4-dimensional compact handlebody without handles of index $>2$, 
$X_1$ denotes its 1-skeleton, 
and $c\in H^2(H,X_1;\zed)$ and $\tau$ on $X_1$ specify an 
almost-complex structure on $H$, then 
there is a compact Stein surface $X$ with boundary, and a homotopy 
equivalence $\varphi:X\to H\rel X_1$ preserving the intersection pairing 
and with $\varphi^* (c)=c_1(X,\tau)$. 
To see this, simply note that for a Legendrian knot $K$, $tb(K)$ can be 
increased by any even number without changing $r(K)$, 
by forming the connected sum with a 
knot with sufficiently large $tb$ and $r=0$. 
For example, we can realize any finitely presented group as $\pi_1(X)$ or 
any finite rank, symmetric $\zed$-bilinear form as the intersection 
pairing of such a compact $X$. 

\proclaim{Corollary 3.2} 
Any smooth, closed, connected, oriented 4-manifold contains 
a smooth, finite wedge 
of circles whose complement is homeomorphic to a Stein surface. 
Any smooth (resp.\ topological) 
open, oriented 4-manifold contains a smooth (resp.\ locally flat) 
1-complex whose complement is 
homeomorphic to a Stein surface. 
\endproclaim 

\demo{Proof} 
In the topological case, the manifold admits a smooth 
structure \cite{Q}, \cite{FQ}. 
Now in either case, there is a handle decomposition, and the required 
1-complex is dual to the 3- and 4-handles.\qed
\enddemo 

Since the Stein surfaces constructed in the proof of Theorem~3.1 are built 
with Casson handles, the underlying smooth manifolds will typically be 
``exotic'' in some sense. 
It seems likely that any $X$ will admit uncountably many diffeomorphism 
types of such homeomorphic Stein surfaces, and that no such smooth 
manifold can admit a proper Morse function with finitely many critical 
points (provided that $X$ is not homeomorphic to 
the interior of $B^4 \cup 1$-handles). 
We illustrate this behavior with a concrete example. 

\proclaim{Theorem 3.3} 
For each integer $n$, let $H_n$ denote the compact 2-disk bundle over $S^2$ 
with Euler number $n$. 
Then $\rint H_n$ is (orientation-preserving) homeomorphic to a Stein 
surface $V_n$ that contains no smoothly embedded sphere generating its 
homology and realizes any preassigned almost-complex structure. 
For $n=\pm1$, there are uncountably many diffeomorphism types of 
such manifolds $V_{\pm1}$, none of which admit proper Morse functions with 
finitely many critical points. 
\endproclaim  

In contrast, the manifolds $\partial H_0 = S^1\times S^2$, $\partial H_{\pm1} 
= S^3$ and $\partial H_{\pm2} =\real P^3$ admit unique tight contact 
structures, and these are uniquely fillable (up to blowing up) by 
$S^1\times B^3$, $B^4$ and $H_{-2}$, respectively (by work of Gromov \cite{Gro} 
and Eliashberg, see \cite{E3}). 
Thus, $H_0 = S^2\times D^2$, $H_1= \complex P^2 -\rint B^4$ and $H_2$ are 
not diffeomorphic to Stein surfaces with boundary (or even symplectic 
manifolds with convex boundaries), although their interiors admit exotic 
smooth structures that are Stein. 
Clearly, we could construct many other examples of manifolds with these 
boundaries, whose interiors are homeomorphic to Stein surfaces. 
(Consider closed manifolds minus $\rint B^4$, $\rint S^1\times B^3$ or 
$\rint H_{\pm2}$.) 

\demo{Proof} 
The manifold $H_n$ is obtained from $B^4$ by gluing a 2-handle to an unknot 
in $\partial B^4$ with framing $n$. 
For any Casson handle $CH$, let $V_n(CH)$ be the interior of the manifold 
obtained by gluing $CH$ to $B^4$ along an $n$-framed  unknot. 
Then by the proof of Theorem~3.1, $V_n(CH)$ will be realized as a Stein 
surface homeomorphic to $\rint H_n$ and realizing the preassigned 
almost-complex structure, provided that $CH$ has a suitable 
excess of positive self-plumbings. 

There are many known examples of smooth, simply connected 
4-manifolds $M$ whose intersection 
pairing contains a subspace with pairing 
$$\left[ \matrix n&1\cr 1&\text{even}\endmatrix \right]\ ,$$ 
such that the class $\alpha$ with square $n$ cannot be represented by a 
smoothly embedded sphere. 
(For $n\ge0$ see, for example, \cite{FM} Chapter~6, Corollary~4.2. 
For $n<0$, simply reverse orientation.) 
By Casson's Embedding Theorem \cite{C}, however, we can find a Casson handle 
$CH$ such that $V_n (CH)$ embeds in $M$ representing $\alpha$. 
We can always add self-plumbings and layers of kinky handles to $CH$, 
so that it becomes suitably positive and $V_n(CH)$ admits a Stein structure 
as above. 
It cannot contain a smooth sphere generating its homology, since this would 
also represent $\alpha$ in $M$. 

Now suppose $n=1$. 
Then we may choose $M$ so that the orthogonal complement of $\alpha$ in the 
intersection pairing is nonstandard and negative definite \cite{G1}.  
As in Freedman \cite{F}, we can construct a nested family $\{CH_c\}$ of 
Casson handles inside $CH$, indexed by a Cantor set. 
At each stage of the construction, we add positive self-plumbings wherever  
necessary so that each $V_1(CH_c)$ admits a Stein structure as above. 
If any two of the nested manifolds $V_1(CH_c)$ in $M$ were diffeomorphic, 
then a standard argument \cite{G2} would allow us to contradict the 
Periodic End Theorem of Taubes. 
If any one had a proper Morse function with finitely many critical points, 
then its end would be smoothly collared by a 3-manifold crossed with $\real$, 
and the same argument would apply. 
For $n=-1$, the same argument works in the manifold $\overline{M}$.\qed 
\enddemo 

We reformulate our result about gluing Casson handles in the context of 
Legendrian link presentations of Stein surfaces (Section~2). 
A kinky handle with $k_+$ positive and $k_-$ negative self-plumbings 
can be attached along a Legendrian knot $K$ in $\# nS^1\times S^2$ 
with framing $m$ (after allowing a $C^0$-small perturbation of $K$ to lower 
$tb$ if necessary), provided that $m\le tb (K)-1+2(k_+- k_-)$. 
This kinky handle can be extended to a Casson handle, and the only 
restriction on this extension is that for each of the additional kinky 
handles we must have $k_+ >k_-$. 
The rotation number $r(K)$ (after the perturbation of $K$ to achieve
equality in the above formula) contributes to the Chern class as if the 
Casson handle were an ordinary 2-handle. 
We now obtain the following theorem, that some exotic $\real^4$'s admit 
Stein structures. 

\proclaim{Theorem 3.4} 
There are uncountably many diffeomorphism types of Stein surfaces 
homeomorphic to $\real^4$. 
There is a Stein exotic $\real^4$ that can be built with two 
1-handles, one 2-handle and a single Casson 
handle $CH$ with only one kinky handle at each stage (Figure~25). 
\cfig{2.5}{rg-fig23.eps}{23}
\endproclaim 

\demo{Proof} 
In \cite{BG}, Bi\v zaca  and the author exhibit a particularly simple 
exotic $\real^4$. 
This is the interior $R$ of the manifold shown in Figure~23, which is 
essentially Figure~1 of \cite{BG}. 
The circles with dots represent 1-handles (cf.\ \cite{K},\cite{GS}), 
the solid curve 
represents a 2-handle, and the dashed curve is where we attach the Casson 
handle $CH$ with a single, positive self-plumbing in each kinky handle. 
Both framings are 0, as indicated. 
It is routine to verify that Figure~24 represents the same manifold 
$R$, where 
we are now attaching the handle and Casson handle to a 0-framed Legendrian 
link. 
(Change Figure~24 in the obvious way to represent the 1-handles by 
circles with dots. 
Then isotope to Figure~23.) 
Since the dashed curve has $tb=0$, we can attach the Casson handle $CH$ 
with framing 0 as required. 
However, the solid curve has $tb=-2$, so we must increase its canonical 
framing by 3 units by a smooth isotopy. 
This is easily accomplished by passing two strands around 1-handles 
(variations of the trick in Figure~20), resulting in Figure~25 of $R$, a 
Stein manifold homeomorphic but not diffeomorphic to $\real^4$. 
In \cite{BG} it was also observed (as in \cite{DF}) that $R$ contains 
an uncountable family of exotic  $\real^4$'s $\{R_c\}$ indexed by a Cantor 
set, produced using a nested family of Casson handles $\{CH_c\}$ 
in $CH$ as in 
the proof of Theorem~3.3. 
As before, we can arrange the family $\{CH_c\}$ so that the manifolds $R_c$ 
(given by Figure~25 with $CH_c$ in place of $CH$) are all Stein. 
As in \cite{DF}, the family $\{R_c\}$ represents uncountably many 
diffeomorphism types.\qed
\cdubfig{1.7}{rg-fig24.eps}{24}{2.0}{rg-fig25.eps}{25}
\enddemo 

The Stein exotic $\real^4$'s described here can all be smoothly embedded in 
the standard $\real^4$. 
There is another type of exotic $\real^4$ that cannot be so embedded 
\cite{G2}. 
It is still an open (and apparently difficult) question whether any of these 
larger exotic $\real^4$'s admit Stein structures. 

\subhead 4. Invariants of 2-plane fields on 3-manifolds\endsubhead 

In this section, we define a complete set of invariants for distinguishing 
homotopy classes of oriented 2-plane fields 
(or equivalently, nowhere zero vector 
fields or ``combings'') on oriented 3-manifolds. 
We show how to compute the invariants for the boundary of a compact Stein 
surface presented in standard form as in Section~2. 
We obtain various corollaries, including invariance of the rotation number (up 
to sign) and Thurston-Bennequin invariant for contact 3-manifolds obtained 
by surgery on Legendrian knots (Corollary~4.6). 

At first glance, the classification of oriented 2-plane fields $\xi$ on an 
oriented 3-manifold $M$ seems to be easy with modern techniques. 
If we fix a trivialization of the tangent bundle $TM$, the problem 
becomes equivalent to classifying maps $\varphi :M\to S^2$ up to homotopy, 
and this latter problem was solved by Pontrjagin around 1940 \cite{P}. 
The difficulty is that the invariants depend on the choice of trivialization 
of $TM$. 
In fact, one obtains both a 2-dimensional obstruction (which is not 
necessarily determined by the Chern class  $c_1(\xi)$) and a 3-dimensional 
obstruction. 
However, if we classify only up to isomorphisms of $TM$, the problem reduces 
to classifying abstract, stably trivial 
2-plane  bundles on $M$, and this is accomplished by 
the Chern class. 
Thus, if we allow the trivialization of $TM$ to vary, all of our invariants 
except for the Chern class will be lost. 
The main problem, then, is how to construct invariants that capture the 
information contained in the homotopy class of $\xi$, but can be manipulated 
without the awkward task of keeping track of a trivialization of $TM$. 
We solve this problem by expressing the 3-dimensional invariant 
in terms of a 4-manifold bounded by $M$, and the 2-dimensional invariant 
in terms of a spin structure. 
Our strategy is first to understand the invariants relative to a fixed 
trivialization (Proposition~4.1). 
Then we discuss invariants $\Theta$ and $\theta$ that capture the 
3-dimensional obstruction up to (at most) a $2:1$ ambiguity (4.2--4.6 and 
preceding text). 
To resolve the ambiguity, we must understand the 2-dimensional obstruction 
$\Gamma$, which we analyze in 4.7--4.14. 
We give an explicit formula for $\Gamma$ that applies to the contact 
boundary of any compact Stein surface in standard form (Theorem~4.12) 
and give several applications. 
Finally, we exhibit the full 3-dimensional obstruction $\widetilde\Theta$, 
with further properties of $\theta$ and a sample of applications (4.15--4.20). 
We show that for Stein boundaries, the invariants $\Gamma$ and  
$\widetilde\Theta$ (or $\theta$) are independent from each other and from 
$c_1$, by exhibiting holomorphically fillable contact  structures 
with $c_1=0$ on a fixed 3-manifold, and showing that they are distinguished 
by $\theta$ but not $\Gamma$ (Corollary~4.6) or vice versa (Example~4.14). 

The first step in constructing the invariants is to understand the 
classification for a fixed trivialization of $TM$. 
For our purposes, it is convenient to use differential topology rather 
than Pontrjagin's original obstruction-theoretic method of proof \cite{P}. 
(See also Kuperberg \cite{Ku} for a modern treatment via homotopy theory.) 
For an oriented 2-plane field $\xi$, let $d(\xi)\in\zed$ denote the 
divisibility of the Chern class, so that $c_1(\xi)$ equals $d(\xi)$ times 
a primitive class in $H^2(M;\zed)$ modulo torsion, and $d(\xi) =0$ 
if $c_1(\xi)$ is of finite order. 

\proclaim{Proposition 4.1} 
Let $M$ be a closed, connected 3-manifold. 
Then any trivialization $\tau$ of the tangent bundle of $M$ determines  
a function $\Gamma_\tau$ sending homotopy classes of 
oriented 2-plane fields $\xi$on $M$ into $H_1(M;\zed)$, and for any $\xi$, 
$2\Gamma_\tau(\xi)$ is Poincar\'e dual to $c_1(\xi)\in H^2 (M;\zed)$. 
For any fixed $x\in H_1(M;\zed)$, the set 
$\Gamma_\tau^{-1}(x)$ of classes of 2-plane 
fields $\xi$ mapping to $x$ has a canonical $\zed$-action and is isomorphic 
to $\zed/d(\xi)$ as a $\zed$-space. 
\endproclaim 

\noindent Note that $d(\xi)$ equals twice the divisibility of $x$, so it 
is independent of the choice of $\xi\in\Gamma_\tau^{-1}(x)$. 

\demo{Proof} 
Since $\tau$ identifies 
each tangent space of $M$ with $\real^3$ (with the 
standard orientation and inner product), 
oriented 2-plane fields $\xi$ on $M$ correspond 
to their orthogonal unit vector fields, or to maps $\varphi_\xi :M\to S^2$. 
By the Thom-Pontrjagin construction \cite{M2}, 
homotopy classes of such maps correspond 
bijectively to framed cobordism classes of framed links in $M$, and the 
correspondence sends $\varphi_\xi$ to $\varphi_\xi^{-1}(p)$ for any regular 
value $p\in S^2$, framed by pulling back an oriented basis of $T_pS^2$. 
Since $M$ is oriented by $\tau$, $\varphi_\xi^{-1}(p)$ is an oriented cycle. 
We define $\Gamma_\tau(\xi)$ to be the class 
$[\varphi_\xi^{-1}(p)] \in H_1(M;\zed)$, 
which is independent of $p$ and depends on $\xi$ and $\tau$ only 
through their homotopy classes. 
Clearly, as a 2-plane bundle over $M$, $\xi\cong \varphi_\xi^* (TS^2)$, 
so $c_1(\xi)$ is Poincar\'e dual to $2\Gamma_\tau (\xi)$. 

Now for fixed $x\in H_1(M;\zed)$, the set $\Gamma_\tau^{-1}(x)$ 
is identified with 
the set of framed cobordism classes of framed links representing $x$. 
Clearly, $\Gamma_\tau^{-1}(x)$ is nonempty and has 
a canonical $\zed$-action, where 
$n\in\zed$ acts by adding $n$ right twists to the framing. 
The $\zed$-action on $\Gamma_\tau^{-1}(x)$ is obviously transitive. 
To verify that the stabilizer of a class is $d(\xi)\zed$, fix a nonempty 
framed link $L$ representing $x$, and let $L'$ denote $L$ with $n$ 
twists added to its framing. 
Suppose there is a framed cobordism in $I\times M$ between $L$ and $L'$,  
with $L\subset 1\times M$ and $L'\subset 0\times M$. 
By gluing $1\times M$ to $0\times M$, we get a closed surface in 
$S^1\times M$ with self-intersection number $n$. 
Let $\alpha \in H_2(S^1\times M;\zed)$ denote the corresponding homology 
class with $\alpha^2=n$, and let $\lambda = [S^1]\times x \in H_2 (S^1\times 
M;\zed)$. 
Then $\alpha-\lambda$ intersects $0\times M$ trivially, so it pulls back to 
$H_2(M)$. 
Thus, $(\alpha-\lambda)^2=\lambda^2=0$, and $n= ((\alpha-\lambda)+\lambda)^2 
= 2(\alpha-\lambda)\cdot\lambda$ in $S^1\times M$. 
But this equals $2(\alpha-\lambda)\cdot x = \langle \alpha-\lambda, c_1(\xi) 
\rangle$ in $M$. 
Hence, $n$ is divisible by $d(\xi)$. 
Conversely, we can find a class $\beta\in H_2(M)$ with $\langle \beta,c_1 
(\xi)\rangle = d(\xi)$, and construct a framed cobordism as above with 
$\alpha-\lambda = \beta$ and $n=d(\xi)$. 
We conclude that the stabilizer of the framed cobordism 
class of $L$ is $d(\xi)\zed$, completing the proof.\qed 
\enddemo 

\demo{Remark} 
This proposition can be interpreted in terms of $\spinc$-structures. 
A 2-plane field $\xi$ as above determines a complex structure on 
$TM\oplus\real$ (by splitting it into a pair of complex line bundles), 
and hence, a $\spinc$-structure on $M$ (since $\Spinc (3) = U(2)$ acting 
on $\Lambda^- (\real^4) \cong \real^3$). 
A trivialization $\tau$ identifies the set of $\spinc$-structures with 
$H^2(M;\zed)$, and the class assigned to the $\spinc$-structure
given by $\xi$ is Poincar\'e dual to 
$\Gamma_\tau(\xi)$. 
To see this, note that since $\xi \cong \varphi_\xi^* (TS^2)$, 
this procedure for mapping $\xi$ into $H^2(M;\zed)$ is the pull-back under  
$\varphi_\xi$ of the corresponding procedure on $TS^2\oplus\complex$, 
which has $c_1=2[S^2]$ and hence determines the unique $\spinc$-structure 
on $T\real^3|S^2$ corresponding  to $[S^2] \in H^2(S^2;\zed)$. 
An immediate corollary is that the canonical map from 2-plane fields on $M$ 
to $\spinc$-structures is surjective, and the preimage of each 
$\spinc$-structure is a $\zed$-space isomorphic to $\zed/d$, where $d$ is the 
divisibility of the Chern class of the structure. 
The grading of the 
Seiberg-Witten-Floer theory of the $\spinc$-structure is also $\zed/d$. 
Our invariant $\widetilde\Theta$ provides a connection between these 
$\zed$-spaces, since it is defined via the obstruction $\frac14 (c_1^2 
-2\chix -3\sigma)$ to extending an almost-complex structure over a closed 
4-manifold $X$ from $X-p$. 
The same quantity gives the Seiberg-Witten index of the corresponding 
$\spinc$-structure on $X$. 
Recently, Kronheimer and Mrowka \cite{KM} have defined an explicit 
bijection between homotopy classes of 2-plane fields and Seiberg-Witten-Floer 
groups. 
\enddemo 

We wish to define the invariant picking out the 3-dimensional obstruction 
(from $\zed/d(\xi)$ up to translation) by making $(M,\xi)$ bound an 
almost-complex 4-manifold $X$, defining the obstruction to be $c_1^2(X) 
- 2\chix (X)-3\sigma (X)$ (where $\chix$ and $\sigma$ denote the topological 
Euler characteristic and signature), and proving invariance by observing 
that this quantity vanishes for closed, almost-complex 4-manifolds. 
However, we need additional structure to define the first term, since there 
is no natural quadratic form on $H^2(X) \cong H_2(X,\partial X)$ in general. 

The proof of Proposition 4.1 actually shows that it makes sense to talk 
about a framing on a homology class $x\in H_1(M;\zed)$ in an oriented 
3-manifold $M$ (not necessarily connected), modulo twice the 
divisibility of $x$, by picking a framed cobordism class mapping to $x$. 
Equivalently, one frames an oriented link representing $x$, then observes 
that any other such nonempty representative of $x$ is uniquely framed (modulo 
twice the divisibility of $x$ and allowing twists to transfer between link 
components in each component of $M$) via a framed cobordism. 
Now consider a {\it smooth 1-cycle\/} $\gamma$ carried by an  oriented 
link $L= \coprod \gamma_i$ in $M$, which we define to be an integer linear 
combination $\sum k_i\gamma_i$ of components of $L$. 
Any framing on $L$ determines a framing on $[\gamma]\in H_1(M;\zed)$ (but 
not conversely, in general) by replacing each $\gamma_i$ with $|k_i|$ 
parallel copies of $\gamma_i$, determined by the given framing on $\gamma_i$. 
Now if $z$ is a relative rational 2-cycle in an oriented 4-manifold pair  
$(X^4,\partial X^4)$, with $\partial z$ a smooth 1-cycle 
in $\partial X$, we can define the square of $z$, $Q_f(z)\in \que$, relative 
to a framing $f$ on the link $L$ carrying $\partial z$, by adding 2-handles 
to $X$ along $L$ with framings given by $f$. 
Then $z$ will extend canonically to a rational cycle $\hat z$ in $\hatX= 
X\cup \text{2-handles}$, and we can define $Q_f(z)$ to be $\hat z^2$ 
in the intersection pairing on $H_2(\hatX;\que)$. 
This quantity does not change if we break up each $k_i\gamma_i$ into 
$|k_i|$ parallel copies of $\gamma_i$ using the framing $f$,  
since the 2-handles attached to the 
parallel copies will lie inside the 2-handle attached to $\gamma_i$. 
Similarly, $Q_f(z)$ is preserved if we change $z$ and $(L,f)$ 
by attaching a framed 
cobordism in $I\times \partial X$ to $z$ in $X$. 
If $z$ is an integer cycle and we change it within its class in $H_2(X, 
\partial X;\zed)$ keeping $\partial z$ fixed, then $[\hat z]$ changes  
by a class in $H_2(\partial X;\zed)$, so $Q_f(z)\in\zed$ reduces 
modulo twice the divisibility of $[\partial z]$ to a class $Q_f[z]$ 
depending only on $[z]\in H_2(X,\partial X;\zed)$ and $f$ on $[\partial z]$. 
Similarly, if $z$ is rational and $[\partial z]$ vanishes in $H_1(\partial X;
\que)$ then $Q_f(z)\in\que $ is determined by 
$[z]\in H_2(X,\partial X;\que)$ and $f$ on 
$[\partial z]\in H_1(\partial X;\zed)$. 

Now if $x\in H_1(M;\zed)$ is rationally trivial, then we can assign numbers 
to framings $f$ on $x$ just as we do for framings on nullhomologous knots. 
Simply define $q_f(x)\in\que$ to be $Q_f(z)$ for any rational cycle $z$ in 
$(I\times M,\{1\}\times M)$ with $[\partial z]=x$. 
As above, this only depends on $x$ and $f$. 
Clearly, if $\gamma$ is a $\zed$-nullhomologous knot, then $q_f[\gamma]$ 
is the usual integer assigned to the framing $f$ on $\gamma$. 
In general, $q_f(kx) = k^2 q_f(x)$. 
Adding a right twist to $f$ (on a link component with multiplicity~1) adds 1 
to $q_f(x)$ (and it adds $k_i^2 $ if the multiplicity is $k_i$), so  
$q_f(x)$ mod~1 is independent of $f$. 
In fact, it equals the square of $x$ under the {\it linking pairing\/}, 
the well-known $\que/\zed$-valued symmetric bilinear form on the torsion 
subgroup of $H_1(M)$. 
It follows immediately that $q_f(x)\in\frac1k \zed$, where $k$ is the 
order of $x$. 

For some classes in $H_2(X^4,\partial X^4)$, there is a canonical square,  
independent of a choice of framing. 
By the long exact sequence for $(X,\partial X)$ and the fact that the image 
of $H_2(\partial X)$ in $H_2(X)$ is annihilated by the intersection pairing, 
we see that $\ker (\partial_*:H_2(X,\partial X;\que)\to H_1(\partial X;\que))$ 
inherits a pairing. 
Thus, if $z$ is a rational 2-cycle in 
$(X,\partial X)$ whose boundary is a 
rationally nullhomologous smooth 1-cycle 
in $\partial X$, then $z^2\in\que$ is canonically 
defined. 
To compare $z^2$ with $Q_f(z)$ for any framing $f$ on a link 
carrying $\partial z$, we form 
$\hat z$ in $\hatX$ as above. 
If $z_1$ is a rational 2-chain in $\partial X$ with $\partial z_1=\partial z$, 
then $\hat z= (z-z_1) +\hat z_1$. 
Squaring, we obtain $\hat z^2 = (z-z_1)^2 + \hat z_1^2$ (the cross-term 
vanishes since $z-z_1$ can be pushed into $\rint X$), or $Q_f(z) = z^2 + 
q_f[\partial z]$. 
In particular, the right-hand side is  integral if $z$ is. 

We can now define invariants $\Theta_f$ and $\theta$ that partially capture 
the 3-dimensional uniqueness obstruction of plane fields. 
We denote any Poincar\'e duality isomorphism by $PD$. 

\definition{Definition 4.2} 
Let $\xi$ be an oriented 2-plane field on a closed,  oriented 
3-manifold $M$ (not necessarily connected). 
We say that $(M,\xi)$ is the {\it almost-complex boundary\/} of a 
compact, almost-complex 4-manifold $X$ if $\partial X=M$ 
(as an oriented manifold) and $\xi$ is the field of complex lines in $TM 
\subset TX|M$. 
If so, then for any framing $f$ on $PDc_1(\xi)\in H_1(M;\zed)$, 
let $\Theta_f(\xi) = Q_f(PD\,c_1(X))-2\chix (X) -3\sigma (X) \in\zed/2d(\xi)$. 
If $c_1(\xi)$ is a torsion class, let  
$\theta (\xi) = (PD\, c_1 (X))^2 -2\chix (X) - 3\sigma (X)\in\que$. 
\enddefinition

\proclaim{Proposition 4.3} 
When $\theta (\xi)$ is defined, we have 
$\Theta_f (\xi)=\theta(\xi)+q_f (PDc_1(\xi))\in\zed$. 
In particular, $\theta(\xi)$ is congruent {\rm mod~1} 
to $-q_f(PDc_1(\xi))$, which is an integer divided by the 
(finite) order of $c_1(\xi)$.\qed 
\endproclaim 

For any $(M,\xi)$ as in Definition 4.2, the required $X$ can be constructed 
by applying the following lemma to each component of $M$. 
The invariants $\Theta_f(\xi)$ and $\theta (\xi)$ are independent of $X$ 
by the next theorem. 
Although we do not need $X$ to be spin or a 2-handlebody, we can arrange 
these additional conditions with no extra work. 

\proclaim{Lemma 4.4} 
Let $\xi$ be an oriented 2-plane field on a closed, connected, oriented 
3-manifold $M$.  
Fix a spin structure $s$ on $M$. 
Then $(M,\xi,s)$ is the almost-complex, spin boundary of a compact, 
almost-complex, spin 
4-manifold $X$ consisting of a 0-handle and 2-handles. 
\endproclaim 

\demo{Proof} 
By \cite{Ka} (for example), we can realize 
$(M,s)$ as the spin boundary of a spin handlebody $X'$ with only 
0- and 2-handles. 
Define a complex structure on $TX'|M$ by declaring $\xi$ and a complementary 
trivial bundle to be complex line bundles. 
Since $SO(4)/U(2)\approx S^2$ is simply connected, we can 
extend the almost-complex 
structure over the cocores of the 2-handles, and then 
over the complement of the center $p$ of the 0-handle. 
Now $c_1(X') \in H^2(X';\zed) \cong H^2 (X'-\{p\};\zed)$ is defined. 
According to \cite{HH}, the 4-dimensional obstruction to defining an 
almost-complex structure on a closed, oriented 4-manifold $W$ is given by 
$\frac14 (c_1^2 (W)-2\chix (W) - 3\sigma (W))$, and this can be written 
as a sum of local obstructions at isolated singularities. 
But $S^2\times S^2$ admits a complex structure with $c_1^2=8$, and an almost 
complex structure in the complement of a point with $c_1=0$. 
Thus, we can add $\pm1$ to the index of a singularity by forming the 
connected sum with $S^2\times S^2$ there. 
Now for a suitable connected sum $X=X'\# n S^2\times S^2$, the obstruction at 
$p$ will vanish, yielding our required almost-complex spin manifold.\qed 
\enddemo 

\proclaim{Theorem 4.5} 
Let $\xi$ be an oriented 2-plane field on a closed, oriented 3-manifold 
$M$ (not necessarily connected), and let $f$ be a framing on $PDc_1(\xi)$. 
Then $\Theta_f(\xi)\in \zed/2d(\xi)$ depends only on $M,f$ and the homotopy 
class  $[\xi]$ of $\xi$. 
If $c_1(\xi)$ is a torsion class, then $\theta (\xi)\in\que$ depends only 
on $M$ and $[\xi]$. 
Both invariants are independent of the orientation of $\xi$, reverse sign 
if the orientation of $M$ is reversed, and add under disjoint union. 
Adding a right twist to $f$ increases $\Theta_f(\xi)$ by~1. 
If $\varphi :M\to M'$ is an orientation-preserving diffeomorphism, then 
$\Theta_{\varphi_*f} (\varphi_*\xi) = \Theta_f(\xi)$ and $\theta (\varphi_*\xi)
=\theta (\xi)$ (when defined).
\endproclaim 

\demo{Proof} 
Given $(M,\xi_0)$ as in the theorem and $\xi_1$ homotopic to $\xi_0$, 
let $X_0$ and $X_1$ be any compact, almost-complex manifolds 
with almost-complex boundaries $(M,\xi_0)$ and 
$(\oM,\xi_1)$, respectively, where $\oM$ 
denotes $M$ with reversed orientation. 
Glue $X_0$ and $X_1$ to opposite boundary components of $I\times M$ to 
obtain a closed, oriented manifold $W$,  
and extend the almost-complex structures across $I\times M$ 
to make $W$ almost-complex, 
using a homotopy from $\xi_0$ to $\xi_1$. 
Clearly, $\chix(W) = \chix (X_0) + \chix (X_1)$ and similarly for $\sigma$. 
Also, $c_1^2 (W) - 2\chix (W) - 3\sigma (W)=0$. 
To analyze $\theta$, we assume that $c_1(\xi_0)$ 
is a torsion class and let $\theta_i = \theta (\xi_i)$ be as in 
Definition~4.2, defined using the manifold $X_i$. 
Now rationally, $PDc_1 (W)$ can be written  as $\alpha_0+\alpha_1$, where 
$\alpha_i\in H_2(W;\que)$ pulls back to $H_2(X_i;\que)$ 
and then maps to $PD c_1(X_i)$ in $H_2(X_i,\partial X_i;\que)$. 
It follows that $c_1^2 (W) = (PD\,c_1 (X_0))^2 + (PD\,c_1(X_1))^2$, so the 
above formula on $W$ implies that $\theta_0+\theta_1=0$. 
Since $X_0$ and $X_1$ were chosen independently, $\theta_0$ depends only on 
$M$ and the homotopy class of $\xi_0$, and $\theta_1= -\theta_0$ is the 
result of reversing the orientation of $M$. 
Reversing the orientation of $\xi_i$ corresponds to conjugating the 
almost-complex structure on $X_i$, which fixes $(PD\,c_1(X_i))^2$. 
A similar argument applies to $\Theta_f$ in the general case, 
once we observe that for any framing $f$ on $PD(c_1(\xi))$, $c_1^2(W)$ 
reduces $\mod 2d(\xi)$ to $Q_f(PD\,c_1(X_0)) + Q_f(PD\,c_1(X_1))$. 
The last sentence of the theorem follows immediately from the observation 
that if $(M,\xi)$ is the almost-complex boundary of $X$, then $(M',\varphi_*
\xi)$ is the almost-complex boundary of $X$ by the gluing map $\varphi$.\qed 
\enddemo 

\proclaim{Corollary 4.6} 
For $i=1,2$, let $(M_i,\xi_i)$ be the holomorphically fillable contact 
3-manifold obtained by contact surgery on a Legendrian knot $K_i$ in $S^3$. 
Suppose there is an orientation-preserving diffeomorphism $\varphi:M_1\to 
M_2$ such that $\varphi_* (\xi_1)$ is homotopic to $\xi_2$ (with 
either orientation). 
Then $tb(K_1) = tb(K_2)$ and $|r(K_1)| = |r(K_2)|$. 
\endproclaim 

For example, if $K$ is any Legendrian knot in $S^3$ with $tb(K)\ge2$ 
(respectively, 4), then by adding zig-zags to $K$ we can obtain smoothly 
equivalent Legendrian knots $K_i$ with $tb(K_i)=0$ 
(respectively 2) but distinct values of $|r(K_i)|$. 
Contact surgery on these knots will yield noncontactomorphic 
holomorphically fillable 
contact structures on the same homology sphere. 
(In fact, they bound diffeomorphic Stein manifolds with different Chern 
classes.) 
Their homotopy classes will be distinguished by the 3-dimensional 
obstruction $\theta$, but not by the 2-dimensional one 
$\Gamma_\tau\in H_1(M;\zed)=0$ --- in particular, not by the Chern class, 
providing new counterexamples to Conjecture~10.3 of \cite{E3}. 
Similarly, one can realize such examples for any finite cyclic $H_1$ by 
arranging $|r|$ to take more values than the order of $H_1$. 

\demo{Proof} 
Let $r_i$ denote $r(K_i)$ and let $n_i=tb(K_i)-1$ denote the surgery 
coefficient. 
Then $|n_1| = |n_2|$ is the order of the cyclic group $H_1(M_i;\zed)$. 
If $n_i=0$, then this group is infinite, $|r_i|= d(\xi_i)$, and we are done. 
Otherwise, the group is finite, $\theta (\xi_i)\in \que$ is defined, and 
$\theta (\xi_1)= \theta (\xi_2)$. 
We compute $\theta (\xi_i)$ using the Stein surface $X_i$ obtained by 
adding a 2-handle to $B^4$ along $K_i$. 
Let $\alpha_i$ be a generator of $H_2(X_i;\zed)\cong \zed$. 
Then $\alpha_i^2=n_i$, so the cocore disk $D_i$ of the 2-handle represents 
${1\over n_i} \alpha_i \in H_2(X_i,\partial X_i;\que) \cong\que$, and 
$PD c_1(X_i) = {r_i\over n_i} \alpha_i$ over $\que$. 
Thus, $\theta (\xi_i) = {r_i^2\over n_i} - 4-3\text{ sign }n_i$. 
Now it clearly suffices to show that $n_1=n_2$. 
If not, then $n_1=-n_2$, and setting $\theta (\xi_1) = \theta (\xi_2)$ shows 
that $r_1^2 + r_2^2 = 6|n_1|$. 
Let $x_i\in H_1 (M_i;\zed)$ be the generator $[\partial D_i]$. 
Then $q_f(x_i)\equiv -{1\over n_i}$ (mod~1). 
(For example, let $z= D_i - {1\over n_i} \alpha_i$, which can be pushed into 
$I\times M_i$.) 
But for some $k\in \zed$ we have $\varphi_* (x_1)= kx_2$, so $q_{f_1} (x_1) 
= q_{\varphi_*f_1} (\varphi_*(x_1)) \equiv k^2 q_{f_2} (x_2)$ (mod~1), 
and therefore $k^2\equiv -1$ mod~$|n_1|$. 
Reducing the equation $r_1^2 + r_2^2 = 6|n_1|$ modulo~3, we see that $r_1$ 
and $r_2$ must be divisible by 3, and so $n_1$ is also divisible by 3. 
But the equation $k^2\equiv -1$ (mod~3) has no solutions, so we have the 
required contradiction.\qed
\enddemo 

Like $\theta$ in the case $d(\xi)=0$, the invariant $\Theta_f$ can easily 
be computed for the boundary of any Stein surface $X$ in standard form. 
In fact, $PDc_1(X)$ is represented by $\sum r(K_i)D_i$, where $D_i$ is 
the cocore of the handle of $X$ attached to $K_i$. 
Thus, $Q_f(PDc_1(X))=0$ when $f$ is the 0-framing on the union of 
meridians carrying $\sum r(K_i)\partial D_i$. 
To compare $\Theta$ for two different Stein surfaces with a diffeomorphism 
preserving $c_1(\xi)$ between the boundaries, it now suffices to compare 
the respective 0-framings by constructing a framed cobordism between 
the corresponding 1-cycles. 

We will see that the canonical generator of the $\zed$-action of 
Proposition~4.1 subtracts 4 from each of $\Theta_f$ and $\theta$. 
Thus, when $c_1(\xi)$ is a torsion class, either invariant captures the 
3-dimensional obstruction. 
When $c_1(\xi)$ has infinite order, however, $\theta (\xi)$ is undefined 
and $\Theta_f(\xi)\in\zed/2d(\xi)$ only captures the obstruction up to a 
$2:1$ ambiguity. 
We will lift $\Theta_f(\xi)$ to an invariant $\tilde\Theta (\xi,s,f)\in 
\zed/4d(\xi)$ depending on a choice of spin structure $s$, and this will 
capture the 3-dimensional obstruction. 
This lifting is intimately related to the 2-dimensional obstruction 
$\Gamma_\tau (\xi)\in H_1(M;\zed)\cong H^2(M;\zed)$. 
To clarify the relationship between $\Gamma$ and $\tilde\Theta$, we now 
give an alternate definition of $\Gamma$. 
The new definition is readily computable (e.g., Theorem~4.12) 
, and it clarifies the dependence 
of $\Gamma$ on the trivialization $\tau$. 
Note that $\Gamma_\tau (\xi)$ should only depend on the restriction 
of $\tau$ to the 2-skeleton of $M$, or equivalently, on a spin structure. 
Here and elsewhere, it is convenient to use Milnor's definition \cite{M1} 
(see also \cite{GS}) of a spin structure on an oriented manifold $M$ as a 
homotopy class of positively oriented trivializations of $TM$ over the 
2-skeleton of some cell decomposition of $M$, where if $\dim M\le 2$ we 
first stabilize $TM$ by summing with a trivial bundle. 
Note that there is a canonical way to reverse the orientation of a 
trivialization, allowing us to identify spin structures on $\oM$ with those 
on $M$. 

\definition{Definition 4.7} 
Let $\xi$ be an oriented 2-plane field on a closed, oriented 3-manifold $M$ 
(not necessarily connected). 
Let $v$ be a vector field in $\xi$ whose zero locus, counted with 
multiplicities, has the form $2\gamma$ for some smooth 1-cycle $\gamma$ 
carried by a link $L$ in $M$. 
The vector field $v$ in $\xi$ 
determines a spin structure $s$ on $M-L$, and this 
extends uniquely over $M$ since $v$ vanishes with even multiplicity on $L$. 
Define $\Gamma (\xi,s)$ to be the class $[\gamma]\in H_1(M;\zed)$.
\enddefinition 

The next proposition shows that $\Gamma (\xi,s)$ is well-defined and 
depends in a simple way on $s$. 
To understand this dependence, recall that the group $H^1(M;\zed_2)\cong 
H_2(M;\zed_2)$ acts freely and transitively on the set $\S(M)$ of spin 
structures of $M$. 
The same group acts on $H_1(M;\zed)$ via the Bockstein homomorphism $\beta$ 
in homology induced by the coefficient sequence $\zed\to\zed\to\zed_2$. 

\proclaim{Proposition 4.8} 
Let $\xi$ be an oriented 2-plane field on a closed, oriented 3-manifold $M$. 
Then the map $\Gamma(\xi,\cdot):\S(M) \to H_1(M;\zed)$ is well-defined 
(depending only on $M$ and $[\xi]$) and $H^1(M;\zed_2)$-equivariant. 
For $M$, $\tau$ and $\Gamma_\tau$ as in Proposition~4.1, 
$\Gamma_\tau (\xi)$ equals $\Gamma(\xi,s)$ where $s$ is the spin structure 
induced by $\tau$. 
\endproclaim 

\demo{Proof} 
For $i=0,1$, let $v_i$, $L_i$, $\gamma_i$ and $s_i$ be as in Definition~4.7. 
Let $\Delta (v_0,v_1) \in H^1(M-(L_0\cup L_1);\zed)$ be the difference 
class of the nonzero vector fields $v_i$ in $\xi|M- (L_0\cup L_1)$. 
Then the mod~2 reduction of $\Delta (v_0,v_1)$ extends uniquely over 
$H^1(M;\zed_2)$ as the difference class $\Delta (s_0,s_1)$ of the spin 
structures. 
Now $\beta PD\Delta (s_0,s_1)= \frac12\partial PD\Delta (v_0,v_1) = 
[\gamma_0] - [\gamma_1]$, so $\Gamma (\xi,\cdot)$ is well-defined 
on spin structures induced by vector fields as 
in Definition~4.7, and it has the required equivariance on them.  

Now any spin structure $s\in \S(M)$ is induced by a trivialization $\tau$ 
of $TM$, since $\pi_2(SO(3))=0$. 
Using $\tau$, define $\varphi_\xi :M\to S^2$ with $\xi \cong\varphi_\xi^* 
(TS^2)$ as in the proof of Proposition~4.1. 
Let $w$ be a vector field on $S^2$ with a unique zero, occurring at a 
regular value $p$ of $\varphi_\xi$. 
Then $v=\varphi_\xi^* (w)$ is a vector field in $\xi$ vanishing with 
multiplicity 2 on $\gamma=\varphi_\xi^{-1}(p)$. 
Applying Definition~4.7, we see that $\Gamma (\xi,s') = [\varphi_\xi^{-1}(p)] 
= \Gamma_\tau (\xi)$, where $s'$ is the spin structure on $M$ induced by $v$. 
But $s'$ is the pull-back under $\varphi_\xi$ of the unique spin structure 
on $S^2$, which comes from the canonical trivialization of $\real^3$. 
Since this trivialization pulls back to $\tau$, we have $s'=s$. 
In particular, any $s$ is induced by a suitable vector field $v$,  
so $\Gamma(\xi,s)$ is defined, and $\Gamma(\xi,s) = \Gamma_\tau (\xi)$.\qed 
\enddemo 

\proclaim{Corollary 4.9} 
Under the above hypotheses, the map $\Gamma(\xi,\cdot) :\S (M) \to 
H_1(M;\zed)$ is determined by its value on any one $s\in \S$,  
hence, by any vector field $v$ as in Definition~4.7. 
Its image is $\{x\in H_1(M;\zed)\mid 2x = PDc_1(\xi)\}$, and  
$\Gamma (\xi,s_0) = \Gamma(\xi,s_1)$ if and only if 
$\Delta (s_0,s_1)$ lifts to $H^1(M;\zed)$. 
Reversing the orientation of either $\xi$ or $M$ reverses the sign of 
$\Gamma (\xi,s)$ for fixed $s$. 
$\Gamma$ is preserved by surgeries on $S^0$ 
(addition of 4-dimensional 1-handles), 
and adds in the obvious way under disjoint union (hence, connected sum). 
If $\varphi : M\to M'$ is an orientation-preserving diffeomorphism, then 
$\Gamma (\varphi_*\xi,\varphi_*s) = \varphi_*\Gamma(\xi,s)$. 
Under finite coverings, $PD\Gamma$ lifts in the obvious way. 
For $M$ connected, with fixed $s$ and $q\in M$, $\Gamma(\cdot,s)$ 
classifies 2-plane fields on $M-\{q\}$ up to homotopy. 
\endproclaim 

\demo{Proof} 
This is immediate from Definition  4.7 and Proposition 4.8. 
Note that by the long exact coefficient sequence, $\Im \beta$ is the 
subgroup of elements of $H_1(M;\zed)$ with order at most~2, and 
$\ker (\beta\circ PD)$ consists of those classes in $H^1(M;\zed_2)$ with 
integer lifts. 
The behavior of $\Gamma$ under orientation reversals is the same 
as that of $PDc_1(\xi)$. 
(Reversing $[M]$ reverses the sign of $PD$.) 
For the last sentence, extend $s$ to a trivialization $\tau$ and note that 
any 2-plane field on $M-\{q\}$ extends over $M$, since the corresponding map 
$S^2\to S^2$ is nullhomologous. 
Clearly, $PD\,\Gamma (\cdot,s) = PD\, \Gamma_\tau(\cdot)$ is the 
2-dimensional uniqueness obstruction for oriented 2-plane fields on $M$, 
so it classifies 2-plane fields on $M-\{q\}$.\qed
\enddemo 

\demo{Remarks} 

1) Kuperberg (\cite{Ku} Section~2.1) constructs an invariant $c$ equivalent 
to $\Gamma (\xi,s)$ directly from homotopy theory, and essentially also 
obtains the above invariance, kernel and image of $\Gamma (\xi,\cdot)$. 

2) A spin structure $s$ on $M$ determines a $\spinc$-structure by the 
inclusion $\Spin (3) \subset \Spinc (3)$. 
From our viewpoint, this $\spinc$-structure  is the one canonically 
associated to a trivialization $\tau$ extending $s$. 
By the remark following Proposition~4.1, it follows that $\Gamma(\xi,s)$ is 
the difference class of the $\spinc$-structures associated to $\xi$ and $s$. 

3) Here is yet another interpretation of $\Gamma (\xi,s)$, suggested 
by D.~Freed. 
The plane field $\xi$ determines an $SO(2)$-subbundle of the tangent 
$SO(3)$-principal bundle of $M$. 
A spin structure $s$ lifts the latter to an $SU(2)$-bundle, in which $\xi$ 
determines a $U(1)$-subbundle $\tilde\xi$. 
Then $\Gamma (\xi,s) = PDc_1 (\tilde\xi)$. 
To verify this, fix a trivialization $\tau$ extending $s$, write 
$\xi\cong \varphi_\xi^* (TS^2)$ as above, then note that $c_1(\tilde\xi) = 
\varphi_\xi^* c_1 (\widetilde{TS}{}^2) = \varphi_\xi^* [S^2]$ is dual 
to $[\varphi_\xi^{-1}(p)] = \Gamma (\xi,s)$, as required. 
\enddemo 

\proclaim{Corollary 4.10} 
Let $\xi$ be an oriented 2-plane field on a closed,  oriented 
3-manifold $M$. 
Then $\xi$ is homotopic to itself with reversed orientation if and only 
if $c_1(\xi)=0$. 
\endproclaim 

This rules out the case where $c_1(\xi)$ has order 2. 

\demo{Proof} 
Given such a homotopy, then for a fixed spin structure $s$ on $M$ we have 
$\Gamma (\xi,s) = -\Gamma (\xi,s)$, 
so $PDc_1(\xi) = 2\Gamma (\xi,s)=0$. 
The converse is obvious by direct construction, trivializing $\xi$.\qed 
\enddemo 

To compute $\Gamma(\xi,s)$ for the almost-complex boundary $(M,\xi)$ of a 
compact Stein surface $X$, it is convenient to express spin structures in 
terms of characteristic sublinks. 
(See, for example, \cite{GS}.) 
We can assume that $X$ is presented in standard form (Definition~2.1). 
Then there is a canonical way to surger out the 1-handles of $X$ to obtain 
a 4-manifold $X^*$ with $\partial X^*=M$, such that $X^*$ is obtained from 
$B^4$ by adding 2-handles along a framed link $L$ in $S^3$. 
This link is obtained by stretching the box in the plane, forming an annulus 
by gluing together its lateral edges in the obvious way (Figure~26). 
Each 2-handle of $X$ becomes a 2-handle of $X^*$ with the same framing, 
and the remaining 2-handles of $X^*$ form a sublink $L_0\subset L$ 
consisting of a 0-framed unknot for each 1-handle $h$. 
(Curves that formerly ran over $h$ will now link the unknot.) 
Clearly, $H_2(X^*;\zed_2)$ corresponds bijectively to the set of all 
sublinks $L'$ of $L$, with each $L'$ mapping to the homology class determined 
by the cores of the corresponding 2-handles. 
Fix an orientation on $L$. 
For components $K,K'$ of $L$, we define the {\it linking number\/} 
$\ell k(K,K')$ to be the usual linking number if $K\ne K'$ and to be the 
framing of $K$ if $K=K'$. 
We use bilinearity to extend the definition to $\ell k(\gamma,\gamma')$ for 
$\gamma,\gamma'$ smooth 1-cycles carried by $L$, and use the same 
formalism mod~2 for sublinks of $L$. 
These pairings correspond to the intersection pairings of $X^*$ over 
$\zed$ and $\zed_2$, respectively. 
\cfig{2.00}{rg-fig26.eps}{26}

\definition{Definition 4.11 \cite{Ka}} 
A sublink $L'$ of $L$ is {\it characteristic\/} if for each component 
$K$ of $L$, the framing of $K$ is congruent modulo~2 to $\ell k(K,L')$.
\enddefinition 

By the Wu formula, a sublink $L'$ is characteristic if and only if it 
corresponds to a class in $H_2(X^*;\zed_2)$ mapping to $PDw_2(X^*)\in 
H_2(X^*,M;\zed_2)$. 
It is not hard to show (for example, \cite{GS}) that $\S(M)$ maps 
bijectively to the set of characteristic sublinks of $L$ by sending a 
spin structure $s$ to the link $L'$ corresponding to $PDw_2(X^*,s)\in H_2
(X^*;\zed_2)$, the dual of the Stiefel-Whitney class of $X^*$ relative to $s$. 
The spin structure $s$ corresponding to $L'$ 
is then characterized by the fact that if we attach 
a 2-handle to $X^*$ along any knot $K$ in $S^3-L$ with framing $n$, the 
structure $s$ will extend over the new handle if and only if $n\equiv 
\ell k(K,L')$ mod~2. 
Note that the spin structure near $K$ that extends over the 2-handle is 
the one {\it not\/} induced by the normal framing and tangent vectors 
to $K$, since the tangent vector field to $\partial D^2\subset \real^2$ 
has odd degree. 
The difference class $\Delta (s_0,s_1)\in H^1(M;\zed_2) \cong H_2(\partial X^*;
\zed_2) \subset H_2(X^*;\zed_2)$ corresponds to the difference of the 
characteristic sublinks for $s_0$ and $s_1$. 
There is a simple procedure (see \cite{GS}, for example) for following spin 
structures through sequences of Kirby moves via their characteristic sublinks.  

\proclaim{Theorem 4.12} 
Let $X$ be a compact Stein surface in standard form, 
with $\partial X=(M,\xi)$,  
and $X^*$, $L=K_1\cup\cdots\cup K_n$ and $L_0$ as above. 
Let $\{\alpha_1,\ldots,\alpha_n\}\subset H_2(X^*;\zed)$ be the basis 
determined by $\{K_1,\ldots,K_n\}$. 
Let $s$ be a spin structure on $M$, represented by a characteristic sublink 
$L'$ of $L$. 
Then $PD\Gamma (\xi,s)$ is the restriction to $M$ of the class $\rho\in H^2 
(X^*;\zed)$ whose value on each $\alpha_i$ is the integer 
$$\langle \rho,\alpha_i\rangle = \frac12 (r(K_i)+\ell k(K_i, L_0+L'))\ .$$ 
\endproclaim 

Here we define $r(K_i)$ to be $0$ if $K_i$ is in $L_0$, and otherwise to be 
the rotation number of the corresponding Legendrian knot in the 
diagram for $X$. 
The theorem still holds if we replace $L_0+L'$ in the formula
by any smooth 1-cycle 
(independent of $i$) carried by $L$ and agreeing with $L_0+L'$ mod~2. 
This is because 
adding $2\gamma$ to $L_0+L'$ changes $\langle \rho,\alpha_i\rangle$ by 
$\ell k(K_i,\gamma) = \alpha_i\cdot\alpha_\gamma$, 
where $\alpha_\gamma\in H_2(X^*;\zed)$  is the class determined by $\gamma$, 
so $\rho$ changes by $PD\alpha_\gamma$ in $H^2(X^*,M;\zed)\cong 
H_2(X^*;\zed)$, leaving $\rho|M$ unchanged. 
Similarly, $2\rho$ restricts to $c_1(\xi)$ 
as required. 
(Furthermore, the equivariance of Proposition~4.8 can be verified directly.) 

\proclaim{Corollary 4.13} 
Let $K$ be a Legendrian knot in $\# S^1\times S^2$, presented in standard 
form as in Definition~2.1. 
Then $tb(K) +r(K)+1$ is congruent modulo~2 to the number of times $K$ 
crosses 1-handles.
\endproclaim 

\demo{Proof} 
Add a 2-handle along $K$ to obtain a Stein surface $X$. 
By Theorem~4.12, $r(K)+\ell k(K,L_0+L')$ must be even for any characteristic 
sublink $L'$ of $L$ as above. 
By Definition~4.11, $\ell k(K,L')$ is congruent mod~2 to the framing of $K$, 
$tb(K)-1$. 
Since $\ell k(K,L_0)$ is congruent to the number of times $K$ crosses 
1-handles, the corollary follows immediately.\qed
\enddemo 

\example{Example 4.14} 
Consider the Stein surface $X_p$ $(p\ge1)$ shown in Figure~27(a). 
This is obtained by adding a 2-handle to $S^1\times D^3$ along a Legendrian 
knot $K$ that runs $2p$ times over the 1-handle. 
Thus, $w(K)=2p-1$. 
There are $2p-2$ left cusps, 
half oriented upward and the other half downward, so $\lambda_+ = 
\lambda_- = \rho_+ = \rho_- = p-1$. 
Thus, the framing of the 2-handle is $tb(K)-1=0$, and $r(K)=0$. 
The 2-handlebody $X_p^*$ (Figure~27(b)) admits  a unique spin structure,  
whose restriction $s$ to $M_p = \partial X_p^*=\partial X_p$ 
is given by the empty characteristic sublink. 
Using Theorem~4.12, it is easy to calculate that $\Gamma (\xi,s)=p\mu$, 
where $\mu$ is the meridian of $K$ in $H_1(M_p;\zed)\cong \zed/2p \oplus 
\zed/2p$ (which is generated by the two meridians in $\partial X_p^*$). 
Now observe that $X_p^*$ admits an involution $\varphi$ that interchanges 
the two 2-handles but preserves the orientation and $s$. 
Since $\varphi_* \Gamma(\xi,s) \ne \Gamma (\xi,s)$, we conclude that the 
holomorphically fillable positive contact structures $\xi$ and $\varphi_*\xi$ 
on $M_p$ are not homotopic as 2-plane fields. 
In particular, they are not isotopic (although they are contactomorphic 
via $\varphi$). 
However, they are not distinguished by the other homotopy invariants, since 
$c_1(\xi)=c_1(\varphi_*\xi)=0$ and $\theta (\xi) = \theta (\varphi_*\xi)=-2$. 
Figure~27(c) shows that $M_p$ is Seifert fibered over $S^2$ with 3 multiple 
fibers (cf.\ Section~5). 
In particular, it has no incompressible tori, hence, no Giroux torsion. 
The manifold $X_1$ is the disk bundle over $\real P^2$ with (twisted) Euler 
number $-2$, so $M_1$ is diffeomorphic to the  projectivization 
of $T^*\real P^2$. 
Equivalently, $M_1$ is $S^3$ in the quaternions modulo left 
multiplication by the order~8 subgroup generated by $\{i,j\}$. 
\endexample 
\cfig{2.00}{rg-fig27.eps}{27}

\demo{Proof of Theorem 4.12} 
For convenience, we identify $M=\partial X$ with $\partial X^*$, and use 
the same symbol to denote $K_i \subset L-L_0$ and the corresponding 
knot in the diagram for $X$. 
We let $\nu L$ denote the tubular neighborhood of $L$ in which the 
2-handles of $X^*$ intersect $S^3$. 
As in the proof of Proposition~2.3, we use the vector field $\partial \over 
\partial x$ to define a complex trivialization of $TX$ over the 0- and 
1-handles. 
The obstruction to extending this over $X$ is the relative Chern class, dual to 
$\sum r(K_i)D_i$, where $i$ ranges over the 2-handles of $X$ and $D_i 
\subset \nu L$ is a cocore disk of the handle of $X$ attached along $K_i$. 
The trivialization determines a spin structure on $X-\coprod D_i$ 
whose restriction  to $M-  \coprod \partial D_i$ we denote by $s_0$. 
Since $s_0$ is determined by $\partial\over\partial x$ and 
${\partial\over\partial z}$ in the box 
containing the diagram for $X$, it is characterized in $\partial X^*$ 
by extending over a 
2-handle attached to a knot $K\subset S^3-\nu L$ if and only if the framing is 
congruent mod~2 to $\ell k(K,L_0)$. 
(Note that a 0-framed meridian $\mu \subset\partial X^*$ 
of a component of $L_0$ corresponds 
to a horizontal curve crossing the corresponding 1-handle in the original 
diagram, with framing induced by $\partial\over\partial x$. 
The structure $s_0$ will not extend over a 2-handle attached to $\mu$ with 
this framing, due to the odd degree of a tangent vector field to 
$\partial D^2$.) 
Now add a collar $I\times M$ to $X$, gluing $0\times M$ to $\partial X$. 
If we put the above spin structure on $X-\coprod D_i$ and the given structure 
$s$ on $1\times M$, then the obstruction to extending to a spin structure 
on $I\times M\cup_{\partial X} X$ will be dual to a $\zed_2$-cycle 
that lifts to an integer chain of 
the form $c=\sum r (K_i)D_i+z_1$, where $z_1$ is 
a chain in $I\times M$.
If we identify $z_1$ with a chain in $M= \partial X^*$, we can assume that 
its intersection with the 2-handles is a subchain $z_2$ that 
is a linear combination of core disks. 
Since the spin structure on $0\times (M-\coprod \partial D_i)$ is given by 
$s_0$,  the chain $z_1$ represents $PD\Delta (s,s_0)\in H_2 (M,\coprod 
\partial D_i;\zed_2)$. 
After inclusion in $H_2(X^*,B^4;\zed_2)\cong H_2(X^*;\zed_2)$, this class 
corresponds to the sublink $L_0+L'$. 
To verify this, recall that $s_0$ and $s$ extend over a 2-handle attached 
to a knot $K$ in $S^3-\nu L$ if and only if the framing is congruent 
mod~2 to $\ell k(K,L_0)$ or $\ell k(K,L')$, respectively. 
Thus, $s_0$ and $s$ agree on $K$ if and only if $\ell k(K,L_0+L')$ 
vanishes mod~2. 
By definition, the obstruction $\Delta (s,s_0)$ then satisfies 
$\langle \Delta (s,s_0),K\rangle = \ell k(L_0+L',K)$ over $\zed_2$ for all $K$. 
We conclude that in $H_2(X^*;\zed_2)$, the class corresponding to 
$L_0+L'$ agrees with $PD\Delta (s,s_0)$ when paired with any element of 
$H_2(B^4,S^3 -\nu L;\zed_2)\cong H_2(X^*,M;\zed_2)$, so the classes are equal. 

To relate this construction to Definition~4.7 of $\Gamma$, we observe 
that $s$ is determined by a vector field $v$ in $\xi$ that we construct 
from $\partial\over\partial x$ by twisting along $z_1$ --- that 
is, we define $v$ by the formula  
$PD\Delta (v,{\partial\over\partial x}) = [z_1] \in H_2(M,\partial z_1;\zed)$. 
Since $c$ reduces mod~2 to a cycle in $I\times M\cup_{\partial X}X$, 
we have $\partial c=2\gamma$ 
for some smooth 1-cycle $\gamma$ in $S^3-L$. 
Thus, $\partial z_1 = 2\gamma-\sum r(K_i)\partial D_i$. 
Since $\partial\over\partial x$ extends to a vector field in $\xi$ on $M$ 
with zero locus $\sum r(K_i)\partial D_i$, $v$ is defined on $M$ with zero 
locus $2\gamma$. 
By Definition~4.7, $\Gamma (\xi,s) = [\gamma] \in H_1(M;\zed)$. 
But $z_1-z_2$ is a chain in $S^3-L$ with  
$\partial (z_1-z_2)=2\gamma-\sum r(K_i)\partial D_i- \partial z_2$, so 
the right-hand side is nullhomologous in $S^3-L$, and in the 
free abelian group $H_1(S^3-L;\zed)$, $[\gamma]$ is given by the uniquely 
defined class $\frac12 (\sum r(K_i)[\partial D_i]+[\partial z_2])$. 
Now we apply the map $h:H_1(S^3-L;\zed) \cong H_2 (B^4,S^3-L;\zed) 
\cong H^2 (B^4,L;\zed) \cong H^2(X^*,\text{2-handles};\zed)\to H^2(X^*;\zed)$ 
induced by $\partial$, $PD$, excision and inclusion, respectively. 
Since the latter group is also free abelian, $h[\gamma]$ is the uniquely 
defined class $\frac12 (\sum r(K_i)h[\partial D_i]+h[\partial z_2])$. 
Setting $\rho_0 = h[\gamma]$, we see from the definition of $h$ that $\rho_0|M 
= PD\Gamma (\xi,s)$. 
Similarly, $h[\partial D_i] = PD[D_i]$ in $X^*$, 
so $\langle h[\partial D_i],\alpha_j \rangle = \delta_{ij}$. 
Now $\partial z_2 = \partial z_3$ for some cycle $z_3$ in $(B^4,S^3-L)$, 
and $h[\partial z_2] = PD[z_3]$ in $X^*$. 
Thus, $\langle h[\partial z_2],\alpha_i\rangle = [z_3]\cdot\alpha_i 
= [z_3-z_2]\cdot \alpha_i$ since $z_2$ lies in $M=\partial X^*$, and so 
$\langle \rho_0,\alpha_i\rangle = \frac12 (r(K_i) + [z_3-z_2]\cdot \alpha_i)$ 
for each $\alpha_i$. 
But $z_3-z_2$ is an integer cycle in $X^*$ whose image in 
$H_2(X^*,B^4;\zed_2)$ agrees with that of $z_1$. 
Thus $[z_3-z_2]$ is an integer lift of the $\zed_2$-class represented by 
$L_0+L'$, and the observation  following the theorem shows that the 
given class $\rho$ satisfies $\rho|M = \rho_0|M= PD\Gamma (\xi,s)$.\qed
\enddemo

We are now ready to define the 2-fold lift $\tilde\Theta (\xi,s,f)$ of 
$\Theta_f(\xi) \in\zed/2d(\xi)$, obtaining a complete set of homotopy 
invariants for oriented 2-plane fields. 

\definition{Definition 4.15} 
Let $\xi$ be an oriented 2-plane field on a closed, oriented 3-manifold $M$ 
(not necessarily connected). 
Let $s$ be a spin structure on $M$, and let $f$ be a framing on 
$\Gamma (\xi,s) \in H_1(M;\zed)$. 
Choose a vector field $v$ as in Definition~4.7 determining $s$, with zero 
locus $2\gamma$ carried by $L$, and the framing $f$ determined by a 
framing (also called $f$) on $L$. 
(For example, we can take all multiplicities of $\gamma$ to be 1 by the 
proof of Proposition~4.8.) 
Let $X$ be any compact, almost-complex 4-manifold with almost-complex 
boundary $(M,\xi)$. 
Then $v$ determines a complex trivialization of $TX|M-L$. 
Let $z$ be a 2-cycle in $(X,L)$ that is dual to the relative Chern 
class $c_1(X,v) \in H^2(X,M-L;\zed)$.  
Define $\tilde\Theta (\xi,s,f)$ to be $Q_f (z) -2\chix (X) -3\sigma (X) 
\in\zed/4d(\xi)$. 
\enddefinition 

Note that adding a twist to $f$ on $\Gamma (\xi,s)$ corresponds to adding 
4 twists to the induced framing $\tilde f$ on $PDc_1(\xi) = 2\Gamma (\xi,s)$. 

\proclaim{Theorem 4.16} 
Let $M$, $\xi$, $s$ and $f$ on $\Gamma (\xi,s)$ be as in Definition~4.15. 
Then $\tilde\Theta (\xi,s,f)\in \zed/4d(\xi)$ depends only on $M$, $s$, $f$ 
and the homotopy class of $\xi$, and it reduces to $\Theta_{\tilde f}(\xi)$ 
modulo $2d(\xi)$. 
If $\Gamma (\xi,s_0) = \Gamma (\xi,s_1)$ then $\tilde\Theta (\xi,s_1,f_1)
= \tilde\Theta (\xi,s_0,f_0) + 4(f_1-f_0) + 2\langle c_1(\xi),\lambda\rangle$, 
where $\lambda$ is any integer lift of $PD\Delta (s_0,s_1)$. 
In general, $\tilde\Theta (\xi,s_1,f_1) = \tilde\Theta (\xi,s_0,f_0) + 
Q_{f_0,f_1}(z_M)$, where $z_M$ is any integral 2-cycle in $(I\times M,\{0,1\}
\times M)$ such that $\partial z_M$ has the form   
$2(1\times\gamma_1 - 0\times\gamma_0)$ 
with $[\gamma_i] = \Gamma (\xi,s_i)$ and $[z_M|_2] = PD\Delta(s_0,s_1)$. 
For fixed $s$ and $f$, $\tilde\Theta (\xi,s,f)$ is independent of the 
orientation of $\xi$, reverses sign if the orientation of $M$ is reversed, 
and drops by 4 if $\xi$ is changed on a component of $M$ 
by the canonical generator of the 
$\zed$-action of Proposition~4.1 (for any fixed $\tau$). 
Performing a surgery on $S^0$ in $(M,\xi)$ 
(i.e., a 4-dimensional 1-handle addition 
as in Section~2) adds 2 to $\tilde\Theta (\xi,s,f)$, and $\tilde\Theta$ 
adds in the obvious way under disjoint union. 
If $\varphi :M\to M'$ is an orientation-preserving diffeomorphism, then 
$\tilde\Theta (\varphi_* \xi,\varphi_*s,\varphi_*f) = \tilde\Theta (\xi,s,f)$. 

If $\xi_0$ and $\xi_1$ are oriented 2-plane fields on a 
connected $M$, then they are 
homotopic if and only if for some (hence, any) choice of $s$ and $f$, 
$\Gamma (\xi_0,s) = \Gamma (\xi_1,s)$ and $\tilde\Theta (\xi_0,s,f) = 
\tilde\Theta (\xi_1,s,f)$. 
If $c_1(\xi_0)$ is a torsion class, then the same is true with $\theta (\xi_i)$ or $\Theta_f(\xi_i)$ in place of each $\tilde\Theta (\xi_i,s,f)$. 
If $\Gamma (\xi_0,s) = \Gamma (\xi_1,s)$, then $\theta (\xi_1)-\theta(\xi_0)$ 
(if defined) lifts under finite covering maps $\pi$ by multiplying by 
the degree of $\pi$, and similarly for $\Theta$ and $\tilde\Theta$ (for 
fixed $s$ and $f$). 
\endproclaim 

We can compute $\tilde \Theta$ for a Stein boundary $M=\partial X$ in 
standard form, by reconstructing the setup of the proof of Theorem~4.12. 
Using that notation, let $z_1$ be any 
integral 2-chain in $M$ with $\partial z_1 = 2\gamma - 
\sum r(K_i) \partial D_i$ for some smooth 1-cycle $\gamma$. 
There is a unique spin structure $s$ on $M$ determined by the equation 
$PD\Delta (s,s_0) = [z_1] \in H_2(M,\coprod \partial D_i;\zed_2)$, and its 
characteristic sublink is obtained from the sublink representing 
$[z_1]\in H_2 (X^*,B^4;\zed_2)$ by subtracting the sublink $L_0$ coming 
from the 1-handles of $X$. 
Now $\tilde \Theta(\xi,s,f)$ is obtained for any framing $f$ induced from
$\gamma$ on 
$\Gamma (\xi,s)=[\gamma]$, by adding 2-handles to $X$ 
along $(\gamma,f)$, computing the self-intersection number $Q_f(z)$ of 
$z= \sum r(K_i)D_i+z_1$ suitably extended over the new 2-handles, and 
subtracting $2\chix (X) + 3\sigma (X)$. 
The difference terms $2\langle c_1(\xi),\lambda\rangle$ and 
$Q_{f_0,f_1}(z_M)$ in Theorem~4.16 are easily computed once we recall that
for spin structures $s_0$ and $s_1$ on M, the sublink 
$L_\Delta$ corresponding to $PD\Delta (s_0,s_1)$ is the difference of 
the corresponding characteristic sublinks. 
In fact, $2\langle c_1(\xi),\lambda\rangle$ is given by 
$2\sum r(K_i)\lambda_i$, where $\sum \lambda_i \alpha_i\in H_2(X^*;\zed)$ 
lifts $L_\Delta$ and pulls back to a class $\lambda\in H_2(M;\zed)$. 
To compute $Q_{f_0,f_1}(z_M)$, it suffices to assume that $(\gamma_0,f_0)$ 
is carried by a 0-framed knot in $\partial B^4\subset X^*$, then extend 
$z_M$ by a surface in $B^4$ with boundary $\gamma_0$ and add handles along 
$\gamma_1$ as before. 

Here are a few simple applications of the theorem. 

\proclaim{Corollary 4.17} 
If $\xi$ is an oriented 2-plane field with $c_1(\xi)=0$ 
on a connected $M$, and $H_1(M;\zed)$ 
has no 2-torsion, then all orientation-preserving self-diffeomorphisms 
$\varphi$ of $M$ preserve $\xi$ up to homotopy.
\endproclaim 

\demo{Proof} 
Since $2\Gamma (\xi,s) = PDc_1(\xi)=0$, we must have $\Gamma(\xi,s) =0=
\Gamma (\varphi_*\xi,s)$ for any $s$. 
Since $\theta (\varphi_*\xi) =\theta(\xi)$, the result follows.\qed 
\enddemo 

\proclaim{Corollary 4.18} 
Let $M$ be an integral homology 3-sphere 
with an orientation-reversing self-diffeomorphism $\varphi$. 
If $M$ bounds a smooth, orientable rational homology ball $X$ 
(e.g., if $M=S^3$ or $\Sigma\,\#\, \overline{\Sigma}$ 
for any homology sphere $\Sigma$),  then no 2-plane 
field $\xi$ (or nowhere zero vector field) is invariant up to homotopy 
under $\varphi$.
\endproclaim 

\demo{Proof} 
The rational ball $X$ admits an almost-complex structure. 
(For example, the double of $X$ admits 
an almost-complex structure with isolated singularities \cite{HH}.) 
The resulting complex line field $\xi_0$ on $M$ has $\theta (\xi_0)=-2$. 
Since $H_1(M;\zed)=0$, all 2-plane fields on $M$ can be oriented and 
the $\zed$-action of Proposition~4.1 acts transitively 
on the resulting homotopy classes, so any plane field $\xi$ satisfies 
$\theta (\xi)\equiv 2$ (mod~4). 
But if $\xi$ is preserved by $\varphi$ up to homotopy then $\theta (\xi) = 
- \theta (\xi)=0$.\qed 
\enddemo 

\proclaim{Corollary 4.19} 
Let $M$ be a spherical space form with a (positive) 
contact structure $\xi_0$ whose 
universal cover is the tight structure on $S^3$. 
If $\xi_1$ is any other such contact structure on $M$ with $\Gamma (\xi_1,s)=
\Gamma (\xi_0,s)$, then $\xi_1$ is homotopic to $\xi_0$. 
If $M\ne S^3$ then there is no negative contact structure 
$\xi_2$ on $M$ with $\Gamma (\xi_2,s)=\Gamma(\xi_0,s)$ 
and tight universal cover. 
In particular, any negative contact structure on the Poincar\'e homology 
sphere must have an overtwisted universal cover.
\endproclaim 

\demo{Proof} Since the tight positive contact 
structure on $S^3$ is unique, we must have 
$\theta (\tilde\xi_1) = \theta (\tilde\xi_0)$, where  $\tilde\xi_i$ is  
the universal cover of  $\xi_i$. 
Thus $\theta(\xi_1) = \theta (\xi_0)$ and so $\xi_0$ and $\xi_1$ are homotopic. 
If $\xi_2$ exists, then it is in the same $\zed$-orbit as $\xi_0$, so 
$\theta(\tilde\xi_2) - \theta (\tilde\xi_0)$ must be 
an integer divisible by 4 times the degree of the cover. 
But $\theta (\tilde\xi_0) = -2$, so $\theta (\tilde\xi_2) = +2$ and the 
cover has degree~1. 
The Poincar\'e homology sphere  admits such a $\xi_0$ (since 
the complex line field on $\partial B^4$ in $\complex^2$ is invariant 
under the binary icosahedral group in $SU(2)$), so it admits no 
negative $\xi_2$ with a tight universal cover.\qed 
\enddemo 

\example{Example 4.20} 
Let $X$ be the 4-manifold obtained by attaching a 2-handle to $B^4$ along 
an oriented knot $K$ with framing~0. 
The almost-complex structures on $X$ are classified by the integers, by 
setting $J_k$ equal to the unique structure with $c_1(J_k)$ given by 
$2k$ times the canonical generator of $H^2(X;\zed)$. 
We will compute $\tilde\Theta$ for the 2-plane field $\xi_k$ induced 
on $M= \partial X$ by $J_k$. 
(Note that $\theta (\xi_k)$ is not defined unless $k=0$.) 
When $K$ is the unknot, $M= S^1\times S^2$, and there is a diffeomorphism 
$\varphi :M\to M$ fixing the $S^1$-coordinate and acting on the family 
of 2-spheres by the nontrivial element of $\pi_1 (SO(3))$. 
We will show that $\varphi$ preserves the homotopy class of $\xi_k$ 
if and only if $k$ is either odd or $0$. 
\endexample

There are exactly two spin structures on $M$. 
Let $s$ be the one that extends over $X$, and let $s'$ be the other one. 
Since $H_1(M;\zed) \cong \zed$ has no 2-torsion, we must have 
$\Gamma (\xi_k,s) = \Gamma (\xi_k,s') = \frac12 PDc_1(\xi_k) = k[\partial D]$ 
where $D$ is a cocore of the 2-handle. 
The plane field on $\partial B^4$ determined by $J_k$ contains a nowhere 
zero vector field, which extends to a vector field $v$ in $\xi_k$ with 
zero locus $2k\partial D$, and $v$ induces the spin structure $s$ on $M$. 
We can set $z= 2kD$ in Definition~4.15. 
Let $f_n$ denote the $n$-framing on $\partial D$ (as measured in 
$\partial B^4$), which induces a framing on $\gamma=k\partial D$. 
(Note that we can also obtain framings on $\Gamma(\xi_k,s)= k[\partial D]$ 
between $f_n$ and $f_{n+1}$, by splitting $\gamma$ into $k$ disjoint 
circles.) 
Clearly, $Q_{f_n} (z) = 4k^2n$, so $\tilde\Theta (\xi_k,s,f_n) = 4(k^2n-1) 
\in \zed/8k$. 
Since the generator $\lambda$ of $H_2(M;\zed)$ lifts $PD\Delta (s,s')$, 
changing the spin structure adds $4k$ to $\tilde\Theta$, 
$\tilde\Theta (\xi_k,s',f_n) = 4(k^2 n+k-1) \in \zed/8k$. 
Now if $K$ is the unknot, then $\varphi$ interchanges $s$ and $s'$, and it 
preserves $\partial D$, changing $f_{n-1}$ to $f_n$. 
Thus, $\tilde \Theta (\varphi_* \xi_k,s,f_n) = \tilde\Theta (\xi_k,s',f_{n-1})
= 4(k^2 (n-1) + k-1) = \tilde\Theta (\xi_k,s,f_n) + 4k(1-k) \in\zed/8k$. 
Since $\Gamma(\varphi_* \xi_k,s) =\Gamma (\xi_k,s)$, $\varphi _*\xi_k$ is 
homotopic to $\xi_k$ if and only if $4k(1-k)$ vanishes modulo $8k$, i.e., 
$k$ is either odd or $0$. 
As a check, observe that $\varphi\circ \varphi$ always preserves $\xi_k$ 
as it should (being isotopic to the identity). 
Note that for any $k$ and $n$, $\Theta_{f_n} (\xi_k) = -4\in\zed/4k$, 
so the full invariant $\tilde\Theta$ is required to distinguish 
$\varphi_*\xi$ from $\xi$. 

\demo{Proof of Theorem 4.16} 
We proceed as in the proof of Theorem 4.5. 
Given $(M,\xi_0)$ as in Definition~4.15 and  $\xi_1$ homotopic to $\xi_0$, 
choose almost-complex manifolds $X_0,X_1$ and $W$ as before. 
For $i=0,1$, choose $s_i,f_i,v_i,\gamma_i, L_i$ and $z_i$ for 
$(X_i,\partial X_i)$ as in Definition~4.15, and let $\tilde\Theta_i\in\zed/
4d(\xi_0)$ denote the resulting value for $\tilde\Theta(\xi_i,s_i,f_i)$ on 
$\partial X_i=\pm M$. 
Note that $\partial z_i = 2\gamma_i$. 
Let $z_M$ be a cycle in $(I\times M,0\times L_0\cup 1\times L_1)$ 
dual to the relative 
Chern class $c_1(I\times M,v_0,v_1)\in H^2 (I\times M,0\times (M-L_0) \cup 
1\times (M-L_1);\zed)$ for the given almost-complex structure on  
$I\times M\subset W$. 
Then $\partial z_M = -2(0\times\gamma_0 + 1\times\gamma_1)$ 
in $\partial I\times M$ (recall that 
$PD$ reverses sign under orientation reversal on $M$), and the cycle 
$z=z_0 +z_M+z_1$ in $W$ is dual to $c_1(W)$. 
Reducing modulo~2, we have $[z_M|_2] = PDw_2 (I\times M,s_0,s_1)\in 
H_2(I\times M;\zed_2)$, or projecting into $M$,  $[z_M|_2] = PD\Delta 
(s_0,s_1) \in H_2(M;\zed_2)$. 
Now $c_1^2(W) = z^2 = Q_{f_0} (z_0) + Q_{f_0,f_1} (z_M) + Q_{f_1}(z_1)$, 
so $0=c_1^2(W) - 2\chix (W) - 3\sigma (W) = \tilde\Theta_0 + \tilde\Theta_1 + 
Q_{f_0,f_1} (z_M)$. 
In the case 
$\Gamma (\xi_0,s_0)=\Gamma(\xi_1,s_1)$ (relative to the original orientation 
of $M$), $\gamma_0$ and $-\gamma_1$ will be homologous, so we may compute 
$Q_{f_0,f_1}(z_M)$ by writing $z_M$ as $2(I\times\gamma_0)+z'$ (changing 
$z_M$ and $(\gamma_1,f_1)$ by attaching a framed cobordism if necessary) 
where $z'$ is an 
integer cycle in $I\times M$ with $[z'|_2 ] = PD\Delta (s_0,s_1)$. 
Now $Q_{f_0,f_1}(z_M)=\hat z_M^2 = 4(f_1-f_0) + 4\gamma_0 \cdot z' = 
4(f_1-f_0) + 2\langle c_1 (\xi_0),z'\rangle$. 
Setting $f_0=f_1$ and $s_0=s_1$, we see that $[z'] \in H_2(M;\zed)$ 
is even, so that $Q_{f_0,f_1}(z_M)\equiv 0$ 
mod $4d(\xi_0)$, and 
$\tilde\Theta_1 = -\tilde\Theta_0$. 
Thus, $\tilde\Theta (\xi,s,f)$  depends only on $M,s,f$ and 
the homotopy class of $\xi$, and it flips sign when the orientation of $M$ 
is reversed. 
In the general case, we now have $\tilde\Theta (\xi_0,s_1,f_1) = 
-\tilde\Theta_1 =\tilde\Theta (\xi_0,s_0,f_0) + Q_{f_0,f_1}(z_M)$, and the  
formula remains true if we replace $z_M$ by any other relative 
cycle $z'_M$ with 
the same boundary and with $[z'_M|_2] =PD\Delta (s_0,s_1)$ (since 
$z'_M-z_M$ will be an integral cycle in $M$ with $[z'_M-z_M]\in H_2 (M;\zed)$ 
even). 
It is immediate from the definitions that $\tilde\Theta(\xi,s,f)$ reduces 
mod $2d(\xi)$ to $\Theta_{\tilde f}(\xi)$, and the rest of the first paragraph 
of Theorem~4.16 is also clear, except for the assertion about the 
$\zed$-action.

To compute the effect of the standard generator $g$ of the $\zed$-action 
in Proposition~4.1 on a component of $M$ 
(adding a right twist to the framing of $\Gamma_\tau 
(\xi)$), we fix a trivialization $\tau$ of $TM$ 
and a 2-plane field $\xi$, and form $X$ as before. 
Let $\xi'$ be a plane field on $M$ obtained from $\xi$ via $g$. 
We can assume that $\xi=\xi'$ outside of a ball $B$ in $M$, 
that $\varphi_\xi^{-1}
(p)\cap B= \emptyset$ and that $\varphi_{\xi'}^{-1} (p)\cap B$ is an unknot 
with framing~1. 
Now we alter the almost-complex structure on $X$ in a neighborhood of $B$ 
so that it has a unique singularity and its field of complex lines on $M$ 
is $\xi'$. 
We can remove the singularity by a suitable sum with copies of $S^2\times S^2$ 
as in the proof of Lemma~4.4, and use the new manifold $X'$ to compute 
$\tilde\Theta (\xi',s,f)$. 
Clearly, this construction is local, so the change in $\tilde\Theta$ is 
independent of our initial choices, including $M$, $\xi$, $X$ and $\tau$,  
and $\theta$ (when defined) changes by the same number.  
It is obvious from the Hirzebruch-Hopf formula (see the proof of Lemma~4.4) 
that the change must be divisible by 4, but to determine it exactly we 
compute it in an example. 
Let $M=S^3$ and let $\xi$ and $\xi'$ be the unique 
negative and positive tight contact structures on $M$. 
Then $\theta (\xi') - \theta (\xi) =-4$, and we verify that $\xi'$ is 
obtained from $\xi$ via $g$ (up to homotopy). 
Let $S^3$ be the unit sphere in the quaternions $\HH = \{q=z+wj\mid z,w\in
\complex\}$. 
The complex structure is given by left multiplication by $i$, and $S^3 = 
SU(2)$ acts by right multiplication, determining a canonical trivialization 
$\tau$ of $TM$. 
The standard contact structure $\xi'$ is orthogonal to 
the right-invariant vector 
field $v'(q) = iq$, which corresponds to the constant map $\varphi_{\xi'}(q)
=i$ (into $S^2\subset \Im \HH$) under $\tau$. 
The opposite structure $\xi$ is orthogonal to a vector field $v$ obtained 
from $v'$ by conjugating by the orientation-reversing map $q\mapsto \bar q
= \bar z-wj$. 
Thus, $v(q) = -qi$ and $\varphi_\xi (q) = -qi\bar q$. 
Clearly, $\varphi_\xi (q)=\varphi_\xi(q')$ if and only if $q$ and $q'$ 
lie in the same orbit under right multiplication by $S^1=\complex \cap S^3$. 
Since {\it left\/} multiplication by $S^1$ induces the standard Hopf 
fibration, whose fibers have linking number $+1$, the  fibers of $\varphi_\xi$ 
will have linking number $-1$. 
Thus, $g$ sends $\xi$ to $\xi'$, as required. 

We have now verified that $g$ lowers the invariants $\theta\in\que$ and 
$\tilde\Theta\in  \zed/4d(\xi)$ by 4. 
By Proposition~4.1, $g$ generates a $\zed$-action on 
each set $\Gamma_\tau^{-1}(x)$ of a connected $M$, 
and $\Gamma_\tau^{-1}(x)$ is isomorphic to $\zed/d(\xi)$ as a $\zed$-space. 
Thus, $\tilde\Theta$ (or $\theta$ when defined) 
distinguishes any two nonhomotopic plane fields 
with the same value of $\Gamma_\tau = \Gamma(\cdot,s)$, so $\Gamma$ 
and $\tilde\Theta$ (or $\theta$) are a complete 
set of invariants for distinguishing homotopy classes of plane fields.
Furthermore, for a degree $k$ covering $\tM\to M$ and the trivialization of 
$T\tM$ lifting that of $TM$, $g$ will lift to $kg$ on $\tM$, so 
differences in $\tilde\Theta$ or $\theta$ on 
$\Gamma_\tau^{-1}(x)$ will multiply by $k$ under lifting.\qed  
\enddemo 

\subhead 5. Contact 3-manifolds\endsubhead 

We now examine holomorphically fillable  contact structures on 
oriented 3-manifolds in more detail. 
We show that many such structures have finite covers that are 
overtwisted (Proposition~5.1). 
As an example, we exhibit a pair of fillable structures on an oriented 
3-manifold that are are homotopic as plane fields (Example~5.2), and 
distinguish them by observing that under a certain 2-fold cover one becomes 
overtwisted while the other remains tight. 
We then expand our notation to allow rational surgeries on links 
(Proposition~5.3 and surrounding text) and proceed 
to exhibit fillable structures on several families of oriented 3-manifolds. 
We show that ``most'' oriented Seifert fibered spaces admit fillable 
positive contact structures (Theorem~5.4). 
For example, all such spaces that fiber over (possibly nonorientable)  
surfaces other than $S^2$ bound Stein surfaces realizing both boundary  
orientations, as do many Brieskorn spheres, 
and all Seifert fibered spaces bound Stein surfaces after 
possibly reversing orientation (Corollary~5.5). 
Our construction frequently yields contact structures representing more 
than one homotopy class of 2-plane fields on the manifold; we examine 
the case of circle bundles over surfaces in detail (Corollary~5.7). 
For a family of examples most of which are hyperbolic, we realize ``most'' 
rational surgeries on the Borromean rings as oriented boundaries of 
Stein surfaces (Theorem~5.9 and subsequent corollaries).  
For an arbitrary link, Proposition~5.8 discusses the local topology of 
the set of $n$-tuples of rational coefficients for which surgery yields 
a Stein boundary. 

\proclaim{Proposition 5.1} 
Let $(M,\xi)$ be a contact 3-manifold exhibited as the boundary of a Stein 
surface in standard form (Definition~2.1). 
Suppose that the diagram intersects some disk in $\real^2$ in 
a collection of $m\ge1$ parallel strands in the configuration 
shown in Figure~28, and let $\gamma$ be a loop surrounding the 
strands as shown. 
Let $(\tM,\tilde\xi)\to (M,\xi)$ be any locally contactomorphic 
covering map such that some conjugate 
of $\gamma$ in $\pi_1(M)$ is not in the image of $\pi_1(\tM)$. 
Then $\tilde\xi$ is an overtwisted contact structure. 
\cfig{2.00}{rg-fig28.eps}{28}
\endproclaim 

For example, if $L$ is any Legendrian link in $\# n\,S^1\times S^2$, and 
we modify $L$ by adding zig-zags to some component $K$ to decrease $tb(K)$ 
by 2 without changing $r(K)$, then the proposition applies to the manifold 
$(M,\xi)$ obtained by contact surgery on the modified link, with $\gamma$ 
a meridian of $K$. 
Applying this to the unknot in $S^3$, we see that any lens space of the 
form $L(p,1)$, $p\ge4$, admits fillable structures all of whose covers 
are overtwisted. 
Of course, these same lens spaces are pseudoconvex boundaries of negative 
holomorphic disk bundles over $S^2$, and all covers of the resulting 
contact structures are tight, since they are
fillable by holomorphic disk bundles (by 
branched covering). 

\demo{Proof} 
Let $C$ be the Legendrian curve shown in Figure~29(a), with $tb(C)=-2$. 
As a smooth knot in the complement of the given link, 
$C$ is isotopic to the Whitehead double of $\gamma$,  
so it bounds an immersed disk $D$ disjoint from the 
link, as is clearly visible in Figure~29(b). 
The framing induced by $D$ on $C=\partial D$ is the blackboard framing 
in Figure~29(b), corresponding to the writhe $w(C)=-2=tb(C)$. 
Thus, $D$ induces the canonical framing on $C$. 
The $\pi_1$-condition on $\gamma$ guarantees that some lift $\tD$ of $D$ is 
an embedded disk in $\tM$. 
Since $\tD$ still induces the canonical framing on its boundary, it is an 
overtwisted disk in $\tM$, i.e., its boundary is an unknot with 
$tb(\partial \tD) =0$.\qed
\cfig{2.00}{rg-fig29.eps}{29}
\enddemo 

\example{Example 5.2} 
Let $(M,\xi)$ be the boundary of the Stein surface $X$ obtained by attaching 
a 2-handle  to the Legendrian knot $K\subset S^1\times S^2$ shown in 
Figure~30(a). 
Then $tb(K)=1$ and $r(K)=0$, so the framing on $K$ is $0$ and $c_1(X)=0$. 
By changing the 1-handle to a 2-handle (cf.\ Figure~26), we see that 
$M$ is also obtained by $0$-surgery on the symmetric link $L$ 
shown in Figure~30(b). 
Since $H_1(M;\zed) \cong \zed\oplus\zed$ has no 2-torsion and $c_1(X)=0$, 
we have $\Gamma(\xi,s)=0$ for all spin structures $s$, and $\theta (\xi)$ 
is defined and equals $-2$. 
Since the link $L$ is symmetric, there is an orientation-preserving 
diffeomorphism $\varphi :M\to M$ that interchanges the two surgery tori. 
Clearly, the fillable positive contact structure $\varphi_*\xi$ on $M$ has 
the same invariants $\Gamma$ and $\theta$ as $\xi$, so by Theorem~4.16, 
$\xi$ and $\varphi_*\xi$ are homotopic as 2-plane fields. 
To distinguish the two contact structures, we pass to a double cover. 
Since $\pi_1(X)\cong\zed$, $X$ has a unique double cover $\tX$. 
The corresponding contact manifold $(\tM,\tilde \xi)$ covering $(M,\xi)$ 
bounds $\tX$, so it is tight. 
However, the lift of $\varphi_*\xi$ to $\tM$ is contactomorphic to a 
double cover of $(M,\xi)$ with multiplicity~2 along $K$, so 
Proposition~5.1 shows that it is overtwisted. 
Thus, $\xi$ and $\varphi_*\xi$ cannot be isotopic, although they are both 
contactomorphic (via~$\varphi$) and homotopic as 2-plane fields. 
For similar examples that are not homotopic, consider Example~4.14 
with $p\ge2$. 
These examples can also be distinguished by a new technique of 
Akbulut and Matveyev \cite{AM} relying on gauge theory. 
\cfig{2.00}{rg-fig30.eps}{30}
\endexample

For the remainder of the paper, we will work with {\it rational surgery\/} 
diagrams for 3-manifolds. 
For more details, see \cite{GS} or \cite{Ro}. 
Let $M_0$ be an oriented 3-manifold, and suppose we have established a 
convention for determining $0$-framings of knots in $M_0$. 
In our applications, $M_0$ will usually 
be $\#n\, S^1\times S^2$ $(n\ge0)$ exhibited 
in standard form as in Section~2. 
Let $L$ be a link in $M_0$. 
If we assign a rational number $r_i = {p_i\over q_i} \in \que \cup\{\infty\}$ 
($p_i,q_i$ coprime integers) to each component $K_i$ of $L$, then we 
obtain a 3-manifold $M$ by {\it rational surgery\/} on $L$ with coefficients 
$r_i$, as follows. 
Each $K_i$ has a tubular neighborhood $\nu K_i$. 
Let $(\mu_i,\lambda_i)$ be a positively oriented basis for 
$H_1(\partial\nu K_i;\zed)\cong \zed\oplus \zed$, where $\lambda_i$ is 
determined up to sign as the class of a parallel copy of $K_i$ determined 
by the 0-framing, and $\mu_i$ is determined by a suitably oriented meridian 
(nullhomologous circle in $\nu K_i$); see Figure~31. 
We obtain $M$ by cutting each $\nu K_i$ out of $M_0$ and regluing it by a 
diffeomorphism of $\partial \nu K_i$ sending $\mu_i$ to 
$p_i\mu_i +q_i\lambda_i$. 
This procedure determines $M$ up to orientation-preserving diffeomorphism. 
Thus, setting $r_i=\infty$ corresponds to deleting $K_i$ from $L$. 
If $M_0=\partial$(0-handle $\cup$ 1-handles) in standard form 
and all coefficients are integral, then 
$M$ is the boundary of the 4-manifold obtained by attaching handles along $L$ 
using the same coefficients. 
For this $M_0$, any 1-handle can be replaced by a 0-framed 2-handle without 
changing $M$; recall Figure~26. 
Two other useful moves leaving $M$ unchanged are the {\it Rolfsen twist\/} 
and {\it slam-dunk\/}. 
For $m\in\zed$, an $m$-fold 
Rolfsen twist (Figure~32) is obtained by cutting $M_0$ open along a disk 
bounded by an unknotted component $K_i$ of $L$, and regluing it with $m$ 
right $360^\circ$ twists. 
(A negative number denotes $|m|$ left twists.) 
The reciprocal of the coefficient of $K_i$ increases by $m$ as shown, 
and the other coefficients $r_j$ increase by $m(\ell k(K_i,K_j))^2$. 
A slam-dunk (Figure~33) is obtained from a pair of components $K_i,K_j$ of 
$L$ with $K_j$ a meridian of $K_i$ and $r_i\in\zed$, by pushing $K_j$ 
across $\partial\nu K_i$ into the solid torus glued to $\partial\nu K_i$, 
changing the coefficient of $K_i$ to $r_i-\frac1{r_j}$.
\cfig{2.00}{rg-fig31.eps}{31}
\cxfig{6.3}{rg-fig32.eps}{32}
\cfig{1.8}{rg-fig33.eps}{33}

\proclaim{Proposition 5.3} 
Let $L$ be a Legendrian link in standard form in $\#n\,S^1\times S^2$, 
with a rational coefficient $r_i$ assigned to each component $K_i$. 
If $r_i<tb(K_i)$ or $r_i=\infty$ for each $i$, then the manifold $M$ 
obtained by rational surgery on $L$ with these coefficients is the oriented 
boundary of a Stein surface. 
\endproclaim 

\demo{Proof} 
First, we erase each component of $L$ with coefficient $\infty$, which 
does not change $M$. 
For each remaining component $K_i$, the coefficient $r_i\in\que$ has 
a unique continued fraction expansion of the form 
$$r_i = a_0 - \cfrac 1\\ 
a_1-\cfrac 1\\ 
a_2-\cdots - \cfrac 1\\
a_k\endcfrac\quad ,$$
with each $a_j\in\zed$, and $a_j\le -2$ for $j\ne0$. 
In fact, $a_0$ is the greatest integer $\le r_i$, and the remaining 
coefficients can be obtained by successively solving and taking the greatest 
integer part. 
Modify $L$ by adding a chain of $k$ unknots linked to $K_i$ and 
associating integer coefficients $a_0,\ldots,a_k$ as shown in Figure~34. 
This well-known procedure leaves $M$ unchanged, 
as is easily seen by a sequence of slam-dunks. 
Since $a_0\le tb(K_i)-1$ and $a_j\le-2$ for $j\ne0$, we can add zig-zags 
until $tb-1$ for each component equals the given integer. 
Applying this procedure to each component $K_i$ of $L$, we obtain a diagram 
of a Stein surface whose oriented boundary is $M$.\qed
\enddemo 
\cxfig{6.3}{rg-fig34.eps}{34}

Now we can use diagrams in standard form with rational coefficients 
$r_i<tb(K_i)$ to draw Stein surfaces bounded by the corresponding 
3-manifolds obtained by rational surgery. 
To uniquely specify a Stein surface, we can use Figure~34 and an additional 
convention: 
Fix an orientation on $L$, orient the new link components so that 
the linking numbers in Figure~34 are nonnegative, and require all additional 
zig-zags to be oriented upward, decreasing $r(K)$. 
In general, the resulting Stein structure depends on the choice of 
orientation of $L$, and other Stein structures can be obtained by 
allowing zig-zags oriented in both directions. 
In practice, we will usually disregard the nonuniqueness 
when realizing 3-manifolds as Stein boundaries. 
For a simple example of this construction, observe that any lens space 
$L(p,q)$ is obtained  by $-\frac pq$-surgery on the unknot, with 
$-\frac pq <-1$ (except for $S^3$ and $S^1\times S^2$). 
The proposition realizes any $L(p,q)$ as the oriented boundary of a Stein 
surface, usually with a variety of Stein structures. 
Thus, we see the construction more clearly with contact lenses. 

For a deeper example, we consider oriented Seifert fibered spaces. 
These are oriented, connected  3-manifolds $M$ that are foliated by circles, 
so that the quotient is a closed, not necessarily orientable, surface $F$, 
and the fibers can be described locally (over orientable neighborhoods 
in $F$) as the orbits of a circle action. 
The manifold $M$ is obtained from the circle bundle over $F$ with oriented 
total space and vanishing Euler class ($F\times S^1$ when $f$ is orientable)
by rational surgery on $k\ge0$ 
fibers (whose 0-framings are determined by the projection to $F$) 
with coefficients $r_i=\frac{p_i}{q_i}$, 
$p_i$ a positive integer called the {\it multiplicity\/}, $i=1,\ldots,k$. 
The pairs $(p_i,q_i)$, together with data determining $F$, are called the 
(unnormalized) {\it Seifert invariants\/}, and specify $M$ up to 
orientation- and fiber-preserving diffeomorphism. 
We obtain $\oM$ ($=M$ with its opposite orientation) by reversing the 
sign of each $r_i$. 
The {\it Euler number\/} $e(M) = \sum_{i=1}^k -\frac1{r_i} \in\que$ 
is an invariant of the Seifert fibration that reverses sign under 
orientation reversal. 
Neumann and Raymond (\cite{NR}, Corollary~5.3) showed that for 
$F$ orientable and $e(M)<0$, $M$ is the oriented, 
strictly pseudoconvex boundary of a 
compact complex surface (a negative definite plumbing), hence, of a Stein 
surface \cite{Bo}. 
We show that Seifert manifolds are oriented Stein boundaries under much 
weaker hypotheses. 
For $r\in\que$, let $r=\llbrack r\rrbrack +\text{frac}(r)$ denote the 
decomposition of $r$ with $\llbrack r\rrbrack \in\zed$ and $\text{frac}(r)
\in [0,1)$. 
Then $-\llbrack -r\rrbrack$ is the smallest integer $\ge r$, so 
$\llbrack r\rrbrack +\llbrack -r\rrbrack = -1$ for $r\notin \zed$. 
For $F$ orientable, let $e_0(M)$ denote $\sum_{i=1}^k \llbrack 
-\frac1{r_i}\rrbrack\in\zed$, which is also an invariant of 
the Seifert fibration.  
Clearly, $e_0(M) + e_0(\oM)=-k_0$, where $k_0$ is the number of coefficients 
with $\frac1{r_i}\notin \zed$. 
(Note that a circle bundle over an oriented base is given by the single  
coefficient $r=-{1\over e}$, where $e$ is the Euler number, and so $e_0(M)=e$. 
The integer $b$ arising in the normalized form of the invariants must then 
be $e_0(M)$ or $-e_0(\oM)$, depending on one's normalization convention.) 

\proclaim{Theorem 5.4} 
Let $M$ be an oriented Seifert fibered space with possibly nonorientable 
base $F$. 
If $F\approx S^2$, define $r'_i \in [-\infty,-1)$, $i=1,\ldots,k$, by 
$-\frac1{r'_i} = \text{\rm frac}(-\frac1{r_i})$, with  
$r_i$ as above. 
Then $M$ is the oriented boundary of a compact Stein surface if at least 
one of the following conditions is satisfied: 
\roster
\item"a)" $F\not\approx S^2$ 
\item"b)" $F\approx S^2$ and $e_0(M)\ne -1$ 
\item"c)" $F\approx S^2$ and if $k\ge 3$ 
there is a permutation of the indices of 
$r'_1,\ldots,r'_k$ after which $r'_i<n(r'_1,r'_2)$ for all $i\ge3$, where 
the function $n$ is defined below. 
This last condition is automatically satisfied if $k\le2$ or $r'_i<-2$ 
for all $i$ or $r'_i < - \llbrack {1\over (1/r'_1) +1} \rrbrack -1$ 
for all $i\ge 2$. 
\endroster
\endproclaim 

\proclaim{Corollary 5.5} 
{\rm (a)} Any oriented Seifert fibered space $M$ bounds a compact Stein 
surface after possibly reversing orientation. 

{\rm (b)} An oriented Seifert fibered space $M$ 
bounds compact Stein surfaces realizing both orientations, 
except possibly when $F\approx S^2$ and $e_0(M)$ equals $-1$ or $1-k_0$. 

{\rm (c)} The following Brieskorn homology spheres bound compact Stein 
surfaces realizing both orientations, except possibly for the negative 
orientation on $\Sigma (2,3,5)$, the Poincar\'e homology sphere: 
$\Sigma (p_1,p_2,mp_1p_2\pm 1)$, 
$\Sigma (2,p_2,(2m+1)p_2\pm2)$, 
$\Sigma (2,p_2,p_3)$ with $p_2<11$, 
$\Sigma (3,p_2,p_3)$ with $p_2<7$. 
\endproclaim 

Part (c) of the corollary is included as a test case of condition (c) 
of the theorem, since Brieskorn spheres with three multiple fibers always 
fail conditions (a) and (b) for one choice of orientation. 
(See the remark following the proof of the corollary.) 
The author has encountered only one stubborn 
Seifert fibered homology sphere (cf.\ Corollary~4.19). 

\proclaim{Conjecture 5.6} 
Only one orientation on the Poincar\'e homology sphere $\Sigma (2,3,5)$ 
is realized as the boundary of a Stein surface. 
\endproclaim 

We define the function  $n: Q\times Q\to\zed\cup \{\infty\}$ of Theorem~5.4, 
with $Q= [-\infty,-1)\cap (\{-\infty\}\cup \que)$, as follows. 
(Note that $-\infty=\infty$; the signs are chosen for compatibility with 
the ordering.) 
Given $r'_1,r'_2\in Q$, define $s\in (-\infty,-1]$ by $\frac1s = -1-
\frac1{r'_1}$. 
If $s=r'_2$, set $n(r'_1,r'_2)=0$ (meaning that $M$ is a Stein boundary 
regardless of $r'_3,\ldots,r'_k$). 
Otherwise, there is a projective linear 
transformation $A=\pm \left[\smallmatrix 
a&b\\ c&d\endsmallmatrix\right] \in PSL (2,\zed)$ acting on 
$\que\cup \{\infty\}$ (slopes in $\zed\times\zed$) by $Ar=\frac{c+dr}{a+br}$, 
such that $As\in (-1,0]$ and $Ar'_2 \in [-\infty,-1)$. 
For example, there is a unique such $A$ with $As= 0$. 
Let 
$$t = \cases 
0&\text{if $A0 \in [0,\infty]$}\\
\frac1{As}&\text{if $A0\in [-1,0)$}\\
Ar'_2&\text{if $A0\in (-\infty,-1)$}\endcases\quad ,$$ 
so that $t\in [-\infty,0]$, and let $M=\max (|a|,|c|)$, $m=\min (|a|,|c|)$. 
Let 
$$n_A (r'_1,r'_2) = -m(\llbrack t\rrbrack +1)-M\ .$$ 
(Here, $-\llbrack -\infty\rrbrack = +\infty$.) 
We define $n(r'_1,r'_2)$ to be the supremum of $n_A(r'_1,r'_2)$ 
over all such $A$. 
For example, suppose there is an integer $\ell$ with $r'_2<\ell< s$. 
For the largest such $\ell$, we have $\ell+1 = -\llbrack -s\rrbrack 
= -\llbrack {1\over (1/r'_1)+1}\rrbrack \le-1$. 
Then for $A= \left[\smallmatrix 1&0\\ -\ell-1&1\endsmallmatrix\right]$ 
we have $As = - (\ell+1-s) = -\text{frac}(-s)\in (-1,0]$ and 
$Ar'_2 = r'_2 -\ell-1\in [-\infty,-1)$ as required. 
Thus, $n_A (r'_1,r'_2)$ is defined and equals $\ell$. 
We conclude that condition~(c) of the theorem is satisfied whenever 
$r'_i<\ell = -\llbrack {1\over (1/r'_1)+1}\rrbrack -1$ for all $i\ge2$. 
The case $r'_i <-2$ for all $i$ is precisely the case $\ell=-2$, 
so we have verified the last sentence of the theorem. 

\demo{Proof of Corollary 5.5} 
To prove (a), it suffices to assume that $F= S^2$. 
Then $F\times S^1$ is obtained by 0-surgery on an unknot $K$, with 
fibers $\{x_i\}\times S^1$ corresponding to meridians of $K$, so a surgery 
diagram of $M$ is given by Figure~35. 
We can assume that either $k_0=k$ or $k=1$, by 
slam-dunking meridians with $\frac1{r_i}\in\zed$ 
and restoring the resulting integer 
framing on $K$ to 0 by a Rolfsen twist on a remaining meridian. 
But $M$ is a Stein boundary unless $k\ge3$, in which case $e_0(M) + e_0(\oM)
= -k_0\le -3$, so either $M$ or $\oM$ satisfies $e_0\ne-1$. 
\cfig{2.00}{rg-fig35.eps}{35}

Since (b) follows immediately from Theorem 5.4, we proceed to part (c). 
The Brieskorn homology sphere $M= \Sigma (p_1,p_2,p_3)$ is defined for 
any pairwise coprime integers $p_1,p_2,p_3\ge 2$ and is the unique (up to 
orientation) Seifert fibered homology sphere with multiplicities $p_1,p_2,p_3$. 
The base $F$ must be $S^2$, and the denominators of the 
coefficients $r_i= \frac{p_i}{q_i}$ 
can be chosen arbitrarily subject to the constraint that $c=q_1p_2p_3 
+ p_1q_2p_3 + p_1 p_2 q_3$ should equal $\pm1$. 
Clearly, $c= - p_1p_2p_3 e(M)$, so the case $c=1$ corresponds to the 
usual orientation with $e(M)<0$, when $M$ is oriented as the link of the 
singularity of $z_1^{p_1} + z_2^{p_2} + z_3^{p_3} =0$ in $\complex^3$, and 
$c=-1$ corresponds to the opposite orientation. 
(See \cite{NR} for details.) 
Since $\frac1{r_i}$ is never an integer, we have $e_0(M) +e_0(\oM)=-3$. 
For each of the Brieskorn spheres listed in the corollary, we will exhibit 
coefficients $r_1,r_2,r_3$ for which $e_0(M)=-1$. 
Then $e_0(\oM\,) =-2$, 
so $\oM$ is a Stein boundary and it suffices to analyze $M$, 
which we will do via part (c) of the theorem. 

We begin with the family $\Sigma (p_1,p_2,mp_1p_2\pm1)$. 
Since $p_1$ and $p_2$ are coprime, we can find $\ell_1,\ell_2\in\zed$ 
with $\left| \smallmatrix p_1&\ell_1\\ p_2&\ell_2\endsmallmatrix\right|=1$. 
We choose $\ell_1\in (0,p_1)$, and it follows that $\ell_2\in (0,p_2)$. 
Let $r_1 = \frac{p_2}{\ell_2}$, $r_2 = \frac{p_1}{-\ell_1}$ and 
$r_3 = \frac{mp_1p_2\pm1}{-m}$. 
Clearly, the resulting $M$ satisfies $c=\pm1$ and $e_0(M)=-1$, 
so it is $\Sigma (p_1,p_2,mp_1p_2\pm1)$ with the interesting  orientation. 
Now $s= -\frac{p_2}{\ell_2}$ and $r'_2=r_2$, and the matrix 
$A= \left[ \smallmatrix p_1&\ell_1\\ p_2&\ell_2\endsmallmatrix\right] 
\in SL (2,\zed)$ satisfies $As =0$, $Ar'_2 = -\infty$. 
We have $n_A (r'_1,r'_2)= -p_1-p_2 >r_3 = r'_3$ except in the case of 
$\Sigma (2,3,5)$ oriented with positive Euler number ($c=-1$). 

Table 1 lists some other families of Brieskorn spheres by their 
multiplicities $p_1,p_2,p_3$, together with 
coefficients $r_1,r_2,r_3$ realizing them with orientation such that 
$e_0(M)=-1$, and a reason why each family consists of Stein boundaries.  
The signs are to be chosen consistently throughout each example, and we 
require $\ell,m$ and the signs to be chosen to avoid the previous examples 
$\Sigma (p_1,p_2,mp_1p_2\pm1)$. 
All multiplicities $p_i$ 
must be $\ge2$, but note that $m=0$ is sometimes allowed. 
Most of these Brieskorn spheres are Stein boundaries by the last sentence 
of Theorem~5.4; for the rest we list a suitable matrix $A$ with $r'_3< 
n_A(r'_1,r'_2)$. 
It is easily verified that the list contains all of the remaining homology  
spheres in the corollary.\qed 

$$\vbox{\halign{
$#$\hfill\qquad
&$\displaystyle#$\hfill\qquad
&$\displaystyle#$\hfill\cr 
\multispan3\hfill \hbox{\smc Table 1}\hfill\cr
\noalign{\vskip12pt} 
\qquad \underline{\quad p_1,p_2,p_3\quad} 
	&\qquad \underline{\quad r_1,r_2,r_3\quad} 
	&\quad \underline{\quad\hbox{Reason}\quad} \cr
\noalign{\vskip6pt}
2,4\ell\pm1, 2(4\ell\pm1) m+4\ell\mp 1 
&{2\over1} , {4\ell\pm1\over -\ell}, {2(4\ell\pm1)m+4\ell\mp1\over 
-(2\ell\pm1) m-\ell} 
&r'_1=-2,\ r'_2,r'_3<-3\cr 
\noalign{\vskip6pt}
2,4\ell\pm3, 2(4\ell\pm3)m +4\ell\pm1 
	&{2\over1} , {4\ell\pm3\over -\ell\mp 1} , 
	{2(4\ell\pm3)m + 4\ell\pm1 \over -(2\ell \pm 1)m -\ell} 
	&r'_1 =-2,\ r'_2,r'_3 <-3\cr  
\noalign{\vskip6pt}
2,7,14m \pm3 
	&{2\over1}, {14m\pm3\over -5m\mp1} , {7\over-1} 
	& A= \left[\matrix 3&1\\ 2&1\endmatrix \right]\cr 
\noalign{\vskip6pt}
2,9,18m\pm5 
	&\frac21, \frac{18m\pm5}{-7m\mp2}, \frac9{-1}
	&A=\left[\matrix 3&1\\ 2&1\endmatrix \right]\cr 
\noalign{\vskip6pt} 
3,4,12m\pm5 
	&\frac32, \frac4{-1}, \frac{12m\pm5}{-5m\mp2}
	&\text{All }r'_i<-2\cr 
\noalign{\vskip6pt} 
3,5,15m\pm2
	&\frac3{-1}, \frac5{-1}, \frac{15m\pm2}{8m\pm1} 
	&\text{All }r'_i <-2\cr 
\noalign{\vskip6pt} 
3,5,15m\pm4 
	&\frac32, \frac5{-2}, \frac{15m\pm4}{-4m\mp1} 
	&\text{All }r'_i <-2\cr 
\noalign{\vskip6pt} 
3,5,15m\pm7 
&\frac31, \frac5{-1}, \frac{15m\pm7}{-2m\mp1} 
&r'_1=-\frac32,\ r'_2,r'_3 <-4\cr
}}$$
\enddemo 

\remark{Remark}
Brieskorn spheres with three multiple fibers are relatively difficult 
to analyze, in the sense that any such $M$ fails both conditions (a) 
and (b) of Theorem~5.4 for one choice of orientation. 
To see this,  note that in general we are free to change each $\frac1{r_i}$ 
by adding any $m_i\in\zed$, provided that $\sum_{i=1}^k m_i =0$. 
(Simply Rolfsen twist along the meridians in Figure~35.) 
Thus, for $M=\Sigma (p_1,p_2,p_3)$ we may assume that 
$-\frac1{r_1},-\frac1{r_2}\in (0,1)$. 
But the equation $c=\pm1$ shows that the numbers $q_i$, hence $-\frac1{r_i}$, 
cannot all have the same sign. 
Thus $-\frac1{r_3}<0$ and $e_0(M)\le -1$. 
Similarly, $e_0(\oM\,) \le-1$. 
Since $e_0(M) +e_0(\oM\,)=-3$, 
either $e_0(M)$ or $e_0(\oM\,)$ equals $-1$. 
In contrast, the Brieskorn sphere $M= \Sigma (2,3,5,7)$ is given by 
$r_i = \frac2{-1},\frac31,\frac53, \frac7{-3}$, so $e_0(M)=e_0(\oM\,)=-2$ 
and both $M$ and $\oM$ are Stein boundaries by Theorem~5.4(b). 
\endremark 

\demo{Proof of Theorem 5.4} 
First, we draw a link diagram of a general Seifert fibered space $M$. 
Figure~36(a) shows a disk bundle over a torus 
$T^2$ with Euler number $e$, drawn with two 1-handles and a 2-handle.
The dashed 
lines are circles lifted from $\{x\}\times S^1$ and $S^1\times\{x\}$ in $T^2
= S^1\times S^1$, 
and the fiber circles in the boundary of the bundle 
are meridians of the 2-handle. 
(See \cite{GS}, for example.) 
We put the diagram in standard form by moving a 1-handle  as 
indicated in (a), obtaining a Legendrian knot with $tb=0$~(b). 
This realizes the 4-manifold by a Stein surface, provided that $e\le-1$. 
We can do a bit better (c) using the trick of Figure~20, realizing 
$T^2\times D^2$ ($e=0$). 
Now recall that Figure~35 shows an arbitrary Seifert fibered space 
with base $F=S^2$. 
If $F=T^2$, we obtain the corresponding picture (of an oriented, not 
necessarily contact, 3-manifold) by adding meridians to Figure~36 ($e=0$) 
with coefficients $r_1,\ldots, r_k$. 
If $F$ is an orientable surface of genus~$g$, the picture is obtained 
from Figure~35 by summing $K$ with $g$ copies of Figure~37(a) and 
retaining the 0-framing on $K$. 
If $F$ is nonorientable with genus~$g$, we use $g$ copies of Figure~37(b) 
and the 0-framing. 
(In the nonorientable case, the ratio $p_i\over q_i$ of the 
Seifert invariants  no longer equals the surgery coefficient 
$r_i$ when the coefficient on $K$ is $0$, since 
the 0-framing on $K$ no longer corresponds to $e=0$ in the circle bundle 
over $F$, but rather $e=-2g$ \cite{GS}. 
In fact, we have $e(M) = \sum_{i=1}^k - \frac{q_i}{p_i} = -2g +\sum_{i=1}^k 
- \frac1{r_i}$.) 
\cfig{2.00}{rg-fig36.eps}{36}
\cfig{2.00}{rg-fig37.eps}{37}

To begin constructing  the required Stein surfaces, we first observe that 
we can change each $\frac1{r_i}$ by any integer using a Rolfsen twist, 
at the expense of changing the framing of $K$. 
We use this to replace each $r_i$ by the unique $r'_i$ for which $r'_i \in 
[-\infty,-1)$, which is given by $-\frac{1}{r'_i}=\text{frac}(-\frac1{r_i})$. 
The resulting coefficient on $K$ is $e_0(M)=\sum_{i=1}^k \llbrack 
-\frac1{r_i}\rrbrack$. 
(For $F$ nonorientable, we are defining $e_0(M)$ here.) 
Since we can draw $K$ as Legendrian with $tb(K) =2g-1$ and add Legendrian 
meridians with $tb=-1$, Proposition~5.3 realizes $M$ as a Stein boundary,  
provided that $e_0(M)<2g-1$. 
(The same technique also shows that we can realize any disk bundle over $F$ 
with Euler class $e\le -\chix(F)$ by a Stein surface, provided that we 
improve Figure~37(b) using Figure~20.) 

We next realize $M$ as a Stein boundary whenever $e_0(M)\ge1$, completing 
the proof of part~(a). 
We change the framing on $K$ from $e_0(M)$ 
to $1$ by an inverse slam-dunk, adding a new 
meridian with coefficient $r'_{k+1} = \frac1{1-e_0(M)} \in [-\infty,0)$. 
For each $T^2$- or $\real P^2$-summand of $F$, we surger a 1-handle to a 
0-framed 2-handle, replacing Figure~37 in the construction by Figure~38. 
Now $K$ is a 1-framed unknot in $S^3$, and we stack the other curves along 
it as in Figure~39. 
Note that as these curves pass through the disk bounded by $K$, they all lie on 
one or two ribbons that appear in the diagram with a right half-twist. 
Next, we blow down $K$ (Rolfsen twist on it with $m=-1$), replacing it by a 
full left twist as in Figure~40. 
The half-twists on the ribbons become left-handed. 
The framing on each curve $K_i$ drops by $(\ell k(K,K_i))^2$, so the new 
framing on each meridian is $r'_i-1\in [-\infty,-1)$. 
The framings on the other curves are either 0 or $-4$, depending on 
whether $F$ is orientable ((a) or (b) in the figures). 
Proposition~5.3 shows that Figure~40 represents a Stein surface. 
\cfig{2.00}{rg-fig38.eps}{38}
\cfig{3.0}{rg-fig39.eps}{39}
	\midinsert {\epsfysize=3.75 truein  
	\centerline{\epsfbox{rg-fig40.eps}} }
	\vskip-40pt
	\botcaption{Figure 40} 
	\endcaption
	\endinsert
\cfig{1.3}{rg-fig41.eps}{41}

It now remains to deal with the case $F\approx S^2$ (Figure~35), 
$e_0(M)=-1$ or 0. 
If $e_0(M)=0$, we can replace $K$ by a 1-handle as in Figure~41, obtaining 
a Stein surface. 
This completes the proof of part~(b). 

To prove part (c), it now suffices to assume that $e_0(M)=-1$. 
We can assume that $k\ge2$, by adding coefficients $r'_i=-\infty$ if 
necessary. 
After permuting the indices of $r'_1,\ldots,r'_k$ if desired, we slam-dunk 
the meridian with coefficient $r'_1$ as in Figure~42. 
We denote the new coefficient of $K$ by $\frac1s = -1-\frac1{r'_1}$ and 
the meridian labelled $r'_i$ by $K_i$, $i=2,\ldots,k$. 
Let $T$ denote the boundary of a tubular neighborhood of $K$, containing 
$K_3,\ldots,K_k$ and  identified with  $\real^2/\zed^2$ so that $\mu$ 
and $\lambda$ correspond to the positive $x$- and $y$-axes, respectively. 
Now we can identify the complement of a tubular neighborhood of $K\cup K_2$ 
in $S^3$ with $T\times I$. 
Thus, $M$ is the  unique oriented manifold obtained from $T\times I$ 
by rational surgery with coefficients $r'_3,\ldots,r'_k$ on curves $K_i$ with 
slope $0$ in $T= \real^2/\zed^2$, where the $0$-framing on each $K_i$ 
is determined by the 
product structure on $T\times I$, with solid tori $S^1\times D^2$ glued 
to the two boundary components so that $\{x\}\times \partial D^2$ maps to 
curves with slopes $s$ and $r'_2$, respectively. 
We will complete the proof by changing coordinates in $T$ and then drawing 
a suitable link picture. 
\cfig{2.00}{rg-fig42.eps}{42}

First, suppose that $s=r'_2$. 
Then there is a matrix $A\in SL(2,\zed)$ such that $As=Ar'_2=0$ and $A0 = 
-\frac pq$, $p,q>0$. 
Applying the corresponding diffeomorphism to $T$, we obtain a Hopf link 
with coefficients $\infty$ and $0$, 
and the curves $K_3,\ldots, K_k$ are sent to 
parallel left-handed $(p,q)$ torus knots in $T$, linking the 0-framed 
unknot $p$ times. 
We erase the $\infty$-framed unknot and replace the 0-framed unknot by a 
1-handle. 
The image of each $K_i$, $i\ge3$, is now a copy of the knot $K'$ in 
Figure~43, which is Legendrian with $tb(K')=-pq$. 
The framing $f$ determined on $K'$ by the product structure of $T\times I$ 
has coefficient $-pq = tb(K')$ in Figure~43, as we can easily see by 
pushing a copy of $K'$ off of $T$ and computing its linking number with $K'$. 
Thus, $f$ equals the canonical framing on $K'$. 
Since the images of the curves $K_i$ are parallel copies of $K'$ determined 
by $f$, it follows that they fit together as a Legendrian link with each 
$tb= tb(K')$. 
But each surgery coefficient is $r'_i<-1$ relative to the framing $f$, or 
$r'_i -pq < tb(K')-1$ in Figure~43 relative to our standard convention 
for framing coefficients. 
Thus, Proposition~5.3 exhibits $M$ as the boundary of a Stein surface in 
the case $s=r'_2$. 
\cxfig{5.0}{rg-fig43.eps}{43}

As we observed when defining the function $n$ (after Conjecture~5.6), if 
$s\ne r'_2$ then there is a matrix $A= \left[\smallmatrix a&b\\ c&d
\endsmallmatrix \right] \in SL (2,\zed)$ such that $As \in (-1,0]$ and 
$Ar'_2 \in [-\infty,-1)$. 
Given any such $A$, we apply the corresponding diffeomorphism to $T$, 
again obtaining parallel copies of a torus knot $K'$ in $T$ wrapped 
around the Hopf link, with the latter now having coefficients $\frac1{As}$ 
and $Ar'_2$ in $[-\infty,-1)$. 
If $A0 \in [0,\infty]$, then the torus knot is right-handed as in 
Figure~44. 
In this case $K'$ is Legendrian with $tb(K') =|ac|- |a|- |c|$, and the framing 
$f$ induced by $T$ is $|ac| = tb(K') + |a| + |c|$. 
Since $f$ is obtained from the canonical framing by adding right twists, 
we can realize the images of $K_3,\ldots,K_k$ by a Legendrian link with 
each component satisfying $tb= tb(K')$, as before. 
The required surgery coefficients are $r'_i$ relative to $f$, or 
$r'_i +tb (K') + |a| +|c|$ using the standard convention. 
By Proposition~5.3, $M$ bounds a Stein surface in this case, provided  
that for $i=3,\ldots,k$, $r'_i <-|a| - |c| = n_A (r'_1,r'_2)$. 

In the remaining case, $A0 \in (-\infty,0)$ and the $(|a|,|c|)$ torus knot 
$K'$ is left-handed. 
Now the picture depends on which of $|a|,|c|$ is larger, so we set 
$M=\max (|a|,|c|)$, $m=\min (|a|,|c|)$ and $(t,t')$ equal to $(\frac1{As},
Ar'_2)$ if $m= |c|$ and $(Ar'_2,\frac1{As})$ otherwise. 
Figure~45(a) shows the resulting Legendrian link. 
Since $t,t'<-1$, the Hopf link satisfies the hypothesis of Proposition~5.3, 
provided that we have added no more than $-\llbrack t\rrbrack -2\ge0$ 
extra left cusps to it as shown. 
If $M\le -m(\llbrack t\rrbrack +1)$, then we can draw the $(m,M)$ torus 
knot correctly by winding around these extra zig-zags as shown, ending 
as in (b) of the figure. 
If the inequality fails, we must add the extra $M+m(\llbrack t\rrbrack +1)$ 
twists to $K'$ by the less efficient procedure on the right of the diagram, 
ending as in (c). 
In either case, $tb(K')=-Mm- \max (0,M+m(\llbrack t\rrbrack +1))$ and the 
framing $f$ induced by $T$ is $-Mm\ge tb(K')$, so we can again realize 
the images of $K_3,\ldots,K_k$ by a Legendrian link with $tb=tb(K')$ 
for each component. 
The surgery coefficients are $r'_i + tb(K') +\max (0,M+m( \llbrack t\rrbrack 
+1))$, and the 3-manifold bounds a Stein surface provided that $r'_i<-M-m
(\llbrack t\rrbrack +1) = n_A (r'_1,r'_2)$ for $i=3,\ldots,k$. 
This completes the proof, since the last sentence of part~(c) was 
verified after the definition of $n$.\qed
\cfig{2.3}{rg-fig44.eps}{44}
\cfig{3.0}{rg-fig45.eps}{45}
\enddemo 

Note that many Seifert fibered spaces can be realized as boundaries of 
Stein surfaces in many ways. 
We have already observed that Proposition~5.3 typically produces a 
4-manifold that admits Stein structures with a variety of Chern classes. 
In addition, the proof of Theorem~5.4 sometimes provides several methods 
for analyzing a given Seifert manifold --- for example, if the base $F$ 
has genus $g>1$, then the cases of large and small $e_0(M)$ overlap. 
The author has made no systematic study of this nonuniqueness, but here 
is one simple family of examples. 
Recall that for any tight contact structure $\xi$ on a 3-manifold $M$,  
we have $|\langle c_1(\xi),F_0 \rangle| \le \max (-\chi (F_0),0)$ 
for any  closed, connected, oriented surface $F_0$ in $M$. 
We will call an oriented 2-plane field {\it allowable\/} if it satisfies 
this inequality for all such surfaces $F_0$. 

\proclaim{Corollary 5.7} 
Let $M$ be the circle bundle with Euler  number $e$ over a closed, 
connected, oriented 
surface $F$ of genus~$g$, and let $q\in M$. 
Then the oriented 3-manifold 
$M$ admits at least $\llbrack g-{e\over2}\rrbrack$ holomorphically 
fillable contact structures whose homotopy classes remain distinct after 
allowing orientation-preserving self-diffeomorphisms of $M$. 
If $g\ne0$ and $e<g$, then these contact structures (with both orientations 
allowed on the 2-planes) 
realize all homotopy classes of allowable 2-plane fields on $M-\{q\}$.  
For any $g$ and $e$,  holomorphically 
fillable structures on $M$ and $\oM$ together realize all allowable 
homotopy classes on $M-\{q\}$. 
\endproclaim 

\demo{Proof} 
As in the previous proof, Figure 37(a) realizes $M$ as a Stein boundary 
if $e\le 2g-2 = -\chi (F)$, via an oriented Legendrian knot $K$ in 
$\# 2g\, S^1\times S^2$ whose meridian we denote by $\mu$. 
In the case of equality we have $r(K)=0$, and in general we can realize 
any $r(K)$ congruent to $e\mod 2$ such that $|r(K)| \le 2g-2-e$ (by adding 
zig-zags to the case $e=2g-2$). 
By a computation similar to that of Corollary~4.6, the resulting contact 
structure $\xi$ has $PDc_1(\xi) = r(K)\mu$ and (for $e\ne0$) 
$\theta (\xi) = {r(K)^2\over e} -2\chi (F)-3\text{ sign }e$, so for the 
$\llbrack g-{e\over2}\rrbrack$ different values of $|r(K)|$, the structures 
$\xi$ will not become homotopic after any orientation-preserving 
self-diffeomorphisms of $M$. 
To determine when these contact structures realize all allowable homotopy 
classes over $M-\{q\}$, recall from Corollary~4.9 
that 2-plane fields $\xi$ on $M-\{q\}$ are 
classified by the invariant $\Gamma$, and $2\Gamma (\xi,s) = PDc_1(\xi)$. 
If $\xi$ is an allowable plane field, then $\Gamma (\xi,s)$ 
pairs trivially with each of the tori in $M$ projecting to circles in $F$. 
Thus, $\Gamma(\xi,s)$ is a multiple $\ell\mu$ of the fiber $\mu$, and it 
suffices to understand which multiples we can realize. 
For the above contact structures, we define $X^*$ as in Theorem~4.12,  
let $s$ be the spin structure on $M$ 
corresponding to the characteristic sublink 
$K$, and compute that $\Gamma(\xi,s) = \frac12 (r(K)+e)\mu$. 
If $0= e< g$, we can realize $\ell\mu$ this way for any $\ell\in \zed$ 
with $2|\ell| \le -\chi (F)$. 
But in this case, $F$ lifts to a section $\tF$ of $M$, so any allowable 
2-plane field satisfies 
$-\chi (F) \ge |\langle c_1(\xi),\tF\rangle | = 
|2\ell\mu\cdot\tF| = 2|\ell|$. 
Thus, we have finished the case $e=0$. 
For $e\ne0$, $\mu$ has order $|e|$ in $H_1(M;\zed)$. 
Since the above contact structures realize $2g-1-e$ consecutive values of 
$\ell$, we can realize all multiples of $\mu$ provided that $2g-1-e\ge|e|$, 
or $e<g$, $g\ne0$. 

To prove the final assertion, we also consider contact structures on $\oM$. 
If $g\ne0$ then we are done because either $M$ or $\oM$ has $e\le0<g$. 
Similarly, if $g=0$, it suffices to assume $e(M)<0$. 
(The case $e=0$ is easy and $e(\oM\,) = -e(M)$.) 
But we have already realized $\Gamma(\xi,s) = \frac12 (r(K)+e)\mu$ for 
$r(K)\equiv e \mod 2$, $|r(K)|\le |e|-2$. 
Since $\mu$ has order $|e|$, we are only missing the case $\Gamma (\xi,s)=0$. 
But $\oM=|e|$-surgery on the unknot is also realized by 
${|e|\over 1-|e|}$-surgery on the unknot via a Rolfsen twist, or by a chain of 
$|e|-1$ unknots with framing $-2$ 
as in Figure~34, which is Stein with vanishing Chern class. 
Now $s$ is represented by the empty characteristic sublink (because of the 
Rolfsen twist) and so $\Gamma (\xi,s)=0$ as required.\qed 
\enddemo 

For comparison, a closed, orientable 3-manifold with a taut foliation must 
be covered by $\real^3$ (except in the cases of $S^1\times S^2$ and 
$\real P^3 \,\#\,\real P^3)$, so spherical space forms such as lens spaces 
and the Poincar\'e homology sphere do not admit such foliations. 
In addition, many Seifert fibered spaces with universal cover $\real^3$ 
admit no taut foliations (or even essential laminations) \cite{Br}, although 
they do admit fillable contact structures, frequently with both orientations. 
Taut foliations on Seifert fibered spaces are typically 
horizontal, i.e., transverse to the fibers 
after a suitable isotopy (cf.\ \cite{Br}). 
In contrast, many holomorphically fillable contact structures are neither 
horizontal nor vertical (parallel to the fibers after isotopy). 
In fact, there is a unique homotopy class of 2-plane fields transverse to a 
given Seifert fibration, and when $F$ is orientable 
any vertical 2-plane field $\xi$ has $c_1(\xi)=0$ 
(since the fibers determine a section of $\xi$). 
Thus, most of the examples in the above corollary will satisfy neither 
condition. 

Since Seifert fibered spaces form a somewhat exceptional class of 3-manifolds, 
it seems useful to examine a more generic family. 
Every oriented 3-manifold can be realized as integer surgery on a link 
in $S^3$, so 
it seems natural to consider a family consisting of all integral or rational 
surgeries on a fixed link. 
We will consider the family of all rational surgeries on the Borromean rings, 
which includes other families such as rational surgeries on all twist knots 
and Whitehead links. 
Before proceeding with this, we observe that the set of 
$n$-tuples of rational coefficients 
for which surgery on a fixed link is a Stein boundary will be open in 
$(\que \cup \{-\infty\})^n$, provided that we use the lower limit topology 
on $\que\cup \{-\infty\}$. 
Note that the Stein surfaces involved may have no obvious relation to 
the original link. 

\proclaim{Proposition 5.8}  
Let $L=K_1\cup \cdots \cup K_n$ be a link in $S^3$. 
Suppose that the oriented manifold $M$ obtained by surgery on $L$ with 
coefficients $r_1,\ldots,r_n\in\que \cup \{-\infty\}$ bounds a Stein 
surface $X$. 
Then there are numbers $R_i>r_i$ in $\que$ $(i=1,\ldots,n)$ such that 
whenever $r_i\le r'_i <R_i$, the manifold $M'$ obtained by 
$(r'_1,\ldots,r'_n)$-surgery on $L$ also bounds a Stein surface. 
\endproclaim 

\demo{Proof} 
It suffices to assume that each coefficient $r_i$ is in 
$\zed\cup \{-\infty\}$, for in the general case $M$ can be expanded as an 
integer surgery as in the proof of Proposition~5.3, and then a small 
increase  in $r_i$ can be realized by a small increase in the last 
coefficient $a_k$   of its continued fraction expansion. 
Write $L=L_0\coprod L_\infty$, where $K_i\subset L_\infty$ if and only 
if $r_i = -\infty$. 
Let $L'\subset S^3-L_0\subset M$ denote $L_\infty$ union a meridian $\mu_i$ 
of each $K_i\subset L_0$. 
After a small perturbation, we may assume that $L'$ is Legendrian in 
$(M,\xi)=\partial X$. 
Choose coefficients $s_1,\ldots, s_n\in \que \cup\{-\infty\}$ for $L'$, 
with each $s_i<tb(K'_i)$, where $K'_i=K_i$ or $\mu_i$ in $L'$ and 
the 0-framing on $K'_i$ is determined via the inclusion $S^3-L_0\subset S^3$. 
By the method of proof of Proposition~5.3, 
we can add handles to $X$ to obtain a 
Stein surface whose boundary $M'$ is obtained from $M$ by 
$(s_1,\ldots,s_n)$-surgery on $L'$. 
By slam-dunking each $\mu_i$ in the original picture of $L\subset S^3$, 
we see that $M'$ is also obtained from $S^3$ by surgery on $L$, 
with coefficients $s_i$ on $L_\infty$ and $r_i-\frac1{s_i}$ on $L_0$. 
Thus, we can take $R_i=tb(K_i)$ if $r_i=-\infty$ and $R_i = 
r_i + {1\over\max (1,-tb(\mu_i))}$ otherwise, where $tb$ is measured relative 
to the usual framing convention 
on $S^3$ and the contact structure $\xi$ induced by $X$.\qed
\enddemo 

Now let $M=M(r_1,r_2,r_3)$ denote the oriented 3-manifold obtained by surgery 
on the Borromean rings with coefficients $r_1,r_2,r_3\in\que \cup \{\infty\}$ 
(Figure~46). 
Note that $M$ is independent of the order of the coefficients, and that 
reversing the orientation of $M$ is the same as reversing the sign of 
each $r_i$. 
Since the Borromean rings form a hyperbolic link, $M$ will be hyperbolic 
for ``most'' values of $r_1,r_2,r_3$. 
\cfig{2.00}{rg-fig46.eps}{46}

\proclaim{Theorem 5.9} 
The oriented manifold $M(r_1,r_2,r_3)$ bounds a Stein surface, except 
possibly when $(r_1,r_2,r_3)$ lies in one of the subsets $A_0,A_2,A_3$ of 
$\que^3$ defined below. 
Each point in $A_j$ has exactly $j$ negative coordinates, so $M(r_1,r_2,r_3)$ 
is always a Stein boundary if exactly one $r_i$ is negative. 
\endproclaim 

Here, $A_0=\{(r_1,r_2,r_3)\mid r_1,r_2,r_3\in [1,4)\}$ and $A_2$ is the union 
of $\{(r_1,r_2,r_3)\mid r_1\in [0,\infty)$, $r_2\in [-\frac13,0)$, 
$r_3\in [-2\llbrack -\frac1{r_2}\rrbrack -1,-6)\}$ with its images under 
permutations of the coordinates. 
The set $A_3$ is obtained from $\{(r_1,r_2,r_3)\mid r_1,r_2\in (-\infty,0)$, 
$r_3\in [-2(\llbrack -\frac1{r_1}\rrbrack + \llbrack -\frac1{r_2}\rrbrack+1),
0)\}$ by intersecting it with its images under permutation and deleting 
$\{(r_1,r_2,r_3)\mid r_1,r_2,r_3\in [-6,0)$, two $r_i\in [-1,0)\}$. 

\proclaim{Corollary 5.10} 
{\rm a)} All integer surgeries on the Borromean rings are Stein boundaries 
except possibly for the finite collection whose coefficients $r_1,r_2,r_3$ 
satisfy {\rm (i)}~$r_1,r_2,r_3\in \{1,2,3\}$ 
or {\rm (ii)}~one $r_i$ equals $-1$ and the other two are in $\{-2,-3,-4\}$ 
or {\rm (iii)}~all $r_i=-2$. 

{\rm b)} All homology spheres of the form $M(r_1,r_2,r_3)$ are Stein 
boundaries except possibly if all $r_i=1$, when $M$ is the Poincar\'e 
homology sphere with reversed orientation. 
\endproclaim 

\proclaim{Corollary 5.11} 
For fixed integers 
$\ell,m\in\zed$, let $K(\ell,m)$ be the knot shown in Figure~47, 
generalizing the case of a twist knot, $\ell=\pm1$. 
Then all but finitely many integer surgeries on $K(\ell,m)$ are Stein 
boundaries. 
All rational surgeries on $K(\ell,m)$ are Stein boundaries provided that 
$\ell,m,\ell + m\le 2$ and $\ell,m$ are not both $-1$. 
More specifically, for $r\in \que\cup \{\infty\}$, 
$r$-surgery on $K(\ell,m)$ is a Stein boundary except 
possibly in the following cases: 
{\rm (i)}~$\ell=m=-1$, $1\le r<4$, 
{\rm (ii)}~$\ell<0$, $m\ge 3$, $-2m-1\le r<-6$ (or the same with $\ell$ and 
$m$ interchanged), 
{\rm (iii)}~$\ell,m>0$, $-2(\ell+m+1) \le r<-6$. 
\cdubfig{1.75}{rg-fig47.eps}{47}{1.75}{rg-fig48.eps}{48}
\endproclaim 

\proclaim{Corollary 5.12} 
For $m\in \zed$, let $L(m)$ denote the symmetric link in Figure~48, 
generalizing the Whitehead links $m=\pm1$. 
If $m\le -2$ then all integer surgeries on $L(m)$ are Stein boundaries. 
If $m\le 2$ then all but finitely many such surgeries are Stein boundaries. 
For any $m$ there is a finite subset $A$ of $\zed$ such that  any surgery
on $L(m)$ with coefficients in $\zed-A$ is a Stein boundary. 
Surgeries on $L(m)$ with coefficients $r_1,r_2\in\que\cup \{\infty\}$ are 
Stein boundaries with the following possible exceptions: 
{\rm (i)}~$m=-1$, $r_1,r_2\in [1,4)$, 
{\rm (ii)}~$m<0$, $-\frac13 \le r_i <0$, $-2\llbrack -\frac1{r_i}\rrbrack-1
\le r_j <-6$ for $(i,j) = (1,2)$ or $(2,1)$, 
{\rm (iii)}~$m\ge3$, $r_i\ge0$, $-2m-1\le r_j<-6$ 
for $(i,j)=(1,2)$ or $(2,1)$,  
{\rm (iv)}~$m>0 $, $r_1,r_2<0$, $r_1$ or $r_2<-6$ or both $<-1$, 
$r_1\ge -2(\llbrack -\frac1{r_2}\rrbrack +m+1)$, 
$r_2\ge -2(\llbrack -\frac1{r_1}\rrbrack +m+1)$. 
\endproclaim

\demo{Proof of Corollaries 5.10--5.12} 
Corollary 5.10(a) follows immediately by looking for integer points in 
$A_0\cup A_2 \cup A_3$. 
Part~(b) follows from the observation that
$H_1(M(\frac{p_1}{q_1},\frac{p_2}{q_2},\frac{p_3}{q_3});\zed)\cong 
\zed_{p_1} \oplus \zed_{p_2}\oplus \zed_{p_3}$, so $M(r_1,r_2,r_3)$ is a 
homology sphere precisely when each $\frac1{r_i}\in\zed$. 
Blowing down $M(1,1,1)$ twice exhibits it as $+1$-surgery  on the right 
trefoil, which is well-known to be the Poincar\'e homology sphere oriented as 
the boundary of a {\it positive\/} definite $E_8$-plumbing. 
Similarly, $r$-surgery on the knot $K(\ell,m)$ is the same as 
$M(-\frac1{\ell},-\frac1m,r)$, by a pair of Rolfsen twists. 
Now Corollary~5.11 follows by comparing conditions (i)--(iii) with the sets 
$A_0$, $A_2$ and $A_3$, respectively. 
Surgery on $L(m)$ with coefficients $r_1,r_2$ is the same as 
$M(-\frac1m, r_1,r_2)$, 
and conditions (i)--(iv) of Corollary~5.12 correspond to 
$A_0,A_2,A_2$ and $A_3$, respectively. 
(Note that when $r_1,r_2\in\zed$, (ii) is vacuous and (iv) implies 
$r_1,r_2\ge -2(m+2)$.)\qed
\enddemo 

\cdubfig{1.80}{rg-fig49.eps}{49}{1.80}{rg-fig50.eps}{50}
\demo{Proof of Theorem 5.9} 
First, we exhibit some Stein surfaces bounded by $M(r_1,r_2,r_3)$ for various 
values of the coefficients $r_i$. 
Note that if any $r_i=\infty$ then $M$ is a 
connected sum of lens spaces, hence, a Stein boundary. 
Figure~49 is Stein for $r_1,r_2\ge0$, $r_3<1$, cf.\ Figure~36. 
To see that its boundary is $M(r_1,r_2,r_3)$, surger the 1-handles to  
0-framed 2-handles and slam-dunk the meridians. 
Figure~50 is Stein for $r_1\ge0$, $r_2<0$, $r_3<-2\llbrack-\frac1{r_2}
\rrbrack -1$ (note that the twist is right-handed), 
and its boundary is $M$ by a 
$-\llbrack -\frac1{r_2}\rrbrack$-Rolfsen twist. 
Similarly, Figure~51 solves the case $r_1,r_2<0$, $r_3<-2(\llbrack -\frac1{r_1}
\rrbrack + \llbrack -\frac1{r_2}\rrbrack +1)$. 
Figure~52 is Stein when $r_1\ge1$, $r_2>0$ and $r_3\ge4$. 
(Note that the central $(\llbrack -\frac1{r_2}\rrbrack +1)$-twist makes 
sense, since there are $-2\llbrack -\frac1{r_2}\rrbrack -2\ge0$ left cusps.) 
Finally, Figure~53 is Stein for $r_1\ge-1$, $-1\le r_2<0$, $r_3\ge-6$. 
To verify that Figure~52 depicts $M(r_1,r_2,r_3)$, apply a 
$-\llbrack -\frac1{r_2}\rrbrack$-Rolfsen twist to one component to change its  
coefficient to $r_2$, surger the 1-handle to a 2-handle and change its 
framing to $r_3-4$ by a slam-dunk, then blow down the $-1$-framed unknot. 
Finally, a $+1$-Rolfsen twist on one component 
to restore its coefficient to $r_1$ recovers 
Figure~46. 
(This figure arose as a generalization of the case $r_1=1$, $\frac1{r_2}\in
\zed^+$, the twist knot $K(-1,-\frac1{r_2})$.) 
To recognize Figure~53 as $M$, surger and slam-dunk to change the 1-handle 
to an $(r_3+6)$-framed 2-handle, $+3$-Rolfsen twist to eliminate the 
$-\frac13$-framed unknot, then blow down the $+1$-framed unknot (which only 
links the $(r_3+9)$-framed curve). 
A pair of $-1$-twists to recover coefficients $r_1$ and $r_2$ yields 
Figure~46. 
(Figure~53 was produced by generalizing the case $r_1=r_2=-1$, or 
$r_3$-surgery on the left trefoil, which is $\mp \sum (2,3,6n\pm1)$ if 
$r=\pm \frac1n$.) 
\cfig{2.00}{rg-fig51.eps}{51}
\cfig{3.0}{rg-fig52.eps}{52}
\cfig{2.3}{rg-fig53.eps}{53}

These Stein surfaces are sufficient to prove the theorem. 
Given $r=(r_1,r_2,r_3)\in\que^3$ with at most one $r_i<0$, Figure~49 or 52 
will exhibit $M(r_1,r_2,r_3)$ as a Stein boundary unless $r\in A_0$. 
If exactly two coefficients are negative, we permute so that $r_1\ge 0>
r_2\ge r_3$. 
If $r_2<-1$ then Figure~50 realizes $M$ as a Stein boundary, so we assume 
$r_2\in [-1,0)$. 
Then Figure~50 or 53 works unless $-2\llbrack -\frac1{r_2}\rrbrack -1 
\le r_3 <-6$. 
This last condition is vacuous unless $r_2\ge -\frac13$, in which case 
$r\in A_2$. 
Finally, if all coefficients $r_i$ are negative, then Figure~51 or 53 
will realize $M$ unless $r\in A_3$.\qed
\enddemo

\enddocument